\documentclass[12pt,a4paper,reqno]{amsart}
\pdfoutput=1
\usepackage{amssymb}
\usepackage{fullpage}
\usepackage{ae,aecompl}
\usepackage[english]{babel}
\usepackage{pictex,dcpic}

\title{Anti-selfdual Connections on the\\[10pt] Quantum Projective Plane: Monopoles \\[20pt]}
\date{5 February 2010; v2 18 May 2010}

\author{Francesco D'Andrea} 
\address{D{\'e}partement de Math{\'e}matique, Universit\'e Catholique de Louvain,
Chemin du Cyclotron~2, B-1348, Louvain-La-Neuve, Belgium}
\email{francesco.dandrea@uclouvain.be}

\author{Giovanni Landi \\[20pt]}
\address{Dipartimento di Matematica e Informatica, Universit\`{a} di Trieste,
Via A.~Valerio~12/1, I-34127 Trieste, Italy, and INFN, Sezione di Trieste, Trieste, Italy}
\email{landi@univ.trieste.it}


\allowdisplaybreaks[4]
\numberwithin{equation}{section}

\newtheorem{prop}{Proposition}[section]
\newtheorem{lemma}[prop]{Lemma}
\newtheorem{cor}[prop]{Corollary}
\newtheorem{thm}[prop]{Theorem}

\theoremstyle{definition}
\newtheorem{rem}[prop]{Remark}

\newcommand{\dotimes}{\,\dot{\otimes}\,}
\newcommand{\A}{\mathcal{A}}
\newcommand{\B}{\mathcal{B}}
\newcommand{\E}{\mathcal{E}}
\newcommand{\U}{\mathcal{U}}
\newcommand{\HH}{\mathcal{H}}

\newcommand{\N}{\mathbb{N}}
\newcommand{\Z}{\mathbb{Z}}
\newcommand{\R}{\mathbb{R}}
\newcommand{\C}{\mathbb{C}}
\newcommand{\CP}{\mathbb{C}\mathrm{P}}
\newcommand{\Aq}{\mathcal{A}(\mathbb{C}\mathrm{P}^2_q)}
\newcommand{\Oq}{\mathcal{A}(\mathrm{SU}_q(3))}
\newcommand{\Sq}{\mathcal{A}(S^5_q)}
\newcommand{\Uq}{\mathcal{U}_q(\mathfrak{su}(3))}
\newcommand{\Kq}{\mathcal{U}_q(\mathfrak{u}(2))}

\newcommand{\id}{\textup{id}}

\newcommand{\nn}{\nonumber}
\newcommand{\az}{\triangleright}
\newcommand{\za}{\triangleleft}
\newcommand{\adj}{\stackrel{\mathrm{ad}}{\za}}
\newcommand{\mL}[1]{\mathcal{L}_{#1}}
\newcommand{\ket}[1]{\left|#1\right>}
\newcommand{\inner}[1]{\left<#1\right>}
\newcommand{\sqbn}[1]{\textrm{\footnotesize$\left[\!\!\rule{0pt}{14pt}
                      \smash[b]{\smash[t]{\begin{array}{c}#1\end{array}}}\!\!\right]$}}
\newcommand{\ma}[1]{\left(\!\begin{array}{cccc}#1\end{array}\!\right)}
\newcommand{\maa}[2]{\textrm{\footnotesize$\left(\!\rule{0pt}{#1}\smash[b]{\smash[t]{\begin{array}{ccc}#2\end{array}}}\!\right)$}}
\newcommand{\oh}{\smash[t]{\smash[b]{\tfrac{1}{2}}}}
\newcommand{\aaz}{\,\textrm{\footnotesize$\blacktriangleright$}\,}
\newcommand{\tr}{\mathrm{Tr}}
\newcommand{\nint}{\int\mkern-19mu-\;}
\newcommand{\wprod}{\wedge_q\mkern-1mu}
\newcommand{\de}{\partial}
\newcommand{\deb}{\bar{\partial}}
\newcommand{\dd}{\mathrm{d}}

\begin{document}

\begin{abstract}
We present several results on the geometry of the quantum projective plane. 
They include: explicit generators for the K-theory and the K-homology; 
a real calculus with a Hodge star operator; 
anti-selfdual connections on line bundles with explicit 
computation of the corresponding `classical' characteristic classes (via Fredholm modules);
complete diagonalization of gauged Laplacians on these line bundles;
`quantum' characteristic classes 
via equivariant K-theory and $q$-indices. 
\end{abstract}

\maketitle

\tableofcontents

\section{Introduction}
Quantum spaces are revealing rich geometrical structures and are the subject of intense research activities. In this paper we present a construction of monopoles on the quantum projective plane $\CP^2_q$. By a monopole we mean a line bundle over $\CP^2_q$, that is to say a `rank 1' (in a sense to be made precise) finitely generated projective module over the coordinate algebra $\Aq$, endowed with a connection having anti-selfdual curvature; these are described in Sect.~\ref{se:mono}. Necessary for the anti-selfduality was a differential calculus and a Hodge star operators on forms. In Sect.~\ref{se:diff} we give a full differential $*$-calculus on $\CP^2_q$ with an Hodge star operator which satisfies all the required properties.  

Both the K-theory and the K-homology groups of the quantum projective plane are known: $K_0(\Aq)\simeq\Z^3\simeq K^0(\Aq)$ while $K_1(\Aq) =0= K^1(\Aq)$. Thus, a finitely generated projective module over $\Aq$ is uniquely identified by three integers.

Hermitian vector bundles over a homogeneous space -- actually the corresponding modules of sections --  can be equivalently described as vector valued functions over the total space which are equivariant for a suitable action of the structure group: that is to say, one thinks of sections as 'equivariant maps'. This construction still makes sense for quantum homogeneous spaces. In Sect.~\ref{se:hb} we describe general Hermitian `vector bundles' on $\CP^2_q$ as modules of equivariant elements; we then specialize to line bundles, for which we construct explicit projections (they are all finitely generated projective modules). The generators of the K-homology are obtained in Sect.~\ref{se:top} where we also compute their pairings with the line bundle projections; the resulting numbers are interpreted as the rank (always 1 for line bundles) and 1st and 2nd Chern numbers of the underlying `vector bundles'. In particular, we identify three natural generators for the K-theory of $\CP^2_q$ and three dual generators for its K-homology.

Classically, topological invariants are computed by integrating powers of the curvature of a connection, and the result depends only on the class of the vector bundle, and not on the chosen connection. In order to integrate the curvature of a connection on the quantum projective space $\CP^2_q$ one needs `twisted integrals'; moreover, one does not gets integers any longer but rather $q$-analogues of integers. In Sect.~\ref{se:qinv} we give some general result on equivariant K-theory and K-homology and corresponding Chern-Connes characters, and then focus on $\CP^2_q$ with equivariance under the action of the symmetry algebra $\Uq$. In particular, we construct twisted cocycles that, when paired with the line bundle projections, give $q$-analogues of monopole and instanton numbers of the line bundle. As a corollary, we obtain that when the deformation parameter $q$ is transcendental the equivariant $K_0$-group has (at least) a countable number of generators, $K_0^{\Uq}(\Aq)\supset\Z^\infty$.

After this introduction, we start in Sect.~\ref{se:qqp} with some basic results on the geometry of the quantum projective plane $\CP^2_q$ and its symmetry algebra $\Uq$.

One of our motivation for the present work was to study quantum field theories on noncommutative spaces. The construction of monopoles on $\CP^2_q$ that we present here (instanton configurations on $\CP^2_q$ will be reported in \cite{DL09b}) is a step in this direction. In this spirit, in Sect.~\ref{se:gl} we study gauged Laplacian operators acting on modules of sections of monopole bundles. Their complete diagonalization is made possible by the presence of many symmetries. From the point of view of physics, such an operator would describe `excitations moving on the quantum projective space' in the field of a magnetic monopole. In the limit $q \to 1$ it provides a model of quantum Hall effect on the projective plane.

\medskip

\noindent{\bf Notations.}\\
Throughout this paper the real deformation parameter will be taken to be $0<q<1$.
The symbol $[x]$ denotes the $q$-analogue of $x\in\C$,
$$
[x]:=\frac{q^x-q^{-x}}{q-q^{-1}} \; ;
$$
For $n$ a positive integer the $q$-factorial is $[n]!:=[n][n-1]\ldots[2][1]$, with $[0]!:=1$ and  
the $q$-binomial is given by
$$
\sqbn{n \\ m}:=\frac{[n]!}{[m]![n-m]!}\;. 
$$
For convenience we define the $q$-trinomial coefficient as
\begin{equation}\label{eq:q3coef}
[j,k,l]!=q^{-(jk+kl+lj)}\frac{[j+k+l]!}{[j]![k]![l]!} \;.
\end{equation}
We shall use Sweedler notation for coproducts $\Delta(h) = h_{(1)}\otimes h_{(2)}$ with a sum understood.
Finally, all algebras will be unital $*$-algebras over $\C$, and their representations will be implicitly assumed to be unital $*$-representations.

\section{The quantum projective plane and its symmetries}\label{se:qqp}
We recall from \cite{DDL08b} some results on the geometry of the quantum projective plane.

\subsection{The algebra of `infinitesimal' symmetries}\label{se:2.1}~\\
Let $\Uq$ denote the Hopf $*$-algebra generated by elements $\{K_i,K_i^{-1},E_i,F_i\}_{i=1,2}$,
with $*$-algebra structure $K_i=K_i^*$ and $F_i=E_i^*$, and relations
\begin{gather*}
[K_i,K_j]=0\;,\qquad
[E_i,F_i]=\frac{K_i^2-K_i^{-2}}{q-q^{-1}}\;,  \qquad [E_i,F_j]=0 \quad\mathrm{if}\;i\neq j \;, \\
K_iE_iK_i^{-1}=qE_i\;, \qquad K_iE_jK_i^{-1}=q^{-1/2}E_j \quad\mathrm{if}\;i\neq j \;,  
\end{gather*}
and (Serre relations)
\begin{equation}\label{serre}
E_i^2E_j-(q+q^{-1})E_iE_jE_i+E_jE_i^2=0\qquad\forall\;i\neq j\;.
\end{equation}
Additional relations for the $F_i$'s are obtained from the above ones by taking their $*$-conjugated. 
Our $\Uq$ is the `compact' real form of the Hopf
algebra denoted $\breve{U}_q(\mathfrak{sl}(3))$ in Sect.~6.1.2 of~\cite{KS97}.
With  the $q$-commutator defined as
$$
[a,b]_q:=ab-q^{-1}ba \;,
$$
relations \eqref{serre} can be rewritten in the form $\,[E_i,[E_j,E_i]_q]_q=0\,$
or $\,[[E_i,E_j]_q,E_i]_q=0\,$.
\\
Counit, antipode and coproduct are given by 
\begin{gather*}
\epsilon(K_i)=1\;,\qquad
\epsilon(E_i)=\epsilon(F_i)=0\;, \\
S(K_i)=K_i^{-1}\;,\qquad
S(E_i)=-qE_i\;,\qquad
S(F_i)=-q^{-1}F_i\;, \\
\Delta(K_i)=K_i\otimes K_i\;,\quad
\Delta(E_i)=E_i\otimes K_i+K_i^{-1}\otimes E_i\;, \quad 
\Delta(F_i)=F_i\otimes K_i+K_i^{-1}\otimes F_i\;, 
\end{gather*}
for $i=1,2$. An easy check on generators gives for 
the square of the antipode:
\begin{equation}\label{saut}
S^2(h)= (K_1K_2)^4 \,h\, (K_1K_2)^{-4} \,, \qquad \textup{for  all}  \quad h\in\Uq \, . 
\end{equation}

For obvious reasons we denote $\U_q(\mathfrak{su}(2))$ the Hopf $*$-subalgebra
of $\Uq$ generated by the elements $\{K_1,K_1^{-1},E_1,F_1\}$, while 
$\Kq$ denotes the Hopf $*$-subalgebra generated by $\U_q(\mathfrak{su}(2))$ and $L=K_1K_2^2$
and $L^{-1} = (K_1K_2^2)^{-1}$. We shall also use the notation $\U_q(\mathfrak{u}(1))$ for the 
Hopf $*$-subalgebra generated by $L$ and $L^{-1}$.

We need finite-dimensional  highest weight irreducible representations for which the
$K_i$'s are positive operators; these are all the finite-dimensional irreducible representations
that are well defined for $q\to 1$. For $\Uq$ these representation are labelled by a pair 
of non-negative integers $(n_1,n_2)\in\N^2$. Basis vectors of the representation $(n_1,n_2)$ carry
a multi-index $\underline{j}=(j_1,j_2,m)$ satisfying the constraints
\begin{equation}\label{eq:constr}
j_i=0,1,2,\ldots,n_i \;,\quad i=1,2, \qquad \mathrm{and} \qquad \tfrac{1}{2}(j_1+j_2)-|m|\in\N\;.
\end{equation}
The representation $(n_1,n_2)$ is conveniently described as follows.
We let $\ket{\uparrow}$ denote its highest weight vector, i.e.~$E_i\ket{\uparrow}=0$ and
$K_i\ket{\uparrow}=q^{n_i/2}\ket{\uparrow}$ for $i=1,2$, and we let $X_{\underline{i}}^{n_1,n_2}\in\Uq$ be
the element
\begin{equation}\label{exs}
X_{j_1,j_2,m}^{n_1,n_2}:=N_{j_1,j_2,m}^{n_1,n_2}\!\!\sum_{k=0}^{n_1-j_1}\!\!
\frac{q^{-k(j_1+j_2+k+1)}}{[j_1+j_2+k+1]!}\sqbn{\!n_1-j_1\! \\ k}
F_1^{\frac{1}{2}(j_1+j_2)-m+k}[F_2,F_1]_q^{n_1-j_1-k}F_2^{j_2+k} \,,\!
\end{equation}
with
$$
N_{j_1,j_2,m}^{n_1,n_2}:=\sqrt{[j_1+j_2+1]}\sqrt{\frac{[\tfrac{j_1+j_2}{2}+m]!}
{[\tfrac{j_1+j_2}{2}-m]!}\,\frac{[n_2-j_2]![j_1]!}{[n_1-j_1]![j_2]!}\,
\frac{[n_1+j_2+1]![n_2+j_1+1]!}{[n_1]![n_2]![n_1+n_2+1]!}}\;.
$$
Then, the vector space $V_{n_1,n_2}$ underlying the representation is spanned by
the vectors $\ket{\underline{i}}:=X_{\underline{i}}^{n_1,n_2}\ket{\uparrow}$. 
Using the commutation rules of $\Uq$ one proves that in this basis
the representation is given by (cf.~\cite[Sect.~2]{DDL08b})
\begin{align*}
K_1\ket{j_1,j_2,m} &:=q^m\ket{j_1,j_2,m} \;,\\
K_2\ket{j_1,j_2,m} &:=q^{\frac{3}{4}(j_1-j_2)+\frac{1}{2}(n_2-n_1-m)}\ket{j_1,j_2,m} \;,\\
E_1\ket{j_1,j_2,m} &:=\sqrt{[\tfrac{1}{2} (j_1+j_2)-m][\tfrac{1}{2} (j_1+j_2)+m+1]}
   \,\ket{j_1,j_2,m+1} \;,\\
E_2\ket{j_1,j_2,m} &:=\sqrt{[\tfrac{1}{2} (j_1+j_2)-m+1]}\,A_{j_1,j_2}\ket{j_1+1,j_2,m-\oh} \\
     &\qquad\qquad\qquad +\sqrt{[\tfrac{1}{2} (j_1+j_2)+m]}\,B_{j_1,j_2}\ket{j_1,j_2-1,m-\oh}\;,
\end{align*}
with $F_i$ the transpose of $E_i$, and with coefficients given by
\begin{align*}
A_{j_1,j_2}&:=\sqrt{\frac{[n_1-j_1][n_2+j_1+2][j_1+1]}
              {[j_1+j_2+1][j_1+j_2+2]}} \;,\\[5pt]
B_{j_1,j_2}&:=\begin{cases}
\sqrt{\dfrac{[n_1+j_2+1][n_2-j_2+1][j_2]}{[j_1+j_2][j_1+j_2+1]}}
\quad & \mathrm{if}\;j_1+j_2\neq 0\;, \\
1 & \mathrm{if}\;j_1+j_2=0\;.
\end{cases}
\end{align*}
It is a $*$-representation for the inner product
$\inner{\underline{i}|\underline{i}'}=\delta_{\underline{i},\underline{i}'}$.
With our notation, the highest weight vector is $\ket{\uparrow}=\ket{n_1,0,\frac{1}{2}n_1}$.
The representation $\rho^{n_1,n_2}:\Uq\to\mathrm{End}(V_{n_1,n_2})$,
has matrix elements $\rho^{n_1,n_2}_{\underline{i},\underline{j}}(h):=\inner{\underline{i}|h|\underline{j}}$.

Since $L=K_1K_2^2$ commutes with all elements of $\U_q(\mathfrak{su}(2))$,
any irreducible representation $\sigma_{\ell,N}$ of $\Kq$ is the product  of a 
representation of $\U_q(\mathfrak{su}(2))$ (these are labelled by the `spin' $\ell$ with $\ell\in\frac{1}{2}\N$), and a representation of 
charge $N$ of the Hopf $*$-subalgebra $\U_q(\mathfrak{u}(1))$ generated by $L$, that is $\sigma_{\ell,N}(L)=q^N$.
Restricting the representation $\rho^{n_1,n_2}$ of $\Uq$ to $\Kq$ we get: 
$$
\rho^{n_1,n_2}\big|_{\Kq}\simeq\bigoplus_{2\ell=0}^{n_1+n_2}\;
\bigoplus_{\frac{N+n_1-n_2}{3}=\max\{-\ell,\ell-n_2\}}^{\min\{\ell,n_1-\ell\}}
\sigma_{\ell,N} \;,
$$
and thus the representations $\sigma_{\ell,N}$ of $\Kq$ appearing as components in at least one
representation of $\Uq$ are only those for which $\ell+\frac{1}{3}N \in \frac{1}{3}\Z$, or equivalently 
$\ell+N\in\Z$ (since $2\ell\in\Z$).

For later use (in Sect.~\ref{se:gl}) we need the Casimir operator; it is the operator given 
(in a slightly enlarged algebra, cf. \cite{DDL08b}) by
\begin{multline}\label{eq:Cq}
\mathcal{C}_q = (q-q^{-1})^{-2} \Big\{ (H+H^{-1}) \bigl((qK_1K_2)^2+(qK_1K_2)^{-2}\bigr)
 +H^2+H^{-2}-6 \Big\} \\
+\left(qHK_2^2+q^{-1}H^{-1}K_2^{-2}\right)F_1E_1+\left(qH^{-1}K_1^2+q^{-1}HK_1^{-2}\right)F_2E_2 \\
+qH[F_2,F_1]_q[E_1,E_2]_q+qH^{-1}[F_1,F_2]_q[E_2,E_1]_q \; ,
\end{multline}
with $H:=(K_1K_2^{-1})^{2/3}$. In the representation $\rho^{n_1,n_2}$  it has spectrum
\begin{equation}\label{eq:SpCq}
\mathcal{C}_q\bigr|_{V_{n_1,n_2}}=[\tfrac{1}{3}(n_1-n_2)]^2+
[\tfrac{1}{3}(2n_1+n_2)+1]^2+[\tfrac{1}{3}(n_1+2n_2)+1]^2\;.
\end{equation}

\subsection{The quantum $\mathrm{SU}(3)$ group}\label{se:2.2}~\\
The deformation $\Oq$ of the Hopf $*$-algebra of
representative functions of $\mathrm{SU}(3)$ is given in~\cite{RTF90} 
(see also \cite{KS97}, Sect.~9.2).
As a $*$-algebra it is generated by elements $u^i_j$, with $i,j=1,2,3$ ,
having commutation relations
\begin{align*}
u^i_ku^j_k &=qu^j_ku^i_k \;,&
u^k_iu^k_j &=qu^k_ju^k_i \;,&&
\forall\;i<j\;, \\
[u^i_l,u^j_k]&=0 \;,&
[u^i_k,u^j_l]&=(q-q^{-1})u^i_lu^j_k \;,&&
\forall\;i<j,\;k<l\;. 
\end{align*}
There is also a cubic relation
$$
\sum\nolimits_{\pi\in S_3}(-q)^{||\pi||}u^1_{\pi(1)}u^2_{\pi(2)}u^3_{\pi(3)}=1 \;,
$$
with the sum over all permutations $\pi$ of the three elements
$\{1,2,3\}$ and $||\pi ||$ is the length of $\pi$.
The $*$-structure is given by
$$
(u^i_j)^*=(-q)^{j-i}(u^{k_1}_{l_1}u^{k_2}_{l_2}-qu^{k_1}_{l_2}u^{k_2}_{l_1}) \;,
$$
with $\{k_1,k_2\}=\{1,2,3\}\smallsetminus\{i\}$ and
$\{l_1,l_2\}=\{1,2,3\}\smallsetminus\{j\}$, as ordered sets.
Thus for example $(u^1_1)^*=u^2_2u^3_3-qu^2_3u^3_2$.
Coproduct, counit and antipode are the standard ones:
$$
\Delta(u^i_j)=\sum\nolimits_ku^i_k\otimes u^k_j\;,\qquad
\epsilon(u^i_j)=\delta^i_j\;,\qquad
S(u^i_j)=(u^j_i)^*\;.
$$
There is a non-degenerate dual pairing (cf.~\cite{KS97}, Sect.~9.4)
$$
\inner{\,,\,}:\Uq\times\Oq\to\C \, ,
$$ 
which allows one to define left $\az$
and right $\za$ canonical actions of $\Uq$ on $\Oq$,
$$
h\az a=a_{(1)}\inner{h,a_{(2)}} \;, \qquad \textup{and} \qquad  a\za h=\inner{h,a_{(1)}}a_{(2)}.
$$
By using the counit on these equations one gets that
\begin{equation}\label{counitact}
\inner{h,a}=\epsilon(h\az a)
=\epsilon(a\za h) , 
\end{equation}
for all $h \in \Uq$ and $a \in \Oq$. Also, it is known that the left (resp. right) canonical action is dual to the right (resp. left) regular action. For the case at hand this is the statement that 
\begin{equation}\label{le-ri}
\inner{x \,h \,y, a}=\inner{x,a_{(1)}}\inner{h,a_{(2)}}\inner{y,a_{(3)}}
=\inner{h,y\az a\za x} ,
\end{equation}
for all $x,y, h \in \Uq$ and $a \in \Oq$. On generators the actions are given by
\begin{align}\label{expact}
K_i\az u^j_k &=q^{\frac{1}{2}(\delta_{i+1,k}-\delta_{i,k})}u^j_k \;,&
E_i\az u^j_k &=\delta_{i,k} u^j_{i+1}\;, &
F_i\az u^j_k &=\delta_{i+1,k} u^j_i\;, \nn \\
u^j_k\za K_i &=q^{\frac{1}{2}(\delta_{i+1,j}-\delta_{i,j})}u^j_k \;,&
u^j_k\za E_i &=\delta_{i+1,j} u^i_k \;, &
u^j_k\za F_i &=\delta_{i,j} u^{i+1}_k \;.
\end{align}
By the Peter-Weyl theorem~\cite{KS97}
a linear basis $\{t^{n_1,n_2}_{\underline{i},\underline{j}}\}$
of $\Oq$ is defined  implicitly by
$$
\inner{h,\smash[t]{\smash[b]{t^{n_1,n_2}_{\underline{i},\underline{j}}}}}=
\rho^{n_1,n_2}_{\underline{i},\underline{j}}(h) \;,\qquad\forall\;h\in\Uq\;.
$$
{}From the definition it follows that
$\Delta(t^{n_1,n_2}_{\underline{i},\underline{j}})=
\sum_{\underline{k}}t^{n_1,n_2}_{\underline{i},\underline{k}}
\otimes t^{n_1,n_2}_{\underline{k},\underline{j}}$
and $(t^{n_1,n_2}_{\underline{i},\underline{j}})^*=
S(t^{n_1,n_2}_{\underline{j},\underline{i}})$ (that is,
$t^{n_1,n_2}$ is a unitary matrix). Also, this matrix transforms
according to the representation $\rho^{n_1,n_2}$ under
the left/right canonical action:
\begin{equation}\label{eq:azt}
h\az t^{n_1,n_2}_{\underline{i},\underline{j}}=
\sum\nolimits_{\underline{k}}
t^{n_1,n_2}_{\underline{i},\underline{k}}
\rho^{n_1,n_2}_{\underline{k},\underline{j}}(h)
\;,\qquad
t^{n_1,n_2}_{\underline{i},\underline{j}}\za h=
\sum\nolimits_{\underline{k}}
\rho^{n_1,n_2}_{\underline{i},\underline{k}}(h)
t^{n_1,n_2}_{\underline{k},\underline{j}} \;,
\end{equation}
for all $h\in\Uq$.
For $(n_1,n_2)=(0,1)$ the elements $t^{n_1,n_2}_{\underline{i},\underline{j}}$
are just the generators $u^i_j$ (properly reordered). Let us describe them in general.
\begin{prop}
We have
\begin{equation}\label{eq:eqdef}
t^{n_1,n_2}_{\underline{i},\underline{j}}=
X_{\underline{j}}^{n_1,n_2}\az\{(u_1^1)^*\}^{n_1}(u_3^3)^{n_2}
\za(X_{\underline{i}}^{n_1,n_2})^* \;,
\end{equation}
where $X_{\underline{i}}^{n_1,n_2}$ are given in \eqref{exs}.
\end{prop}

\begin{proof}
By the Poincar{\'e}-Birkhoff-Witt theorem \cite[Thm.~$6.24'$]{KS97},
the vector space $\Uq$ is spanned by elements of the form
$h=f\,K_1^rK_2^s\,e$, where $e$ is a product of $E_i$'s and $f$
is a product of $F_i$'s. Since $E_i\ket{\uparrow}=\left<\uparrow\right|F_i=0$,
it follows that $\rho^{n_1,n_2}_{\uparrow,\uparrow}(h)=0$
unless $e$ and $f$ have degree zero. For $e=f=1$ one gets 
$\rho^{n_1,n_2}_{\uparrow,\uparrow}(h)=q^{\frac{1}{2}(rn_1+sn_2)}$, since
$K_i\ket{\uparrow}=q^{n_i/2}\ket{\uparrow}$.

Let us define $u^{n_1,n_2}_{\uparrow,\uparrow}:=\{(u_1^1)^*\}^{n_1}(u_3^3)^{n_2}$.
{}From the explicit formul{\ae} \eqref{expact} it follows that
$E_i\az u^{n_1,n_2}_{\uparrow,\uparrow}=0=u^{n_1,n_2}_{\uparrow,\uparrow}\za F_i$. Thus, by using \eqref{counitact},
$\inner{fK_1^rK_2^se,u^{n_1,n_2}_{\uparrow,\uparrow}}=0$ unless $e=f=1$.
When $e=f=1$, from $K_i\az u^{n_1,n_2}_{\uparrow,\uparrow}=u^{n_1,n_2}_{\uparrow,\uparrow}\za K_i
=q^{n_i/2}u^{n_1,n_2}_{\uparrow,\uparrow}$ it derives that 
$\inner{K_1^rK_2^s,t^{n_1,n_2}_{\uparrow,\uparrow}}=q^{\frac{1}{2}(rn_1+sn_2)}$.
Hence, $\inner{h,u^{n_1,n_2}_{\uparrow,\uparrow}}=
\rho^{n_1,n_2}_{\uparrow,\uparrow}(h)$ for all $h\in\Uq$.
But this implies that
\begin{align*}
\inner{h,\smash[b]{t^{n_1,n_2}_{\underline{i},\underline{j}}}} &=
\rho^{n_1,n_2}_{\underline{i},\underline{j}}(h)
=\inner{\underline{i}|\smash[b]{h|\underline{j}}} 
=\inner{\,\uparrow\!|\smash[b]{(X_{\underline{i}}^{n_1,n_2})^*h
X_{\underline{j}}^{n_1,n_2}|\!\uparrow\,}} \\
& =\rho^{n_1,n_2}_{\uparrow,\uparrow}((X_{\underline{i}}^{n_1,n_2})^*h
X_{\underline{j}}^{n_1,n_2}) 
=\inner{(\smash[b]{X_{\underline{i}}^{n_1,n_2})^*h
X_{\underline{j}}^{n_1,n_2},u^{n_1,n_2}_{\uparrow,\uparrow}}} \\
&=\inner{h\smash[b]{,X_{\underline{j}}^{n_1,n_2}\az u^{n_1,n_2}_{\uparrow,\uparrow}\za (X_{\underline{i}}^{n_1,n_2})^*}} \;,
\end{align*}
using the identity \eqref{le-ri}. Thus $t^{n_1,n_2}_{\underline{i},\underline{j}}=X_{\underline{j}}^{n_1,n_2}\az
u^{n_1,n_2}_{\uparrow,\uparrow}\za (X_{\underline{i}}^{n_1,n_2})^*$,
which is just (\ref{eq:eqdef}).
\end{proof}

\subsection{The quantum 5-sphere and the quantum projective plane}\label{sec:2.3}~\\
The most natural way to arrive to $\CP^2_q$ is via the $5$-sphere $S^5_q$. We shall therefore start from the algebra
of coordinate functions on the latter, defined as
$$
\Sq:=\big\{a\in\Oq\,\big|\,a\za h=\epsilon(h)a\;, \;\;\forall\;h\in\U_q(\mathfrak{su}(2))\big\} \;
$$
and, as such, it is the $*$-subalgebra of $\Oq$ generated by elements $\{u_i^3,\,i=1,\ldots,3\}$ of the last
`row'. In~\cite{VS91} it is proved to be isomorphic, through the identification $z_i=u_i^3$, to the abstract
$*$-algebra with generators $z_i,z_i^*$, $i=1,2,3$, and relations:
\begin{gather*}
z_iz_j=qz_jz_i\quad\forall\;i<j \;,\qquad\quad
z_i^*z_j=qz_jz_i^*\quad\forall\;i\neq j \;, \\
\rule{0pt}{16pt}
[z_1^*,z_1]=0 \;,\qquad
[z_2^*,z_2]=(1-q^2)z_1z_1^* \;,\qquad
[z_3^*,z_3]=(1-q^2)(z_1z_1^*+z_2z_2^*) \;, \\
\rule{0pt}{16pt}
z_1z_1^*+z_2z_2^*+z_3z_3^*=1 \;.
\end{gather*}
These relations will be useful later on.

The $*$-algebra $\Aq$ of coordinate functions on the quantum projective plane $\CP^2_q$ is the fixed point subalgebra of $\Oq$ for the right action of $\Kq$,
\begin{multline}\label{qcp}
\Aq:=\big\{a\in\Oq\,\big|\,a\za h=\epsilon(h)a\,,\;\forall\;h\in\Kq\big\} \\ \cong
 \big\{a\in\Sq\,\big|\,a\za K_1K_2^2=a\big\}\;.
\end{multline}
Clearly, both $\Sq$ and $\Aq$ are left $\Uq$-module $*$-algebras. 

The $*$-algebra $\Aq$ is generated by elements
$p_{ij}:=z_i^*z_j=(u_i^3)^*u_j^3 = p_{ji}^*$ and 
from the relations of $\Sq$ one gets analogous commutation relations for $\Aq$:
\begin{align*}
p_{ij}p_{kl}&=q^{\mathrm{sign}(i-k)+\mathrm{sign}(l-j)}\, p_{kl}p_{ij}
 &\hspace{-1.5cm}\mathrm{if}\;i\neq l\;\mathrm{and}\;j\neq k\;,\\
p_{ij}p_{jk}&=q^{\mathrm{sign}(i-j)+\mathrm{sign}(k-j)+1}\, p_{jk}p_{ij}-(1-q^2)
\textstyle{\sum_{l<j}}\, p_{il}p_{lk} &\mathrm{if}\;i\neq k\;,\\
p_{ij}p_{ji}&=
q^{2\mathrm{sign}(i-j)}p_{ji}p_{ij}+
(1-q^2)\left(\textstyle{\sum_{l<i}}\,q^{2\mathrm{sign}(i-j)}p_{jl}p_{lj}
-\textstyle{\sum_{l<j}}\,p_{il}p_{li}\right)
 &\mathrm{if}\;i\neq j\;,
\end{align*}
with $\mathrm{sign}(0):=0$.
The elements $p_{ij}$ are  the matrix entries of a projection $P= (p_{ij})$, that is $P^2=P=P^*$. 
This projection has  $q$-trace:
$$
\tr_q(P):=q^4p_{11}+q^2p_{22}+p_{33}=1.
$$
This projection gives the `tautological line bundle' on $\CP^2_q$; general line bundles will be discussed in Sect.~\ref{se:lb} below.

\begin{rem}
The two equalities in \eqref{qcp} give algebra inclusions  
$\Aq \hookrightarrow \Oq$ and $\Aq \hookrightarrow \Sq$.
These could be seen as `noncommutative principal bundles' with `structure Hopf algebra' $\mathrm{U}_q(2)$ and $U(1)$ respectively. 
Indeed, the action of $\Kq$ on $\Oq$ dualizes to a coaction of 
$\mathrm{U}_q(2)$ for which $\Aq$ is the algebra of coinvariants. Analogously, the action of $\U_q(\mathfrak{u}(1))$
on $\Sq$ dualizes to a coaction of $U_q(1)\simeq U(1)$ 
for which again $\Aq$ is the algebra of coinvariants. These are noncommutative analogues of the classical $\mathrm{U}(2)$-principal bundle $\mathrm{SU}(3) \to \CP^2$ 
and $U(1)$-principal bundle $S^5 \to \CP^2$.
\end{rem}

\section{Hermitian vector bundles}\label{se:hb}

On the manifold $\CP^2_q$ we shall select suitable `monopole'  bundles. 
These will come as associated bundles to the principal fibrations on $\CP^2_q$ mentioned in the previous section. We start with the general construction of associated bundles: the starting idea is to define their modules of `sections' as equivariant vector valued functions on $\mathrm{SU}_q(3)$.

\subsection{Hermitian bundles of any rank}~\\
With any $n$-dimensional $*$-representation $\sigma:\Kq\to\mathrm{End}(\C^n)$, one associates an
$\Aq$-bimodule of equivariant elements:
\begin{multline*}
\mathfrak{M}(\sigma) = \Oq\!\boxtimes_{\sigma}\!\C^n \\ :=
\big\{v\in\Oq^n ~\big|~ \sigma(S(h_{(1)})) (v\za h_{(2)}) = \epsilon(h)v\,;\;\; \forall\;h\in\Kq\big\} \;,
\end{multline*}
where $v=(v_1,\ldots,v_n)^t$ is a column vector and row by column multiplication
is implied. It is easy to see that $\mathfrak{M}(\sigma)$ is also a left $\Aq\rtimes\Uq$-module.
In particular, $\Aq=\mathfrak{M}(\epsilon)$ is the module associated to the trivial
representation given by the counit.
A natural $\Aq$-valued Hermitian structure on $\mathfrak{M}(\sigma)$ is defined by
\begin{equation}\label{}
\mathfrak{M}(\sigma)\times\mathfrak{M}(\sigma)\to\Aq \;,\qquad
(v,v')\mapsto v^\dag \cdot v' \;,
\end{equation}
where $v^\dag$ is the conjugate transpose of $v$ and again one multiplies row by column. Indeed, if $v,v'\in\mathfrak{M}(\sigma)$, for $h\in\Kq$, and denoting $t:=S(h)^*$, we have that
\begin{align*}
v^\dag \cdot v' \za h
 &=(v^\dag\za h_{(1)})(v'\za h_{(2)}) 
 =(v^\dag\za h_{(1)})\sigma(h_{(2)})\sigma(S(h_{(3)}))(v'\za h_{(4)}) \\ 
 &=(v^\dag\za h_{(1)})\sigma(h_{(2)})\epsilon(h_{(3)})v' 
 =\bigl\{\sigma(h_{(2)}^*)v\za S(h_{(1)})^*\bigr\}^*v' \\
 &=\bigl\{\sigma(S(t_{(1)}))v\za t_{(2)})\bigr\}^*v' 
 =\epsilon(t)^* v^\dag \cdot v' =\epsilon(h) \, v^\dag \cdot v'  \;;
\end{align*}
thus $v^\dag \cdot v' \in\Aq$ as it should be.
We think of each $\mathfrak{M}(\sigma)$ endowed with this Hermitian structure as
the module of sections of an equivariant `Hermitian vector bundle' over $\CP^2_q$. 

\begin{rem}\label{innprod}
Composing the Hermitian structure with (the restriction to $\Aq$ of) the Haar functional $\varphi:\Oq\to\C$, we get a non-degenerate $\C$-valued  
inner product on $\mathfrak{M}(\sigma)$, $\inner{v,v'}:=\varphi( v^\dag \cdot v' )$. This will be used later on in Sect.~\ref{se:hsci} when defining a Hodge operator on forms and gauged Laplacian operators on modules of sections.
\end{rem}

Clearly $\mathfrak{M}(\sigma_1)\oplus\mathfrak{M}(\sigma_2)\simeq\mathfrak{M}(\sigma_1\oplus\sigma_2)$
for any pair of representations $\sigma_1,\sigma_2$. We have also an inclusion
$\mathfrak{M}(\sigma_1)\otimes_{\Aq}\mathfrak{M}(\sigma_2)\subset\mathfrak{M}(\sigma)$,
where $\sigma:=\sigma_1\otimes\sigma_2$ is the Hopf tensor product of the two representations. Indeed,
for $v=(v_1,\ldots,v_n)^t\in\mathfrak{M}(\sigma_1)$ and $v'=(v'_1,\ldots,v'_{n'})^t\in\mathfrak{M}(\sigma_2)$, the product $w=v\otimes_{\Aq}v'$ satisfies
\begin{align*}
\sigma(S(h_{(1)}))\cdot(w\za h_{(2)})
&=\sigma_1(S(h_{(1)})_{(1)})
(v\za h_{(2)(1)})\otimes_{\Aq}\sigma_2(S(h_{(1)})_{(2)})(v'\za h_{(2)(2)}) \\
&=\sigma_1(S(h_{(2)}))(v\za h_{(3)})\otimes_{\Aq}\sigma_2(S(h_{(1)}))(v'\za h_{(4)}) \\
&=\epsilon(h_{(2)})v\otimes_{\Aq}\sigma_2(S(h_{(1)}))(v'\za h_{(3)}) \\
&=v\otimes_{\Aq}\sigma_2(S(h_{(1)}))(v'\za h_{(2)}) \\
&=\epsilon(h)v\otimes_{\Aq}v'=\epsilon(h)w\;,
\end{align*}
where we used anti-comultiplicativity of the antipode,
$$
S(h)_{(1)}\otimes S(h)_{(2)}=S(h_{(2)}\otimes h_{(1)})\;.
$$
Thus, $w\in\mathfrak{M}(\sigma)$.
The inclusion is an isometry for the natural inner product
$$
\bigl<v_1\otimes_{\Aq}v_2,v'_1\otimes_{\Aq}v_2'\bigr>=\inner{v_1,v_1'}\inner{v_2,v_2'}
$$
on $\mathfrak{M}(\sigma_1)\otimes_{\Aq}\mathfrak{M}(\sigma_2)$.

For the representations discussed at the end of Sect.~\ref{se:2.1}, 
$\sigma_{\ell,N}:\Kq\to\mathrm{End}(\C^{2\ell+1})$, with 
$\ell\in\frac{1}{2}\N$ and $\ell+N\in\Z$,  we denote
$\Sigma_{\ell,N}:=\mathfrak{M}(\sigma_{\ell,N})$.
All these $\Sigma_{\ell,N}$ are finitely-generated and projective as one sided (i.e.~both as left or right)
$\Aq$-modules. Now, since there always exists an irreducible representation $\rho^{n_1,n_2}$
of $\Uq$ containing $\sigma_{\ell,N}$ as a summand, and since $\mathfrak{M}(\sigma_{\ell,N})$
is a direct summand (as a bimodule) in the corresponding  
bimodule $\mathfrak{M}(\rho^{n_1,n_2}|_{\Kq})$,
to establish the statement it is enough to show that all the $\mathfrak{M}(\rho^{n_1,n_2}|_{\Kq})$ are
free both as left and right (but not as bimodules!) $\Aq$-modules. But this is easy; indeed,
as left $\Aq$-modules, we have the following isomorphism
\begin{align*}
\phi &:\Aq^{\dim\rho^{n_1,n_2}}\to\mathfrak{M}(\rho^{n_1,n_2}|_{\Kq})\;,
& \phi(a)_{\underline{i}}&:={\textstyle\sum_{\underline{j}}}a_{\underline{j}}t^{n_1,n_2}_{\underline{i},\underline{j}} \;,\\
\phi^{-1} &:\mathfrak{M}(\rho^{n_1,n_2}|_{\Kq})\to\Aq^{\dim\rho^{n_1,n_2}}\;,
& \phi^{-1}(v)_{\underline{i}}&:={\textstyle\sum_{\underline{j}}}v_{\underline{j}}(t^{n_1,n_2}_{\underline{j},\underline{i}})^* \;,
\end{align*}
where $t^{n_1,n_2}_{\underline{i},\underline{j}}\in\Oq$ are the elements defined in \eqref{eq:eqdef}. 
The maps $\phi$ and $\phi^{-1}$ are well defined from (\ref{eq:azt}); and clearly, they are left $\Aq$-module maps. That they are
one the inverse of the other follows from the unitarity of the matrix $t^{n_1,n_2}$.

Similarly,  as right $\Aq$-modules we have the following isomorphism
\begin{align*}
\phi' &:\Aq^{\dim\rho^{n_1,n_2}}\to\mathfrak{M}(\rho^{n_1,n_2}|_{\Kq})\;,
& \phi'(a)_{\underline{i}}&:={\textstyle\sum_{\underline{j}}}t^{n_1,n_2}_{\underline{i},\underline{j}}a_{\underline{j}} \;,\\
\phi'^{-1} &:\mathfrak{M}(\rho^{n_1,n_2}|_{\Kq})\to\Aq^{\dim\rho^{n_1,n_2}}\;,
& \phi'^{-1}(v)_{\underline{i}}&:=
{\textstyle\sum_{\underline{j}}}(t^{n_1,n_2}_{\underline{j},\underline{i}} )^*v_{\underline{j}} \;.
\end{align*}

\subsection{Line bundles}\label{se:lb}~\\
In the construction of the vector bundles over $\CP^2_q$ we have in particular $\Sigma_{0,0}=\Aq$. Moreover $\Sigma_{0,N}$ are `sections' of `line bundles'; before we discuss them in more details in this  section, we note that they can equivalently be obtained out of the sphere $S^5_q$:
\begin{equation}
\Sigma_{0,N} \simeq
 \big\{\eta \in\Sq\, \; \big|\; \eta \za K_1K_2^2 = q^{N} \eta \big\}\;.
\end{equation}
We are ready to compute projections $P_N$ for the modules $\Sigma_{0,N}$, thus realizing them as 
$\Sigma_{0,N}\simeq P_N \Aq^{d_N}$ as right modules, 
and $\Sigma_{0,N}\simeq \Aq^{d_N}P_{-N}$ as left modules
(notice the `$-$' sign; and remember that $N$ labels representations), for a suitable integer $d_N$. 

\begin{lemma}
For $N\geq 0$ we have
\begin{gather*}
\sum\nolimits_{j+k+l=N}[j,k,l]!(z_1^jz_2^kz_3^l)(z_1^jz_2^kz_3^l)^*=1 \;,\\
\sum\nolimits_{j+k+l=N}q^{2(j-l)}[j,k,l]!(z_1^jz_2^kz_3^l)^*(z_1^jz_2^kz_3^l)=q^{-2N} \;,
\end{gather*}
where $[j,k,l]!$ are the $q$-trinomial coefficients in (\ref{eq:q3coef}).
\end{lemma}

\begin{proof}
We seek coefficients $c_{j,k,l}(N)$ such that
$\sum\nolimits_{j+k+l=N}c_{j,k,l}(N)(z_1^jz_2^kz_3^l)(z_1^jz_2^kz_3^l)^*=1$.
{}From the algebraic identity
\begin{align*}
1 & = \sum\nolimits_{j+k+l=N+1}c_{j,k,l}(N+1)(z_1^jz_2^kz_3^l)(z_1^jz_2^kz_3^l)^*  \\
& = \sum\nolimits_{j+k+l=N}c_{j,k,l}(N)(z_1^jz_2^kz_3^l)(z_1z_1^*+z_2z_2^*+z_3z_3^*)(z_1^jz_2^kz_3^l)^* 
\\
& = \sum\nolimits_{j+k+l=N}q^{-2(k+l)}c_{j,k,l}(N)(z_1^{j+1}z_2^kz_3^l)(z_1^{j+1}z_2^kz_3^l)^* \\
& \qquad\qquad 
+\sum\nolimits_{j+k+l=N}q^{-2l}c_{j,k,l}(N)(z_1^jz_2^{k+1}z_3^l)(z_1^jz_2^{k+1}z_3^l)^* \\
& \qquad\qquad 
+\sum\nolimits_{j+k+l=N}c_{j,k,l}(N)(z_1^jz_2^kz_3^{l+1})(z_1^jz_2^kz_3^{l+1})^* \;,
\end{align*}
we get the recursive equations
\begin{equation}\label{eq:tetra}
c_{j,k,l}(N+1)=q^{-2(k+l)}c_{j-1,k,l}(N)+q^{-2l}c_{j,k-1,l}(N)+c_{j,k,l-1}(N) \;,
\end{equation}
where $N+1=j+k+l$ and with `initial datum' $c_{0,0,0}(0)=1$. Thus, $c_{j,k,l}(N)$ are the vertices of a
$q$-Tartaglia tetrahedron. Since
$$
[j+k+l]=q^{-k-l}[j]+q^{j-l}[k]+q^{j+k}[l] \;,
$$
one verifies that $c_{j,k,l}(N)=[j,k,l]!$ is a solution by plugging it into (\ref{eq:tetra}).

Similarly, since $q^4z_1^*z_1+q^2z_2^*z_2+z_3^*z_3=1$, we have the algebraic identity
\begin{align*}
1 & = \sum\nolimits_{j+k+l=N+1}d_{j,k,l}(N+1)(z_1^jz_2^kz_3^l)^*(z_1^jz_2^kz_3^l)  \\
& =\sum\nolimits_{j+k+l=N}d_{j,k,l}(N)(z_1^jz_2^kz_3^l)^*
  (q^4z_1^*z_1+q^2z_2^*z_2+z_3^*z_3)(z_1^jz_2^kz_3^l) \\
& = \sum\nolimits_{j+k+l=N}q^4d_{j,k,l}(N)(z_1^{j+1}z_2^kz_3^l)^*(z_1^{j+1}z_2^kz_3^l)
\\ & \qquad\qquad 
+\sum\nolimits_{j+k+l=N}q^{-2j+2}d_{j,k,l}(N)(z_1^jz_2^{k+1}z_3^l)^*(z_1^jz_2^{k+1}z_3^l)
\\ & \qquad\qquad 
+\sum\nolimits_{j+k+l=N}q^{-2(j+k)}d_{j,k,l}(N)(z_1^jz_2^kz_3^{l+1})^*(z_1^jz_2^kz_3^{l+1})
\;,
\end{align*}
that gives the recursive equations on the coefficients 
$$
d_{j,k,l}(N+1)=q^4d_{j-1,k,l}(N)+q^{-2j+2}d_{j,k-1,l}(N)+q^{-2(j+k)}d_{j,k,l-1}(N) \;,
$$
with `initial datum' $d_{0,0,0}(0)=1$. The solution is $d_{j,k,l}(N)=q^{2N}q^{2(j-l)}[j,k,l]!$,
as one can check by using the identity $[j+k+l]=q^{k+l}[j]+q^{l-j}[k]+q^{-j-k}[l]$.
\end{proof}

\begin{prop}\label{lemma:linebundles}
Define
\begin{align*}
(\psi_{j,k,l}^N)^*&:=\sqrt{[j,k,l]!}\,z_1^jz_2^kz_3^l\;,
   && \textup{if}\;N\geq 0\; \quad\textup{and with}\quad \; j+k+l=N\,,\\
(\psi_{j,k,l}^N)^*&:=q^{-N+j-l}\sqrt{[j,k,l]!}\,(z_1^jz_2^kz_3^l)^*\;,
   && \textup{if}\;N\leq 0\; \quad\textup{and with}\quad \; j+k+l=-N\,.
\end{align*}
Let $\Psi_N$ be the column vector with components $\psi_{j,k,l}^N$
and $P_N$ be the projection -- of size $d_N:=\frac{1}{2}(|N|+1)(|N|+2)$ -- given by
\begin{equation}\label{mon-pro}
P_N:=\Psi_N\Psi_N^\dag \;. 
\end{equation}
Then one has $\Sigma_{0,N}\simeq P_N\Aq^{d_N}$ as right $\Aq$-modules
and $\Sigma_{0,N}\simeq \Aq^{d_N}P_{-N}$ as left $\Aq$-modules.
\end{prop}

\begin{proof}
The column vector $\Psi_N$  has entries in number $d_N=\frac{1}{2}(|N|+1)(|N|+2)$ and $\Psi_N^\dag$ is a row vector of the same size. By the previous lemma $\Psi^\dag_N\Psi_N=1$, so $P_N:=\Psi_N\Psi_N^\dag$ is a projection. Next, consider the right $\Aq$-module map
$$
P_N\Aq^{d_N}\to\Oq\;,\qquad
 v=(v_{j,k,l}) \mapsto \Psi_N^\dag\cdot v=\sum_{i+j+k=N}(\psi^N_{j,k,l})^*v_{j,k,l} \;.
$$
Since each $z_i$ is $\U_q(\mathfrak{su}(2))$-invariant and 
$z_i\za K_2=q^{\frac{1}{2}}z_i$, it follows that  
$\Psi_N^\dag\za h=\epsilon(h)\Psi_N^\dag$ for all $h\in\U_q(\mathfrak{su}(2))$ and
$\Psi_N^\dag\za K_1K_2^2=q^N\Psi_N^\dag$. Hence, the image of $v$ is in
$\Sigma_{0,N}$. The inverse $\Aq$-module map is
$$
\Sigma_{0,N}\to P_N\Aq^{d_N}\;,\qquad a\mapsto \Psi_Na \;,
$$
thus proving that we have an isomorphism of  right  $\Aq$-modules.
The proof for the left module structure is completely analogous,
the two module maps being in this case $\,v\mapsto 
\sum_{i+j+k=N}v_{j,k,l}\psi^{-N}_{j,k,l}\,$ and $a\mapsto a\Psi_{-N}^\dag$. 
\end{proof}

It turns out that what we just computed are the elements $t^{0,N}_{\uparrow,\underline{i}}$ 
for $ N \geq 0$ (resp.~$t^{-N,0}_{\uparrow,\underline{i}}$ if $N\leq 0$). 
The next proposition yields the exact relation with the matrix entries of $\Psi_N^\dag$.

\begin{prop}\label{lemma:psit}
It holds that
$$
\psi^{N}_{j,k,l}=
\begin{cases}
(t^{0,N}_{0,\underline{i}})^*\;,
&\mathrm{with}\quad\underline{i}=\left(0,j+k,\tfrac{1}{2}(k-j)\right)\quad \textup{and for all}\;N\geq 0\,, \\
(t^{-N,0}_{0,\underline{i}})^*\;,
&\mathrm{with}\quad\underline{i}=\left(j+k,0,\tfrac{1}{2}(j-k)\right)\quad \textup{and for all}\;N\leq 0\,.
\end{cases}
$$
\end{prop}
\begin{proof}
For $N\geq 0$ and $\underline{i}=(0,2l,m)$ the generic label of the
representation $(0,N)$, definition \eqref{eq:eqdef} gives
$$
t^{0,N}_{0,\underline{i}}=X_{\underline{i}}^{(0,N)}\az z_3^N=
\frac{1}{[2l]!}\sqrt{\frac{[N-2l]!}{[N]!}\,\frac{[l+m]!}{[l-m]!}}\,
F_1^{l-m}F_2^{2l} \az z_3^N \,.
$$
Induction on $N$  yields
$$
F_2\az z_3^N=q^{-\frac{N-1}{2}}[N]z_2z_3^{N-1}\;,
$$
and induction on $l$  yields
$$
F_2^{2l}\az z_3^N=q^{-\frac{N-1}{2}}q^{-l(2l-1)}\frac{[N]!}{[N-2l]!} \, z_2^{2l}z_3^{N-2l} \;.
$$
Similarly, changing the labels $3\to 2$ and $2\to 1$:
$$
F_1^{l-m}\az z_2^{2l}=q^{-\frac{2l-1}{2}}q^{-\frac{1}{2}(l-m)(l-m-1)}\frac{[2l]!}{[l+m]!} \,z_1^{l-m}z_2^{l+m} \;.
$$
Thus
\begin{align*}
t^{0,N}_{0,\underline{i}}&=
\sqrt{\tfrac{[N]!}{[N-2l]![l+m]![l-m]!}}\,
q^{-\frac{N-1}{2}}q^{-l(2l-1)}q^{-\frac{2l-1}{2}}q^{-\frac{1}{2}(l-m)(l-m-1)}z_1^{l-m}z_2^{l+m}
z_3^{N-2l} \\ &=(\psi^N_{l-m,l+m,N-2l})^* \;,
\end{align*}
which establishes the case $N\geq 0$.
The proof for the case $N\leq 0$ is similar.
\end{proof}

It is computationally useful to introduce the left action $h\mapsto\mL{h}$, of $\Uq$ on $\Oq$,
given by $\mL{h}a:=a\za S^{-1}(h)$; it satisfies $\mL{x}(ab)=(\mL{x_{(2)}}a)(\mL{x_{(1)}}b)$
due to the presence of the antipode. Also, it is a unitary action for the inner product on
$\Oq$ coming from the Haar state $\varphi$. The proof is a simple computation:
\begin{align*}
\varphi\bigl((\mL{h^*}a)^*b\bigr)
&=\varphi\bigl(\{a\za S^{-1}(h^*)\}^*b\bigr) 
=\varphi\bigl(\{a^*\za h\}b\bigr) 
=\varphi\bigl(\{a^*\za h_{(1)}\}\epsilon(h_{(2)})b\bigr) \\
&=\varphi\bigl(\{a^*\za h_{(1)}\}\{b\za S^{-1}(h_{(3)})h_{(2)}\}\bigr) 
=\varphi\bigl(\{a^*(b\za S^{-1}(h_{(2)}))\}\za h_{(1)}\bigr) \\
&=\epsilon(h_{(1)})\varphi\bigl(a^*\{b\za S^{-1}(h_{(2)})\}\bigr)
=\varphi\bigl(a^*\{b\za S^{-1}(h)\}\bigr) \\
&=\varphi\bigl(a^*(\mL{h}b)\bigr) \;,
\end{align*}
for all $a,b\in\Oq$ and $h\in\Uq$.
With this action the bimodule $\mathfrak{M}(\sigma)$ can be viewed as the set of elements
of $\Oq\otimes\C^{\dim\sigma}$ that are invariant under the action $\mL{h_{(1)}}
\otimes\sigma(h_{(2)})$ of $h\in\Kq$.

\section{Characteristic classes}\label{se:top}

Equivalence classes of finitely generated projective (left or right)  modules over an algebra $\A$ -- the algebraic counterpart of vector bundles -- are elements of the group $K_0(\A)$.  Equivalence classes of even Fredholm modules -- the algebraic counterpart of `fundamental classes' -- gives a dual group $K^0(\A)$. The natural non-degenerate pairing between K-theory and K-homology, which pairs projective modules with even  Fredholm modules, gives  index maps $K_0(\Aq)\to\Z$. For the quantum projective plane, $K_0(\Aq)\simeq\Z^3\simeq K^0(\Aq)$. The result for K-theory can be proved viewing the corresponding $C^*$-algebra as the Cuntz--Krieger algebra of a graph \cite{HS02}. The group $K_0$ is given as the cokernel of the incidence matrix canonically associated with the graph; the result for K-homology can be proven using the same techniques:
the groups $K^0$ is now given as the kernel of the transposed matrix \cite{Cun84}. It is worth mentioning that 
$K_1(\Aq) = K^1(\Aq) = 0$ with the group $K_1$ (resp. $K^1$) given as the kernel (resp. the cokernel) of the incidence matrix (resp. the transposed matrix).

Thus a finitely generated projective (left or right)  module over $\Aq$ is uniquely identified by three integers and these are obtained by pairing the corresponding idempotent with the three generators of the K-homology. 
This section is devoted to the explicit construction of three Fredholm modules that, in the next section will be shown to be generators of the K-homology by pairing them with suitable idempotents.

\subsection{Fredholm modules and their characters}\label{se:fm}~\\
We recall that a $(k+1)$-summable even Fredholm module for the algebra $\A$ is a triple $(\pi,\HH,F)$ consisting of a $\Z_2$-graded Hilbert space $\HH=\HH^+\oplus\HH^-$, a graded representation 
$\pi=\pi^+ \oplus \pi^- : \A\to\B(\HH^+)\oplus\B(\HH^-)$ and an odd operator $F$ such that $[F,a_0][F,a_1]\ldots [F,a_k]$ is of traceclass for any $a_0,\ldots,a_k\in\A$. 
With a $(k+1)$-summable even Fredholm module there are canonically associated even cyclic cocycles 
$\mathrm{ch}^n_{(\pi,\HH,F)}$, for $2n\geq k$. The map 
$\mathrm{ch}^n_{(\pi,\HH,F)}:\A^{2n+1}\to\C$ is given by
$$
\mathrm{ch}^n_{(\pi,\HH,F)}(a_0,\ldots,a_{2n+1}):=
\tfrac{1}{2}(-1)^n\tr_{\HH}\left( \gamma F[F,a_0][F,a_1]\ldots [F,a_{2n}] \right)
$$
where $\gamma$ is the grading and the symbol $\pi$ for the representation is understood.
We recall that a cyclic $2n$-cocycle $\tau_n$ is a $\C$-linear map $\tau_n:\A^{2n+1}\to\C$ which is cyclic, i.e.~it satisfies
$$
\tau_n(a_0,a_1,\ldots,a_{2n})=\tau_n(a_{2n},a_0,\ldots,a_{n-1})\;,
$$
and is a Hochschild coboundary, i.e.~it satisfies $b\,\tau_n=0$, with $b$ the boundary operator:
$$
b\,\tau_n(a_0,\ldots,a_{2n+1}):=
\sum_{j=0}^{2n}(-1)^j\tau(a_0,\ldots,a_ja_{j+1},\ldots,a_{2n+1})
-\tau(a_{2n+1}a_0,a_1,\ldots,a_{2n}) \;.
$$

When applied to an idempotent one gets a pairing
\begin{align}
\inner{\,,\,} &: K^0(\A)\times K_0(\A)\to\Z \;, \notag\\
\inner{[(\pi,\HH,F)],[e]} &:= \mathrm{ch}^n_{(\pi,\HH,F)} ([e]) =\tfrac{1}{2}(-1)^n\tr_{\HH\otimes\C^m}(\gamma F[F,e]^{2n+1}), \label{eq:KKpair}
\end{align}
where $m$ is the size of $e$ and matrix multiplication is understood. 
The pairing is integer valued, being the index of the Fredholm operator 
$\pi^-_m(e) F_m \pi^+_m(e)$ from $\pi^+_m(e) \HH^+_m$ to $\pi^-_m(e) \HH^-_m$, with the natural extensions
$$
\HH_m = \HH \otimes \C^m \,, \quad F_m = F \otimes 1 \,,  
$$
forming an even Fredholm module over $M_m=\A\otimes M_m(\C)$.
The result of the pairing depends only on the classes of $e$ and of the Fredholm module 
(for details see \cite{Con94}).

To construct three independent Fredholm modules we need (at least) four representations. Four is exactly the number of \emph{irreducible} representations of the algebra $\Aq$, that we describe in the next section. It is peculiar  of the quantum case that one we needs to consider only irreducible representations: at $q=1$ irreducible representation are $1$-dimensional and give only one of the generators of the K-homology (the trivial Fredholm module). An additional true `quantum effect' is that for $\CP^2_q$ one gets three independent $1$-summable Fredholm modules corresponding to independent traces on $\Aq$. Thus all relevant information leaves in degree zero. In contrast, for the classical $\CP^2$ one needs (there exist) cohomology classes in degree zero, two and four.

Both at the algebraic and at the $C^*$-algebraic levels there is a sequence
\begin{equation}\label{eq:seq}
\Aq \to \A(\CP^1_q) \to \A(\{\mathrm{pt}\})=\C \to \A(\emptyset)=0 \;,
\end{equation}
of $*$-algebra morphisms. In reverse order: the empty space, the space with one point, the quantum projective line and the quantum projective plane.

The first map is the quotient of $\Aq$ by the ideal generated by $p_{1i}$ and $p_{i1}$: roughly speaking, by putting $p_{1i}=p_{i1}=0$ the remaining generators satisfy the relations of $\CP^1_q$, the quantum projective line -- also known as the standard Podle\'s sphere --; and any $*$-representation of $\Aq$ with $p_{1i}$ in the kernel comes from a representation of $\CP^1_q$. The second map is the further quotient by the ideal generated by $p_{2i}$ and $p_{i2}$: roughly speaking, in $\Aq$ we send $p_{ij}\mapsto\chi_0(p_{ij})=\delta_{i3}\delta_{j3}$ and get the algebra $\C$; this is the only non-trivial character $\chi_0$ of the algebra $\Aq$, the `classical point' of $\CP^2_q$. The last map is simply $1\mapsto 0$.

\subsection{The rank and the 1st Chern number of a projective module} ~\\
For $\A(\emptyset)$ there is only one irreducible representation, the trivial one, not enough to construct
a Fredholm module (and indeed, $K^0(0)=0$).

For $\C=\A(\{\mathrm{pt}\})$ there is the representation coming from the morphism $\C\to 0$, and one further faithful irreducible representation given by the identity map $c\mapsto c$. These two are all\footnote{Irreducible $*$-representations of $\C$ are $1$-dimensional (it is abelian). By linearity, they are in $1$ to $1$ correspondence with maps $1\to p$, with $p\in\C$ a projection. Hence $p=0$ or $1$, and the corresponding representations are the trivial map $c\mapsto 0$, and the identity map $c\mapsto c$.} the irreducible $*$-representations of $\C$, and are enough to construct an even Fredholm module, $$
\HH_0=\C\oplus\C\;,\qquad
\pi_0(c)=\maa{14pt}{c & 0 \\ 0 & 0} \;,\qquad
F_0=\maa{14pt}{0 & 1 \\ 1 & 0} \;.
$$
which generates $K^0(\C)=Z$. 
The pairing \eqref{eq:KKpair} with $K_0(\C)=\Z$ (projections in 
$\mathrm{Mat}_\infty(\C)$ are equivalent if{}f they have the same rank)
is given by the matrix trace:
$$
\mathrm{ch}^0_{(\pi_0,\HH_0,F_0)}([e])=\tr(e).
$$

The pullback of this Fredholm module to $\CP^2_q$ just substitutes the representation $\pi_0$ of $\C$ with the character $\chi_0:\Aq\to\C$ of $\Aq$. Using the same symbol for the character, one gets  
a map
\begin{equation}\label{eq:rank}
\mathrm{ch}^0_{(\pi_0,\HH_0,F_0)}:K_0(\Aq)\to\Z \;,\qquad
\mathrm{ch}^0_{(\pi_0,\HH_0,F_0)}([e])=\tr\,\chi_0(e) \;,
\end{equation}
The geometrical meaning is the following: the rank of a vector bundle is the dimension of the fiber at any point $x$ of the space, and this coincides with the trace of the corresponding projection evaluated at $x$. Here we have only one `classical point', and the map in \eqref{eq:rank} computes the rank of the restriction of the vector bundle to this classical point. Notice that if the module is free, then (\ref{eq:rank}) is really the rank of the module.

We pass to the next level of the sequence (\ref{eq:seq}), that is the algebra $\A(\CP^1_q)$. 
There are two irreducible $*$-representations coming from the map in the sequence -- the trivial one and the character of the algebra --, and one further representation that is faithful and irreducible. The latter is the restriction of  the representation $\chi_1:\Sq\to\B(\ell^2(\N))$ given by
\begin{subequations}\label{eq:chione}
\begin{align}
\chi_1(z_1) &=0 \;,\\
\chi_1(z_2)\ket{n} &=q^n\ket{n} \;, \label{eq:one}\\
\chi_1(z_3)\ket{n} &=\sqrt{1-q^{2(n+1)}}\ket{n+1} \;, \label{eq:two}
\end{align}
\end{subequations}
A potential Fredholm module is given by
$$
\HH_1=\ell^2(\N)\oplus\ell^2(\N)\;,\qquad
\pi_1(a)=\maa{14pt}{\chi_1(a) & 0 \\ 0 & \chi_0(a) \,\id_{\ell^2(\N)} } \;,\qquad
F_1=\maa{14pt}{0 & 1 \\ 1 & 0} \;.
$$
However, this is \emph{not} a Fredholm module over $\Sq$, since
$$
\gamma F_1[F_1,\pi(z_3)]=\bigl\{\chi_1(z_3)-\chi_0(z_3)\bigr\}\maa{14pt}{1 & 0 \\ 0 & 1}
$$
is not compact. On the other hand, this \emph{is} a Fredholm module when restricted to $\Aq$ 
(or $\A(\CP^1_q)$). Indeed, due to the tracial relation $q^4p_{11}+q^2p_{22}+p_{33}=1$,
a complete set of generators for $\Aq$ is made of $\{p_{11},p_{12},p_{13},p_{22},p_{23}\}$ and their adjoints. On these generators the representation $\chi_1$ is
\begin{subequations}\label{eq:chionep}
\begin{align}
\chi_1(p_{1i}) &=0 \;,\qquad\forall\;i=1,2,3\,,\\
\chi_1(p_{22})\ket{n} &=q^{2n}\ket{n} \;, \\
\chi_1(p_{23})\ket{n} &=q^{n+1}\sqrt{1-q^{2(n+1)}}\ket{n+1} \;,
\end{align}
\end{subequations}
and $\chi_1(p_{ij})-\chi_0(p_{ij})$ is of trace class for all $i,j$. In particular, this means that the Fredholm module is $1$-summable.
The associated 
character
is
\begin{equation}\label{eq:charge}
\mathrm{ch}^0_{(\pi_1,\HH_1,F_1)}:K_0(\Aq)\to\Z \;,\qquad
\mathrm{ch}^0_{(\pi_1,\HH_1,F_1)}([e])=\tr_{\ell^2(\N)\otimes\C^m}\,(\chi_1-\chi_0)(e) \;,
\end{equation}
where $m$ is the size of the matrix $e$. The value in (\ref{eq:charge}) depends only on the restriction of the `vector bundle' to the subspace $\CP^1_q$, and will be called for this reason the \emph{1st Chern number} (or also the \emph{monopole charge}).

\subsection{The 2nd Chern number of a projective module}\label{se:inpm}~\\
Besides the representations (and the Fredholm modules) coming from the sequence (\ref{eq:seq}), the algebra $\Aq$ has a further irreducible representation and a further Fredholm module, which is independent from the previous two. The representation is faithful  and comes from the representation $\chi_2:\Sq\to\B(\ell^2(\N^2))$ given by
\begin{subequations}\label{eq:chitwo}
\begin{align}
\chi_2(z_1)\ket{k_1,k_2} &:=q^{k_1+k_2}\ket{k_1,k_2} \;,\\
\chi_2(z_2)\ket{k_1,k_2} &:=q^{k_1}\sqrt{1-q^{2(k_2+1)}}\ket{k_1,k_2+1} \;,\\
\chi_2(z_3)\ket{k_1,k_2} &:=\sqrt{1-q^{2(k_1+1)}}\ket{k_1+1,k_2} \;.
\end{align}
\end{subequations}

The construction of the Fredholm module is a bit involved.
Let us use the labels $\ell=\frac{1}{2}(k_1+k_2)$ and $m=\frac{1}{2}(k_1-k_2)$.
A new basis for the Hilbert space is then given by $\ket{\ell,m}$ with
$\ell\in\frac{1}{2}\N$ and $m=-\ell,-\ell+1,\ldots,\ell$. In this basis the representation reads
\begin{align*}
\chi_2(z_1)\ket{\ell,m} &=q^{2\ell}\ket{\ell,m} \;,\\
\chi_2(z_2)\ket{\ell,m} &=q^{\ell+m}\sqrt{1-q^{2(\ell-m+1)}}\ket{\ell+\oh,m-\oh} \;,\\
\chi_2(z_3)\ket{\ell,m} &=\sqrt{1-q^{2(\ell+m+1)}}\ket{\ell+\oh,m+\oh} \;.
\end{align*}
For the third Fredholm module we take as Hilbert space $\HH_2$ two copies of  the linear span of orthonormal vectors $\ket{\ell,m}$, with $\ell\in\frac{1}{2}\N$ and $\ell+m\in\N$.
The grading $\gamma_2$ and the operator $F_2$ are the obvious ones. It remains to describe the representation
$\pi_2=\pi_+\oplus\pi_-$. As subrepresentation $\pi_+$ we choose
$$
\pi_+(a)\ket{\ell,m}:=\begin{cases}
\chi_2(a)\ket{\ell,m}         &\mathrm{if}\;m\leq\ell \;,\\
\chi_0(a)\ket{\ell,m}         &\mathrm{if}\;m>\ell \;.
\end{cases}
$$
One checks that modulo traceclass operators:
\begin{align*}
\pi_+(p_{11}) &\sim\pi_+(p_{12})\sim\pi_+(p_{13})\sim 0 \;,\\
\pi_+(p_{22})\ket{\ell,m} &\sim\begin{cases}
q^{2(\ell+m)}\ket{\ell,m}    &\mathrm{if}\;m\leq\ell \;,\\
0                            &\mathrm{if}\;m>\ell \;,
\end{cases} \\
\pi_+(p_{23})\ket{\ell,m} &\sim\begin{cases}
q^{\ell+m+1}\sqrt{1-q^{2(\ell+m+1)}}\ket{\ell,m+1}     &\mathrm{if}\;m\leq\ell-1 \;,\\
0                                                      &\mathrm{if}\;m\geq\ell \;.
\end{cases}
\end{align*}
We define the subrepresentation $\pi_-$ by adding multiplicities to $\chi_1$. On the generators:
\begin{align*}
\pi_-(p_{11}) &=\pi_-(p_{12})=\pi_-(p_{13})=0 \;,\\
\pi_-(p_{22})\ket{\ell,m} &=q^{2(\ell+m)}\ket{\ell,m} \;,\quad
\pi_-(p_{23})\ket{\ell,m} &=q^{\ell+m+1}\sqrt{1-q^{2(\ell+m+1)}}\ket{\ell,m+1} \;.
\end{align*}
On each invariant subspace with a fixed $\ell$, putting $n=\ell+m$ one recovers the representation $\chi_1$. Since $\sum_{m>\ell}q^{2(\ell+m)}=(1-q^4)^{-2}$ is finite, on the subspace $m>\ell$ the operators $\pi_-(p_{22})$ and $\pi_-(p_{23})$ are trace class, and so $\pi_+(a)-\pi_-(a)$ is of trace class as well for all $a\in\Aq$: the Fredholm module is $1$-summable with corresponding character
\begin{equation}\label{eq:instnumb}
\mathrm{ch}^0_{(\pi_2,\HH_2,F_2)}:K_0(\Aq)\to\Z \;,\qquad
\mathrm{ch}^0_{(\pi_2,\HH_2,F_2)}([e])=\tr_{\HH_2\otimes\C^m}\,(\pi_+-\pi_-)(e) \;,
\end{equation}
where $m$ is the size of the matrix $e$. The above replaces the 2nd Chern class of the module.

Working with generators and relations as done in \cite{HL04} for quantum spheres, it is not difficult to prove 
that any irreducible $*$-representation of $\Aq$ is equivalent to one of the representations described above.
For completeness, we give the proof in App.~\ref{app-IR}.

By iterating the construction, for any positive integer $n$ one obtains the $n+2$ irreducible representations of the quantum projective spaces $\CP^n_q$ (only one of these is faithful, the others coming from the morphism $\A(\CP^n_q)\to\A(\CP^{n-1}_q)$), and the corresponding $n+1$ Fredholm modules; details will be reported elsewhere \cite{DL09a}.

\subsection{Chern numbers of line bundles}~\\
We know from Sect.~\ref{se:lb} that the bimodule $\Sigma_{0,N}:=\mathfrak{M}(\sigma_{0,N})$   
is isomorphic to $\Aq^{d_N}P_{-N}$ as left module and to $P_N\Aq^{d_N}$ as right module, with $d_N=\frac{1}{2}(|N|+1)(|N|+2)$ and $P_N:=\Psi_N\Psi_N^\dag$ given by Prop.~\ref{lemma:linebundles}. We next compute rank, and 1st and 2nd Chern numbers of $P_N$. We focus the discussion on $N\geq 0$, the case $N\leq 0$ being similar.

Since $\chi_0(P_N)_{j,k,l|j',k',l'}=\delta_{j,0}\delta_{k,0}\delta_{l,N}
\delta_{j',0}\delta_{k',0}\delta_{l',N}$, the rank given by (\ref{eq:rank}) results in
$$
\mathrm{ch}^0_{(\pi_0,\HH_0,F_0)}([P_N])=1 \;,
$$
thus justifying the name `line bundles' for the virtual bundles underlying the modules
$\Sigma_{0,N}$. The same result is valid when $N\leq 0$.

Next, we compute
\begin{equation}\label{eq:trpN}
\tr_{\C^{d_N}}(P_N)=\sum\nolimits_{j+k+l=N}q^{-(jk+kl+lj)}\frac{[N]!}{[j]![k]![l]!}
(z_3^l)^*(z_2^k)^*(z_1^j)^*z_1^jz_2^kz_3^l
\end{equation}
and
\begin{align*}
\chi_1\left( \tr_{\C^{d_N}}(P_N) \right) 
    &=\sum\nolimits_{k+l=N}q^{-kl}\frac{[N]!}{[k]![l]!}
    \chi_1(z_3^*)^l\chi_1(z_2^*)^k\chi_1(z_2)^k\chi_1(z_3)^l \;,\\
    ~\\
\chi_0\left( \tr_{\C^{d_N}}(P_N) \right)
    &=\chi_0(z_3^*)^N\chi_0(z_3)^N=1 \;.
\end{align*}
The use of (\ref{eq:one}) and (\ref{eq:two}) leads to 
\begin{align*}
x_n &:=\inner{n\big|\chi_1\smash[t]{\tr_{\C^{d_N}}}(P_N)\big|n} \\
& =\sum\nolimits_{k+l=N}q^{-kl}\frac{[N]!}{[k]![l]!}
q^{2k(n+l)}(1-q^{2(n+1)})(1-q^{2(n+2)})\ldots(1-q^{2(n+l)}) \\
&=\sum\nolimits_{k=0}^Nq^{2kn}+O(q)=
\begin{cases}
\frac{1-q^{2n(N+1)}}{1-q^{2n}}+O(q)=1+O(q) &\mathrm{if}\;n>0\;,\\
N+1 +O(q)  &\mathrm{if}\;n=0\;.
\end{cases}
\end{align*}
As usual we can compute the index for $q\to 0^+$ (see \cite{DD06,DDLW07}); we get
$$
\sum_{n\geq 0}(x_n-1)=\sum_{n\geq 0}\lim_{q\to 0^+}(x_n-1)=(x_0-1)|_{q=0}=N \;.
$$
The very same result holds for $N\leq 0$. So, the 1st Chern number is 
$$
\mathrm{ch}^0_{(\pi_1,\HH_1,F_1)}([P_N])=N.
$$

It remains to compute the last Chern number.
The representations $\pi_\pm$ in the third Fredholm module are the restriction of (homonymous) representations of $\Sq$ given as follows. With $\chi_0$ and $\chi_2$ the representation described in the previous section, one has that $\pi_+(z_i)\ket{\ell,m}=\chi_2(z_i)\ket{\ell,m}$ 
if $m\leq \ell$, and $\pi_+(z_i)\ket{\ell,m}=\chi_0(z_i)\ket{\ell,m}$ if $m>\ell$. On the other hand, for any $m,\ell$ it holds that 
\begin{align*}
\pi_-(z_1) &=0 \;,\\
\pi_-(z_2)\ket{\ell,m} &=q^{\ell+m}\ket{\ell,m} \;,\\
\pi_-(z_3)\ket{\ell,m} &=\sqrt{1-q^{2(\ell+m+1)}}\ket{\ell,m+1} \;.
\end{align*}
The use of  \eqref{eq:trpN} yields:
\begin{align*}
\bigl<\ell,m\big| \pi_-\tr_{\C^{d_N}}(P_N) \big|\ell,m\bigr>
&=\sum\nolimits_{r+s=N} \, q^{-rs}\frac{[N]!}{[r]![s]!}
q^{2r(\ell+m+s)}\prod\nolimits_{k=1}^s(1-q^{2(\ell+m+k)}) \\
&=\sum\nolimits_{r=0}^N \, q^{2r(\ell+m)}\{1+O(q)\}=1+N\delta_{\ell,-m}+O(q)
\;,\\
\end{align*}
for the representation $\pi_-$; for the representation $\pi_+$ one gets instead 
$$
\bigl<\ell,m\big| \pi_+\tr_{\C^{d_N}}(P_N) \big|\ell,m\bigr> = 1 \,  
$$
if $m>\ell $, and   
\begin{align*}
\bigl<\ell,m\big| \pi_+\tr_{\C^{d_N}}(P_N) \big|\ell,m\bigr>
&=\sum\nolimits_{i+j+k=N} \, q^{-(ij+jk+ki)}\frac{[N]!}{[i]![j]![k]!}
q^{2i(2\ell+j+k)}q^{2j(\ell+m+k)} \times \\ 
 &\quad\quad\quad\quad\quad \times\; \prod\nolimits_{r=1}^j(1-q^{2(\ell-m+r)})
\prod\nolimits_{s=1}^k(1-q^{2(\ell+m+s)}) \\
 &=\sum\nolimits_{i+j+k=N} \, q^{4i\ell+2j(\ell+m)}+O(q) \\
&=1+N\delta_{\ell,-m}+\tfrac{1}{2}N(N+1)\delta_{\ell,0}+O(q),
\end{align*}
if $m\leq\ell $. In the computation we have used $\frac{[N]!}{[i]![j]![k]!}=q^{-(ij+jk+ki)}\bigl(1+O(q)\bigr)$.
Putting things together results in
$$
\bigl<\ell,m\big| \pi_+\tr_{\C^{d_N}}(P_N)-\pi_-\tr_{\C^{d_N}}(P_N) \big|\ell,m\bigr>
=\begin{cases}
O(q)  &\mathrm{if}\;m>\ell \;. \\
\tfrac{1}{2}N(N+1)\delta_{\ell,0}+O(q) &\mathrm{if}\;m\leq\ell \;.
\end{cases}
$$
Again, we compute in the limit
for $q\to 0^+$ to get
\begin{align*}
\mathrm{ch}^0_{(\pi_2,\HH_2,F_2)}([P_N]) &=\lim_{q\to 0^+}
\sum\nolimits_{\ell\in\frac{1}{2}\N}\sum\nolimits_{m=-\ell}^\ell
\tfrac{1}{2}N(N+1)\delta_{\ell,0}+O(q) \\
& = \tfrac{1}{2}N(N+1) \;.
\end{align*}
The same formula is valid for $N\leq 0$. We summarize the results in a proposition.

\begin{prop}\label{prop:same}
For any $N\in\Z$, the (right) module $\Sigma_{0,N}$ has `rank' 
$$
\mathrm{ch}^0_{(\pi_0,\HH_0,F_0)}([P_N]) = 1 \; ,
$$
`1st Chern number' 
$$
\mathrm{ch}^0_{(\pi_1,\HH_1,F_1)}([P_N]) = N \; ,
$$ 
and `2nd Chern number' 
$$
\mathrm{ch}^0_{(\pi_2,\HH_2,F_2)}([P_N]) = \tfrac{1}{2}N(N+1) \; .
$$
\end{prop}

\begin{cor}
A complete set of generators $\{e_1,e_2,e_3\}$ of $K_0(\Aq)$ is given by 
\begin{list}{}{%
\addtolength{\leftmargin}{-5mm}
\addtolength{\itemsep}{2pt}%
}
\item
$e_1=[1]$ the class of the rank one free (left or right) $\Aq$-module $\Sigma_{0,0}$, 
\item
$e_2$ the class of the left module $\Sigma_{0,1}$ (or equivalently of the right module $\Sigma_{0,-1}$),
\item
$e_3$ the class of the left module $\Sigma_{0,-1}$ (or equivalently  of the right module $\Sigma_{0,1}$).
\end{list}
Thus, in terms of left modules, $e_2$ is the class of the tautological bundle and $e_3$ the class of the dual vector bundle. 
A complete set of generators of $K^0(\Aq)$ is given by the classes of the Fredholm modules $(\pi_i,\HH_i,F_i)$, $i=0,1,2$, given in Sect.~\ref{se:top}.
\end{cor}

\begin{proof}
We have already mentioned that $K_0(\Aq)\simeq K^0(\Aq)\simeq\Z^3$. 
A set of generators of the abelian group $\Z^3$ is the same as a basis of $\Z^3$ as a $\Z$-module.

Suppose we have three elements $e_i\in\Z^3$ and three elements $\varphi_i:\Z^3\to\Z$ in the dual space  of $\Z$-linear maps; call $g\in\mathrm{Mat}(3,\Z)$ the matrix with elements $g_{ij}=\varphi_i(e_j)$. 
Assume that $\det g\neq 0$ and that the inverse $g^{-1}=((g^{ij}))$ of $g$ is an element of $GL(3,\Z)$. Then for any linear map $\psi$ the difference $\psi-\sum\nolimits_{ij}\psi(e_i)g^{ij}\varphi_j$ vanishes on all the $e_i$'s and, by the linear independence  over $\Z$  of the $e_i$'s, we deduce that any element $\psi\in(\Z^3)^*$ is a sum
$$
\psi=\sum\nolimits_{ij}\psi(e_i)g^{ij}\varphi_j
$$
with integer coefficients $\psi(e_i)g^{ij}\in\Z$. Hence, the $\varphi_j$'s are a basis of $(\Z^3)^*$. Similarly for any $v\in\Z^3$ the difference $v-\sum\nolimits_{ij}e_ig^{ij}\varphi_j(v)$ is in the kernel of all $\varphi_i$, meaning that the $e_i$'s are a basis of $\Z^3$ over $\Z$.

Now, let $e_i\in K_0(\Aq)$ and $\varphi_i=[(\pi_{i-1},\HH_{i-1},F_{i-1})]\in K^0(\Aq)$, $i=1,2,3$, be the classes in the Corollary.
By Prop.~\ref{prop:same} the matrix $g$ is given by
$$
g=\maa{22pt}{1 & 1 & 1 \\ 0 & -1 & 1 \\ 0 & 0 & 1} \;.
$$
It is invertible in $GL(3,\Z)$ with inverse
$$
g^{-1}=\maa{22pt}{1 & 1 & -2 \\ 0 & -1 & 1 \\ 0 & 0 & 1} \;.
$$
This proves that we have generators of K-theory and K-homology.
\end{proof}

\section{The differential calculus}\label{se:diff}

On the line bundles described previously we aim to define and study (anti-)selfdual connections. 
To this end we need the full machinery of a differential calculus on $\CP^2_q$ and additional ingredients 
like a Hodge star operator. We start with forms.

\subsection[The differential calculus]{Holomorphic and antiholomorphic forms}\label{se:hah}~\\
In \cite{DDL08b} we studied the antiholomorphic part of the differential calculus on $\CP^2_q$.
Antiholomorphic forms were defined as 
$$
\Omega^{0,0}(\CP^2_q):=\Aq\, , \quad  
\Omega^{0,1}(\CP^2_q):=\Sigma_{\frac{1}{2},\frac{3}{2}} 
\quad \textup{and} \quad \Omega^{0,2}(\CP^2_q):=\Sigma_{0,3}\,.
$$
We complete the calculus presently. Each bimodule of forms $\Omega^{i,j}(\CP^2_q)$ will be identified with a suitable 
bimodule $\mathfrak{M}(\sigma)$ of equivariant elements as described in the previous section.
That is to say, $\Omega^{i,j}(\CP^2_q)=\mathfrak{M}(\sigma^{i,j})=\Oq\!\boxtimes_{\sigma^{i,j}}
V^{i,j}$ for a (not necessarily irreducible) representation $\sigma^{i,j}:\Kq\to\mathrm{Aut}(V^{i,j})$  of the
Hopf algebra $\Kq$ on  $V^{i,j}\simeq\C^{\dim\sigma^{i,j}}$.
The representations $\sigma^{0,j}$ are known from \cite{DDL08b} and the 
$\sigma^{j,0}$ are obtained by conjugation. The dimension of $\sigma^{i,j}$
does not depend on $q$, and for $q=1$ gives the rank $r=\binom{2}{i}\binom{2}{j}$ of $\Omega^{i,j}(\CP^2_q)$. 
There is a unique sub-representation
of the Hopf tensor product $\sigma^{i,0}\otimes\sigma^{0,j}$ with the
prescribed dimension: this allows one to identify $\sigma^{i,j}$.
The representations relevant for the calculus are listed in the following figure:
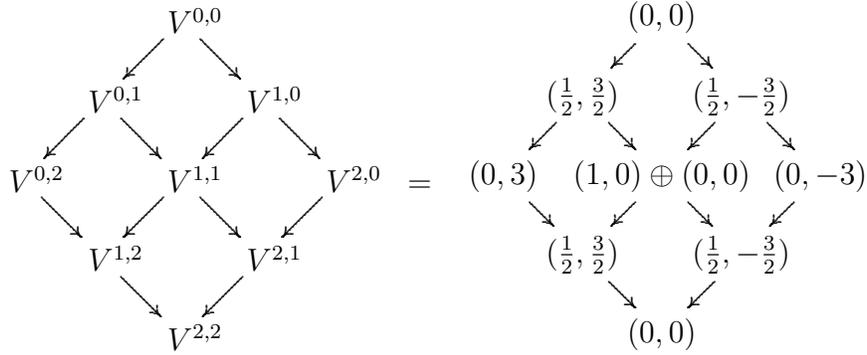
\begin{center}\begin{tabular}{ccc}
\begindc{\commdiag}[30]
 \obj(3,5)[A]{$\;\;\;V^{0,0}$}
 \obj(2,4)[B]{$\;\;\;V^{0,1}$}
 \obj(4,4)[C]{$\;\;\;V^{1,0}$}
 \obj(1,3)[D]{$\;\;\;V^{0,2}$}
 \obj(3,3)[E]{$\;\;\;V^{1,1}$}
 \obj(5,3)[F]{$\;\;\;V^{2,0}$}
 \obj(2,2)[G]{$\;\;\;V^{1,2}$}
 \obj(4,2)[H]{$\;\;\;V^{2,1}$}
 \obj(3,1)[I]{$\;\;\;V^{2,2}$}
 \mor{A}{B}{}
 \mor{A}{C}{}
 \mor{B}{D}{}
 \mor{B}{E}{}
 \mor{C}{E}{}
 \mor{C}{F}{}
 \mor{D}{G}{}
 \mor{F}{H}{}
 \mor{E}{G}{}
 \mor{E}{H}{}
 \mor{G}{I}{}
 \mor{H}{I}{}
\enddc
& \raisebox{2.1cm}{$=$} &
\begindc{\commdiag}[30]
 \obj(3,5)[A]{$(0,0)$}
 \obj(2,4)[B]{$(\frac{1}{2},\frac{3}{2})$}
 \obj(4,4)[C]{$(\frac{1}{2},-\frac{3}{2})$}
 \obj(1,3)[D]{$(0,3)$}
 \obj(3,3)[E]{$(1,0)\oplus (0,0)$}
 \obj(5,3)[F]{$(0,-3)$}
 \obj(2,2)[G]{$(\frac{1}{2},\frac{3}{2})$}
 \obj(4,2)[H]{$(\frac{1}{2},-\frac{3}{2})$}
 \obj(3,1)[I]{$(0,0)$}
 \mor{A}{B}{}
 \mor{A}{C}{}
 \mor{B}{D}{}
 \mor{B}{E}{}
 \mor{C}{E}{}
 \mor{C}{F}{}
 \mor{D}{G}{}
 \mor{F}{H}{}
 \mor{E}{G}{}
 \mor{E}{H}{}
 \mor{G}{I}{}
 \mor{H}{I}{}
\enddc
\end{tabular}\end{center}
In the diamond on the right, in the position $(i,j)$, we give the values of
spin and charge $(\ell,N)$ of the representation $\sigma^{i,j}$.

For later use, we list the representations $\sigma_{\frac{1}{2},N}$ and 
$\sigma_{1,N}$ explicitly: 
\begin{align*}
\sigma_{\frac{1}{2},N}(K_1) &=\maa{14pt}{q^{\frac{1}{2}} & 0 \\ 0 & q^{-\frac{1}{2}}} \;,&
\sigma_{\frac{1}{2},N}(E_1) &=\maa{14pt}{0 & 1 \\ 0 & 0} \;,&
\sigma_{\frac{1}{2},N}(F_1) &=\maa{14pt}{0 & 0 \\ 1 & 0} \;,\\
\rule{0pt}{30pt}
\sigma_{1,N}(K_1) &=\maa{22pt}{q & 0 & 0 \\ 0 & 1 & 0 \\ 0 & 0 & q^{-1}} \;,&
\sigma_{1,N}(E_1) &=[2]^{\frac{1}{2}}\maa{22pt}{0 & 1 & 0 \\ 0 & 0 & 1 \\ 0 & 0 & 0} \;,&
\sigma_{1,N}(F_1) &=[2]^{\frac{1}{2}}\maa{22pt}{0 & 0 & 0 \\ 1 & 0 & 0 \\ 0 & 1 & 0} \;,
\end{align*}
and furthermore, $\sigma_{\frac{1}{2},N}(K_1K_2^2)$ (resp. $\sigma_{1,N}(K_1K_2^2)$)
is $q^N$ times the identity matrix.

\subsection{The wedge product}~\\
We first make  $\Omega^{\bullet,\bullet}(\CP^2_q)=\bigoplus_{i,j}\Omega^{i,j}(\CP^2_q)$ a bi-graded associative algebra.

Let $V^{\bullet,\bullet}=\bigoplus_{i,j}V^{i,j}$, and suppose we have a bi-graded
associative left $\Kq$-covariant product on $V^{\bullet,\bullet}$, denoted $\wprod$.
For $\omega=av$ and $\omega'=a'v'$, with $a,a'\in\Oq$ and $v\in V^{i,j}$, $v'\in V^{i',j'}$,
define
$$
\omega\wprod\omega' := (a\,a')\, (v\wprod v').
$$
Using left $\Kq$-covariance of the product on $V^{\bullet,\bullet}$ and the fact
that $\Oq$ is a $\Kq$-bimodule algebra so that $aa'\za S^{-1}(h)=\{a\za S^{-1}(h_{(2)})\}
\{a\za S^{-1}(h_{(1)})\}$, we get
$$
(\mL{h_{(1)}}\otimes\sigma^{i+i',j+j'}(h_{(2)}))(\omega\wprod\omega')
=(\mL{h_{(2)}}a)
(\mL{h_{(1)}}a')
(\sigma^{i,j}(h_{(3)})v)\wprod
(\sigma^{i,j}(h_{(4)})v') \;.
$$
If $\omega$ is invariant (i.e.~it belongs to $\Omega^{i,j}(\CP^2_q))$, this can be simplified
and becomes
$$
(\mL{h_{(1)}}\otimes\sigma^{i+i',j+j'}(h_{(2)}))(\omega\wprod\omega')
=a(\mL{h_{(1)}}a')v\wprod (\sigma^{i,j}(h_{(2)})v') \;.
$$
If also $\omega'$ is invariant (i.e.~it belongs to $\Omega^{i',j'}(\CP^2_q)$), we
get
$$
(\mL{h_{(1)}}\otimes\sigma^{i+i',j+j'}(h_{(2)}))(\omega\wprod\omega')
= \epsilon(h) \, aa'(v\wprod v')= \epsilon(h) \,\omega\wprod\omega' \;.
$$
Thus, $\wprod$ defines a bilinear map \,$\Omega^{i,j}(\CP^2_q) \,\times\, \Omega^{i',j'}(\CP^2_q)\to\Omega^{i+i',j+j'}(\CP^2_q)$.
Its associativity follows from associativity of both the products in $\Oq$ and $V^{\bullet,\bullet}$.

The datum $(\Omega^{\bullet,\bullet}(\CP^2_q),\wprod)$ is automatically a left $\Aq\rtimes\Uq$-module
algebra, since left and right canonical actions commute and $\Oq$ is a left $\Aq\rtimes\Uq$-module
algebra. This means that
$$
h\az(\omega\wprod\omega')=(h_{(1)}\az\omega)\wprod(h_{(2)}\az\omega')\;,
$$
for all $h\in\Uq$ and $\omega,\omega'\in\Omega(\CP^2_q)^{\bullet,\bullet}$.

All we need is then a graded associative left $\Kq$-covariant
product on $V^{\bullet,\bullet}$. For all $i,j,i',j'$, we shall now construct a left
$\Kq$-module map $\wprod:V^{i,j}\times V^{i',j'}\to V^{i+i',j+j'}$ which is unique
up to some normalization constants. 

When $q=1$, one can fix the normalization, up to some phase factors
and angles, by requiring that these maps are partial isometries. Requiring that 
vectors with real components form
a subalgebra (so that real forms are an algebra), the phases must be $\pm 1$;
the remaining angles and signs are then fixed by the requirement of
associativity and graded commutativity of the product.

For $q\neq 1$, partial isometries do not give an associative product. We determine
in App.~\ref{app:B} the most general value of the normalization constants in
order to have a left $\Kq$-covariant product on $V^{\bullet,\bullet}$ which is
i) associative, ii) graded commutative for $q=1$, and iii) it sends real vectors into
real vectors. Here we just present the result.
\begin{prop}\label{pr:const}
A left $\Kq$-covariant graded associative product $\wprod$ on $V^{\bullet,\bullet}$,
sending real vectors to real vectors and graded commutative for $q=1$, is given by
\begin{align*}
V^{0,1}\times V^{0,1} &\to V^{0,2}\;,
   & v\wprod w &:=c_0\mu_0(v,w)^t\;,\\
V^{0,1}\times V^{1,0} &\to V^{1,1}\;,
   & v\wprod w &:=\bigl(c_1\mu_1(v,w),c_2\mu_0(v,w)\bigr)^t \;,\\
V^{0,1}\times V^{2,1} &\to V^{2,2}\;,
   & v\wprod w &:=c_3\mu_0(v,w)^t\;,\\
V^{0,1}\times V^{1,1} &\to V^{1,2}\;,
   & v\wprod w &:=\frac{c_0}{[2]c_1}\mu_2(v,w)^t-\frac{c_0}{[2]c_2}vw_4 \;,\\
V^{1,0}\times V^{1,0} &\to V^{2,0}\;,
   & v\wprod w &:=c_4\mu_0(v,w)^t\;,\\
V^{1,0}\times V^{0,1} &\to V^{1,1}\;,
   & v\wprod w &:=\bigl(-q^{\frac{1}{2}s}c_1\mu_1(v,w),q^{-\frac{3}{2}s}c_2\mu_0(v,w)\bigr)^t \;,\\
V^{1,0}\times V^{1,2} &\to V^{2,2}\;,
   & v\wprod w &:=\frac{c_3c_4}{c_0}\mu_0(v,w)^t\;,\\
V^{1,0}\times V^{1,1} &\to V^{2,1}\;,
   & v\wprod w &:=-q^{-\frac{1}{2}s}\frac{c_4}{[2]c_1}\mu_2(v,w)^t-q^{\frac{3}{2}s}\frac{c_4}{[2]c_2}vw_4 \;,\\
V^{1,2}\times V^{1,0} &\to V^{2,2}\;,
   & v\wprod w &:=\frac{c_3c_4}{c_0}\mu_0(v,w)^t\;,\\
V^{2,1}\times V^{0,1} &\to V^{2,2}\;,
   & v\wprod w &:=c_3\mu_0(v,w)^t\;,\\
V^{1,1}\times V^{0,1} &\to V^{1,2}\;,
   & v\wprod w &:=-q^{-\frac{1}{2}s}\frac{c_0}{[2]c_1}\mu_3(v,w)^t-q^{\frac{3}{2}s}\frac{c_0}{[2]c_2}v_4w \;,\\
V^{1,1}\times V^{1,0} &\to V^{2,1}\;,
   & v\wprod w &:=\frac{c_4}{[2]c_1}\mu_3(v,w)^t-\frac{c_4}{[2]c_2}v_4w \;,\\
V^{1,1}\times V^{1,1} &\to V^{2,2}\;,
   & v\wprod w&:=-q^{-\frac{1}{2}s}\frac{c_3c_4}{[2]|c_1|^2}\mu_4(v,w)-q^{\frac{3}{2}s}\frac{c_3c_4}{[2]|c_2|^2}v_4 w_4 \;,
\end{align*}
where the maps $\mu_i$'s are 
\begin{align*}
\mu_0:\R^2\times\R^2&\to\R\;, &
\mu_0(v,w)&:=[2]^{-\frac{1}{2}}(q^{\frac{1}{2}}v_1w_2-q^{-\frac{1}{2}}v_2w_1)\;,\\
\mu_1:\R^2\times\R^2&\to\R^3\;, &
\mu_1(v,w)&:=\bigl(v_1w_1,[2]^{-\frac{1}{2}}(q^{-\frac{1}{2}}v_1w_2+q^{\frac{1}{2}}v_2w_1),v_2w_2\bigr)\;,\\
\mu_2:\R^2\times\R^3&\to\R^2\;, &
\mu_2(v,w)&:=\bigl(qv_1w_2-q^{-\frac{1}{2}}[2]^{\frac{1}{2}}v_2w_1,
              q^{\frac{1}{2}}[2]^{\frac{1}{2}}v_1w_3-q^{-1}v_2w_2\bigr)\;,\\
\mu_3:\R^3\times\R^2&\to\R^2\;, &
\mu_3(v,w)&:=\bigl(q^{\frac{1}{2}}[2]^{\frac{1}{2}}v_1w_2-q^{-1}v_2w_1,
              qv_2w_2-q^{-\frac{1}{2}}[2]^{\frac{1}{2}}v_3w_1\bigr)\;,\\
\mu_4:\R^3\times\R^3&\to\R\;, &
\mu_4(v,w)&:=qv_1w_3-v_2w_2+q^{-1}v_3w_1 \;.
\end{align*}
The parameters $c_0,\ldots,c_4\in\R^\times$ and $s=\pm 1$ are not fixed for the time being.
\end{prop}

\medskip

To get an involution we use the fact that the spin $1/2$ (resp.~spin $1$) representation of $\U_q(\mathfrak{su}(2))$ is quaternionic (resp.~real). Rephrased in terms of the representations $\sigma_{\frac{1}{2},N}$ and $\sigma_{1,N}$ of $\Kq$ we have the following lemma, which takes into account the fact that real/quaternionic structures change sign to $N$.

\begin{lemma}\label{lemma:J}
Let $V_{\ell,N}=\C^{2\ell+1}$ be the vector space underlying the representation $\sigma_{\ell,N}$ of $\Kq$.
An antilinear map $J:V_{\ell,N}\to V_{\ell,-N}$ satisfying $J^2=(-1)^{2\ell}$ and such that 
\begin{equation}\label{eq:J}
J\sigma_{\ell,N}(h)=\sigma_{\ell,-N}(S(h)^*)J, 
\end{equation}
for any $h\in\Kq$, is given, for $\ell=0,\frac{1}{2},1$, by
$$
Ja=a^* \;,\quad
J(v_1,v_2)^t=(-q^{-\frac{1}{2}}v_2^*,q^{\frac{1}{2}}v_1^*)^t \;,\quad
J(w_1,w_2,w_3)^t=(-q^{-1}w_3^*,w_2^*,-qw_1^*)^t \;,
$$
for any $a\in V_{0,N}$, $v\in V_{\frac{1}{2},N}$ and $w\in V_{1,N}$ respectively.
Moreover, if $c_0=c_4$, the map
$$
^\star:V^{\bullet,\bullet}\to V^{\bullet,\bullet} \;,\qquad
(v_{i,j}) \mapsto (v^\star)_{i,j}:=(-1)^iJ(v_{j,i}) \;,
$$
is a graded involution, i.e.~it satisfies $ ( v^\star )^\star = v $ and
\begin{equation}\label{eq:claimA}
(v\wprod v')^\star=(-1)^{kk'}v'^\star\wprod v^\star \;,
\end{equation}
for all $v\in V^{i,k-i}$ and $v'\in V^{i',k'-i'}$.
\end{lemma}

\begin{proof}
The property $J^2=(-1)^{2\ell}$ is easily checked in all the above mentioned cases, and the condition 
$\star^2=\mathrm{id}$ is a direct corollary: $\star^2|_{V^{i,j}}=(-1)^{i+j}J^2$, and $V^{i,j}$ is a sum of spaces $V_{\ell,N}$ having $2\ell$ with the same parity of $i+j$. Since $J$ always commutes with 
$\sigma_{\ell,N}(K_1K_2^2)$, and $\sigma_{\ell,N}(K_1K_2^2)=\sigma_{\ell,-N}(S(K_1K_2^2)^*)$, the claim \eqref{eq:J} for $h=K_1K_2^2$ is trivially satisfied. We have to check \eqref{eq:J} for the remaining generators $h=K_1,E_1,F_1$, and in the not-trivial cases $\ell=\frac{1}{2},1$. A direct computation yields:
\begin{align*}
J\sigma_{\frac{1}{2},N}(K_1)J^{-1}&=\maa{14pt}{0 & -q^{-\frac{1}{2}} \\ q^{\frac{1}{2}} & 0}
\maa{14pt}{q^{\frac{1}{2}} & 0 \\ 0 & q^{-\frac{1}{2}}}\maa{14pt}{0 & q^{-\frac{1}{2}} \\ -q^{\frac{1}{2}} & 0}
=\maa{14pt}{q^{-\frac{1}{2}} & 0 \\ 0 & q^{\frac{1}{2}}} \\[5pt]
&=\sigma_{\frac{1}{2},-N}(K_1^{-1}) \;,\\[10pt]
J\sigma_{\frac{1}{2},N}(E_1)J^{-1}&=\maa{14pt}{0 & -q^{-\frac{1}{2}} \\ q^{\frac{1}{2}} & 0}
\maa{14pt}{0 & 1 \\ 0 & 0}\maa{14pt}{0 & q^{-\frac{1}{2}} \\ -q^{\frac{1}{2}} & 0}=\maa{14pt}{0 & 0 \\ -q & 0} \\[5pt]
&=-q\sigma_{\frac{1}{2},-N}(F_1) \;,\\[5pt]
J\sigma_{\frac{1}{2},N}(F_1)J^{-1}&=\maa{14pt}{0 & -q^{-\frac{1}{2}} \\ q^{\frac{1}{2}} & 0}
\maa{14pt}{0 & 0 \\ 1 & 0}\maa{14pt}{0 & q^{-\frac{1}{2}} \\ -q^{\frac{1}{2}} & 0}=\maa{14pt}{0 & -q^{-1} \\ 0 & 0} \\[5pt]
&=-q^{-1}\sigma_{\frac{1}{2},-N}(E_1) \;,\\[10pt]
J\sigma_{1,N}(K_1)J^{-1}&=\maa{22pt}{0 & 0 & -q^{-1} \\ 0 & 1 & 0 \\ -q & 0 & 0}
\maa{22pt}{q & 0 & 0 \\ 0 & 1 & 0 \\ 0 & 0 & q^{-1}}\maa{22pt}{0 & 0 & -q^{-1} \\ 0 & 1 & 0 \\ -q & 0 & 0}
=\maa{22pt}{q^{-1} & 0 & 0 \\ 0 & 1 & 0 \\ 0 & 0 & q} \\[5pt]
&=\sigma_{1,-N}(K_1^{-1}) \;,\\[10pt]
J\sigma_{1,N}(E_1)J^{-1}&=[2]^{\frac{1}{2}}\maa{22pt}{0 & 0 & -q^{-1} \\ 0 & 1 & 0 \\ -q & 0 & 0}
\maa{22pt}{0 & 1 & 0 \\ 0 & 0 & 1 \\ 0 & 0 & 0}\maa{22pt}{0 & 0 & -q^{-1} \\ 0 & 1 & 0 \\ -q & 0 & 0}
\maa{22pt}{0 & 0 & 0 \\ -q & 0 & 0 \\ 0 & -q & 0} \\[5pt]
&=-q\sigma_{1,-N}(F_1) \;,\\[10pt]
J\sigma_{1,N}(F_1)J^{-1}&=[2]^{\frac{1}{2}}\maa{22pt}{0 & 0 & -q^{-1} \\ 0 & 1 & 0 \\ -q & 0 & 0}
\maa{22pt}{0 & 0 & 0 \\ 1 & 0 & 0 \\ 0 & 1 & 0}\maa{22pt}{0 & 0 & -q^{-1} \\ 0 & 1 & 0 \\ -q & 0 & 0}
=\maa{22pt}{0 & -q^{-1} & 0 \\ 0 & 0 & -q^{-1} \\ 0 & 0 & 0} \\[5pt]
&=-q^{-1}\sigma_{1,-N}(E_1) \;.
\end{align*}
Also, by a direct computation one checks that
\begin{align*}
J\mu_0(v,v') &=-\mu_0(Jv',Jv) \;,\qquad
J\mu_1(v,v') =-\mu_1(Jv',Jv) \;,\\
J\mu_2(v,v') &=\mu_3(Jv',Jv) \;,\qquad
J\mu_3(v,v') =\mu_2(Jv',Jv) \;,\qquad
J\mu_4(v,v') =\mu_4(Jv',Jv) \;.
\end{align*}
With these, \eqref{eq:claimA} is straightforwardly established. 
\end{proof}

{}From now on, we take that $c_4=c_0$ for the coefficients of Prop.~\ref{pr:const}.

The composition of the involution on $\Oq$ with the $^\star$ in the previous Lemma, yields a map $\omega\mapsto\omega^*$ sending forms into forms and extending the involution of the algebra $\Aq=\Omega^{0,0}(\CP^2_q)$. Indeed, with $h\in\Kq$, and $t:=S(h)^*$, from \eqref{eq:J}, we get
\begin{align*}
(\mL{h_{(1)}}\otimes\sigma^{i,j}(h_{(2)}))(\omega^*)_{i,j}
&=(-1)^i(\mL{h_{(1)}}\otimes\sigma^{i,j}(h_{(2)}))(*\otimes J)\omega_{j,i} \\
&=(-1)^i(*\otimes J)(\mL{S(h_{(1)}^*)}\otimes\sigma^{j,i}(S(h_{(2)})^*))\omega_{j,i} \\
&=(-1)^i(*\otimes J)(\mL{S^2(t_{(2)})}\otimes\sigma^{j,i}(t_{(1)}))\omega_{j,i} \;.
\end{align*}
Invariance of $\omega_{j,i}$ gives
$$
\mL{h}\omega_{j,i}=
(1\otimes\sigma^{j,i}(S^{-1}(h_{(3)})))(\mL{h_{(1)}}\otimes\sigma^{j,i}(h_{(2)}))\omega_{j,i}
=\sigma^{j,i}(S^{-1}(h))\omega_{j,i} \;,
$$
and in turn, using $\epsilon(S(h)^*)=\epsilon(h)$, 
\begin{multline*}
(\mL{h_{(1)}}\otimes\sigma^{i,j}(h_{(2)}))(\omega^*)_{i,j}  \\
=(-1)^i(*\otimes J)(1\otimes\sigma^{j,i}(t_{(1)}S(t_{(2)})))\omega_{j,i}=\epsilon(h)
(-1)^i(*\otimes J)\omega_{j,i} =\epsilon(h)(\omega^*)_{i,j} \;.
\end{multline*}
Thus, the involution maps invariant elements into invariant elements, i.e.~forms into forms.

As a consequence of \eqref{eq:claimA}, $(\Omega^{\bullet,\bullet}(\CP^2_q),\wprod\,,^*)$ is a graded $*$-algebra:
\begin{equation}\label{eq:sym}
(\omega\wprod\omega')^*=(-1)^{\mathrm{dg}(\omega)\mathrm{dg}(\omega')}\,\omega'^*\!\wprod\omega^* \;,
\qquad \forall\;\omega,\omega'\in\Omega^{\bullet,\bullet}(\CP^2_q)\;.
\end{equation}

\begin{lemma}
The algebra $(V^{\bullet,\bullet},\wprod)$ is generated in degree $1$, that is
any form of degree $\geq 1$  can be written as a sum of products of $1$-forms.
\end{lemma}

\begin{proof}
We need to show that the maps
$$
V^{0,1}\times V^{i,j}\to V^{i,j+1} \;,\qquad (v,w)\mapsto v\wprod w \;,
$$
and
$$
V^{1,0}\times V^{i,j}\to V^{i+1,j} \;,\qquad (v,w)\mapsto v\wprod w \;,
$$
are surjective for all $i,j$. We give the proof for the first map, the second
being analogous. \\
If $w$ is a scalar, the claim is clearly true.\\
If $(i,j)=(0,1)$ (resp.~$(2,1)$) the scalar
$$
(1,0)^t\wprod (0,1)^t=c_0[2]^{-\frac{1}{2}} \;,\qquad
\mathrm{resp.} \quad (1,0)^t\wprod (0,1)^t=c_3[2]^{-\frac{1}{2}} \;,
$$
is a basis $V^{0,2}$ (resp.~$V^{2,2}$), and the map is clearly surjective. \\
If $(i,j)=(1,0)$, the map $\wprod$ is   invertible; 
indeed,  $v\wprod w=\mathrm{diag}(c_1,c_1,c_1,c_2)U(v\otimes w)$ with
$U$ the unitary matrix in \eqref{eq:u}. \\
Finally,  if $(i,j)=(1,1)$ the vectors
$$
(1,0)^t\wprod (0,0,0,1)^t=-c_0c_2^{-1}[2]^{-1} (1,0)^t \;,\quad
(0,1)^t\wprod (0,0,0,1)^t=-c_0c_2^{-1}[2]^{-1} (0,1)^t \;,
$$
form a basis of $V^{1,2}$, and the maps surjective. This concludes the proof.
\end{proof}

\subsection{Hodge star and a closed integral}\label{se:hsci}~\\
Having an (associative, graded involutive) algebra of forms, the next steps consist 
in endowing it with  i) two derivations $\de$ and $\deb$ giving
a double complex, with $\dd:=\de+\deb$ the total differential;  ii) a closed integral; iii) an Hodge star operator. 

We postpone to the next section the explicit construction of the exterior differentials and start here with the integral. Recall from Remark~\ref{innprod} that an inner product on forms is given by composing the natural Hermitian structure on the module $\omega_{i,j}$ of forms with the restriction to $\Aq$ of the Haar state $\varphi:\Oq\to\C$; that is 
$$
\inner{\omega,\omega'}:=\sum\nolimits_{i,j}\varphi\bigl(
\omega_{i,j}^\dag \cdot \omega'_{i,j}\bigr) \;,
$$
where the sum is over all the homogeneous components $\omega_{i,j},\omega'_{i,j}\in\Omega^{i,j}(\CP^2_q)$ 
of $\omega$ and $\omega'$. Both the left canonical action $\az$ and the left action $\mathcal{L}$ of $\Uq$ on $\Oq$ are unitary for this inner product, as well as the action of $\Aq$ by left multiplication.

Now $\Omega^{2,2}(\CP^2_q)=\Aq$ is a free module of rank one, with basis a central element which we denote $1$ (the volume form). Indeed, we call $\mathtt{vol}$ the form with all components equal to zero but for the one in degree $4$, which is $1$. We think of this as the volume form and define an integral by 
\begin{equation}\label{int}
\nint\omega:=
\inner{\mathtt{vol},\omega}
=\varphi(\omega_{2,2}) \;,\qquad\forall\;\omega\in\Omega^{\bullet,\bullet}(\CP^2_q)\;.
\end{equation}
In particular, $\int\mkern-16mu-\;\mathtt{vol}=1$. If the differentials $\de$ and $\deb$ are given via  the (right) action of elements of $\Uq$ which are in the kernel of the counit $\epsilon$, the integral is automatically closed, i.e.
$$
\nint\deb\omega=\nint\de\omega=0 \;,
$$
a simple consequence of the invariance of the Haar state: $\varphi(a\za h)=\epsilon(h)\varphi(a)$.

Using the Hermitian structure on $\Omega^{\bullet,\bullet}$, given by
$$
(\omega,\omega'):=\sum\nolimits_{i,j}\omega_{i,j}^\dag \cdot \omega'_{i,j}
$$
the Hodge star operator is defined on \emph{real} forms via the usual requirement that
$\omega\wprod\omega'=(\ast_H\,\omega,\omega')\texttt{vol}$. 
This can be extended to complex forms both linearly (as e.g.~in 
\cite{Wel80}) or antilinearly (as e.g.~in \cite{Sch63}). In the context of solutions of the Yang-Mills equations, a curvature, being the square of a connection, is always real and it doesn't matter which extension we choose. We choose the former. Recalling
that the Hermitian structure $(\,,\,)$ is linear in the second entry and antilinear in the first, the Hodge star on complex forms is the linear operator
$\ast_H:\Omega^{i,j}(\CP^2_q)\to\Omega^{2-j,2-i}(\CP^2_q)$ defined by
\begin{equation}\label{eq:4.4}
\omega^*\!\wprod\omega'=(\ast_H\,\omega,\omega')\texttt{vol} \;.
\end{equation}
Applying the Haar states to both sides of previous equation we get the usual
relation
\begin{equation}\label{eq:4.3}
\nint\omega^*\!\wprod\omega'=\inner{\ast_H\,\omega,\omega'}\;.
\end{equation}
With a $*$-calculus, a closed integral and the graded Leibniz rule for the differential, the equality (\ref{eq:4.3}) implies that 
\begin{align*}
\inner{\ast_H\,\dd\omega,\omega'} &=\nint(\dd\omega)^*\!\wprod\omega'=-\nint\dd(\omega^*)\wprod\omega' \\
&=-\nint\dd(\omega^*\wprod\omega')+(-1)^{\mathrm{dg}(\omega)}\nint\omega^*\!\wprod\dd\omega' 
=0+(-1)^{\mathrm{dg}(\omega)}\inner{\ast_H\,\omega,\dd\omega'} \;,
\end{align*}
which becomes
\begin{equation}\label{eq:astast}
\dd^\dag\omega=\ast_H\,\dd\ast_H\omega \;,
\end{equation}
if $\ast_H^2 \omega =(-1)^{\mathrm{dg}(\omega)}\omega$. To obtain this last property, which is automatic when $q=1$,  one needs suitable constraints on the parameters $c_i$'s in Prop.~\ref{pr:const}.

\begin{prop}\label{lemma:4.3}
On any form $\omega$ one has \, $\ast_H^2 \omega =(-1)^{\mathrm{dg}(\omega)} \omega$\, if
$$
c_1=\pm q^{-\frac{1}{4}s}[2]^{-\frac{1}{4}}\sqrt{|c_0|} \;,\qquad
c_2=\pm q^{\frac{3}{4}s}[2]^{-\frac{1}{4}}\sqrt{|c_0|} \;,\qquad
c_3=\pm [2]^{\frac{1}{2}} \;,
$$
with arbitrary signs. For this choice of parameters, on (anti-)holomorphic $2$-forms
the Hodge star is the identity, and on $(1,1)$ forms it is the linear map
\begin{equation}\label{nothstar}
\ast_H:(w,w_4)\mapsto\mathrm{sign}(c_3)(-w,w_4) \;.
\end{equation}
\end{prop}

\begin{proof}
If $\ast_H\,\omega$ and $\omega'$ are homogeneous with different degree,
both sides of (\ref{eq:4.4}) are zero. It is then enough to
consider the case $\omega\in\Omega^{j,i}(\CP^2_q)$, $\omega'\in\Omega^{2-i,2-j}(\CP^2_q)$.
{}From the definition of the involution on forms,
for the possible values of the labels, one gets
\begin{align*}
(i,j)=(0,0),(0,2),(2,0),(2,2):\quad
     \omega^*\! &\wprod \omega'=  \omega^\dag \cdot \omega' \;,\\
(i,j)=(0,1),(1,0),(1,2),(2,1):\quad
     \omega^*\! &\wprod \omega'=c_3\mu_0(\omega^*,\omega')=
     (-1)^j[2]^{-\frac{1}{2}} c_3\, \omega^\dag \cdot \omega'  \;,\\
(i,j)=(1,1):\quad
     \omega^*\! &\wprod \omega'=
     \frac{c_3c_4}{[2]}
     \bigl(-q^{-\frac{1}{2}s}|c_1|^{-2}w^\dag,q^{\frac{3}{2}s}|c_2|^{-2}w_4^\dag\bigr)
     \omega' \;.
\end{align*}
Condition \eqref{eq:4.4} is satisfied if
$$
(\ast_H\,\omega)_{2-i,2-j}=
\begin{cases}
   \omega_{j,i} & \mathrm{if}\;(i,j)=(0,0),(0,2),(2,0),(2,2)\;,\\
  (-1)^j[2]^{-\frac{1}{2}}c_3\omega_{j,i}
 & \mathrm{if}\;(i,j)=(0,1),(1,0),(1,2),(2,1)\;, \\
  \frac{1}{[2]}c_3c_0 (-q^{-\frac{1}{2}s}|c_1|^{-2}w,q^{\frac{3}{2}s}|c_2|^{-2}w_4)
 & \mathrm{if}\;(i,j)=(1,1)\;.
\end{cases}
$$
The square of $\ast_H$ on $\omega$ is verified to be $(-1)^{\mathrm{dg}(\omega)}\omega$ if the $c_i$ are those given in the statement of the proposition.
With these, \eqref{nothstar} is immediately checked. 
\end{proof}

\noindent With the previous lemma, all parameters are fixed, but for some arbitrary signs 
and a global rescaling\footnote{ Notations are simplified by choosing $c_0=[2]^{\frac{1}{2}}$. 
In the notations of \cite{DDL08b} we would have $c_0=2[2]^{-\frac{1}{2}}$.} 
encoded in $c_0$. On the other hand, fixing the sign of $c_3$ corresponds to fixing an orientation, as flipping the the orientation results in exchanging selfdual with anti-selfdual forms. From now on we assume that $c_3<0$, so that by \eqref{nothstar} selfdual $(1,1)$-forms are of the type $(w,0)$, and anti-self dual ones are of the type $(0,w_4)$.

\begin{cor}
The Hodge star operator is an isometry.
\end{cor}
\begin{proof}
Notice that since $\varphi(a^*)=\overline{\varphi(a)}$, by graded involutivity of
the conjugation of forms we have
\begin{align*}
\inner{\ast_H\,\omega,\ast_H\,\omega'} &=\overline{\inner{\ast_H\,\omega',\ast_H\,\omega}}
=\left(\nint\omega'^*\!\wprod\ast_H\,\omega\right)^*
=\nint\bigl(\omega'^*\!\wprod\ast_H\,\omega\bigr)^* \\
&=(-1)^{\mathrm{dg}(\omega)}\nint (\ast_H\,\omega)^*\wprod\omega'
=(-1)^{\mathrm{dg}(\omega)}\inner{\ast_H^2\omega,\omega'}
=\inner{\omega,\omega'} \;,
\end{align*}
where we used the fact that $\omega$ and $\omega'$ have the same degree,
and $\mathrm{dg}(\omega)$ and $\mathrm{dg}(\omega)^2$ have the same parity.
\end{proof}

\begin{rem}
There is an equivalent definition of the Hodge star operator.
Firstly, one can define the ``exterior multiplication''
$\mathfrak{e}:\Omega^{\bullet,\bullet}(\CP^2_q)\to\mathrm{End}(\Omega^{\bullet,\bullet}(\CP^2_q))$
and the dual ``contraction'' $\mathfrak{i}:\Omega^{\bullet,\bullet}(\CP^2_q)\to\mathrm{End}(\Omega^{\bullet,\bullet}(\CP^2_q))$
by the formul{\ae}:
$$
\mathfrak{e}(\omega)\omega':=\omega\wedge\omega'\;,\qquad
\mathfrak{i}(\omega):=\mathfrak{e}(\omega)^\dag\;,\qquad
\forall\;\omega,\omega'\in\Omega^{\bullet,\bullet}(\CP^2_q)\;.
$$
Then, from \eqref{eq:4.3} and \eqref{int} one has 
$$
\inner{\ast_H\,\omega,\omega'} = \inner{\mathtt{vol},  \omega^*\!\wedge \omega'}=
\inner{\mathtt{vol},  \mathfrak{e}(\omega^*)\omega'}=\inner{\mathfrak{i}(\omega^*)\mathtt{vol},\omega'}
$$ 
which, from the non degeneracy of the scalar product, yields 
$$
\ast_H\,\omega=\mathfrak{i}(\omega^*)\mathtt{vol} 
$$
for any $\omega\in\Omega^{\bullet,\bullet}(\CP^2_q)$.
\end{rem}

\subsection{The exterior derivatives}\label{se:exde}~\\
We are left with the definition of the exterior derivative $\dd$. 
As we recall in App.~\ref{ap:cal},  in order to have a real differential calculus for 
$\Omega^\bullet(\CP^2_q) = \bigoplus_k \Omega^k(\CP^2_q)$ one needs a derivation, that is a map 
$\dd:\Aq\to\Omega^1(\CP^2_q)$ obeying the  Leibniz rule and such that 
$\Aq\left(\dd\Aq\right)=\Omega^1(\CP^2_q)$ 
and $\dd a=-(\dd a^*)^*$. Then the exterior derivative $\dd$ is extended uniquely
to forms of higher degree. If $\Omega^k(\CP^2_q):=\bigoplus_{i+j=k}\Omega^{i,j}(\CP^2_q)$, 
this is equivalent to write $\dd = \de + \deb$ with two derivations $\de:\Aq \to\Omega^{1,0}(\CP^2_q)$ and
$\deb:\Aq\to\Omega^{0,1}(\CP^2_q)$, such that $\deb\omega=-(\de\omega^*)^*$;  
again both $\de$ and $\deb$ are extended uniquely to forms of higher degree.

We write $X=\sum_iX_i\otimes e^i\in\Uq\otimes V^{1,0}$, with $e^1=(1,0)^t$ and $e^2=(0,1)^t$ the basis vectors of $V^{1,0}$. Then, we set  
$$
\de(\,\cdot\,):=\mL{X}\wprod(\,\cdot\,)=\textstyle{\sum_k}\mL{X_k}\otimes e^k\wprod(\,\cdot\,)
$$
and determine the conditions on the elements $X_i$'s that yield a $\de$ mapping $(i,j)$-forms, that is any 
$\omega\in\Omega^{i,j}(\CP^2_q) =\Oq\!\boxtimes_{\sigma^{i,j}}\!V^{i,j}$, to $(i+1,j)$-forms; namely 
we impose 
that $\de \omega \in \Omega^{i+1,j} (\CP^2_q)=\Oq\!\boxtimes_{\sigma^{i+1,j}}\!V^{i+1,j}$. 

Since the assignment $h\to\mL{h}$ is a representation, for all $h,x\in\Uq$ one gets that 
$\mL{h}\mL{x}=\mL{x\adj S^{-1}(h_{(2)})}\mL{h_{(1)}}$, 
with the right adjoint action given by
$$
x\adj h=S(h_{(1)})xh_{(2)} \;.
$$
In turn, left covariance of the wedge product yields   
\begin{align*}
\bigl\{\mL{h_{(1)}}\otimes\sigma^{i+1,j}(h_{(2)})\bigr\}\, \de (\,\cdot\,)
&=\textstyle{\sum_k}\,\mL{h_{(1)}}\mL{X_k}\otimes\sigma^{i+1,j}(h_{(2)})e^k\wprod(\,\cdot\,) \\
&=\textstyle{\sum_k}(\mL{X_k\adj S^{-1}(h_{(2)})}\otimes
\sigma^{1,0}(h_{(3)})e^k\wprod)(\mL{h_{(1)}}\otimes\sigma^{i,j}(h_{(4)})) (\,\cdot\,) \;,
\end{align*}
for any $h\in\Kq$. Then, for $X$ an element  which is invariant for the tensor product of the right adjoint action with $\sigma^{1,0}$, i.e.~for $X$ such that
$$
\textstyle{\sum_k}X_k\adj S^{-1}(h_{(1)})\otimes\sigma^{1,0}(h_{(2)})e^k=\epsilon(h)\textstyle{\sum_k}X_k\otimes e^k,
$$
we conclude that
\begin{equation}\label{pinco}
\bigl\{\mL{h_{(1)}}\otimes\sigma^{i+1,j}(h_{(2)})\bigr\} \circ \de
=\de \circ \bigl\{\mL{h_{(1)}}\otimes\sigma^{i,j}(h_{(2)})\bigr\} \;.
\end{equation}
This means that $\de$ maps invariant elements of $\Oq\otimes V^{i,j}$ into invariant elements of
$\Oq\otimes V^{i+1,j}$, that is forms to forms. Since both $\mL{h}$ and $\sigma^{i,j}$
are $*$-representation, taking the adjoint of equation \eqref{pinco} yields 
$$
\de^\dag\bigl\{\mL{h_{(1)}^*}\otimes\sigma^N_{k+1}(h_{(2)}^*)\bigr\}
=\bigl\{\mL{h_{(1)}^*}\otimes\sigma^N_k(h_{(2)}^*)\bigr\}\de^\dag \;,
$$
which means that $\de^\dag$ maps $(i,j)$-forms to $(i-1,j)$-forms. 

Similarly, by replacing $V^{1,0}$ 
with the representation $V^{0,1}$, an invariant element $Y$ of $\Uq\otimes V^{0,1}$ would give an operator $\deb=\mL{Y}\wprod(\,\cdot\,)$ mapping $(i,j)$-forms to $(i,j+1)$-forms and the adjoint 
$\deb^\dag$ mapping $(i,j)$-forms to $(i,j-1)$-forms.

On $\CP^2_q$ the only elements of $\Uq$ which act  on the right as derivations of $\Aq$ are the operators $E_2,F_2,[E_1,E_2]_q$ and $[F_2,F_1]_q$. These must be
applied in such a way to get an element of $\Omega^1(\CP^2_q)$.
\begin{prop}
One defines two exterior derivations $\de:\Aq\to\Omega^{1,0}(\CP^2_q)$ and $\deb:\Aq\to\Omega^{0,1}(\CP^2_q)$ by
\begin{equation}\label{eq:ddb}
\de a:=q^{-\frac{1}{2}}(\mL{q^{-1}K_2E_2}a,\mL{-K_1K_2[E_1,E_2]_q}a)^t \;, \qquad
\deb a:=(\mL{K_1K_2[F_2,F_1]_q}a,\mL{K_2F_2}a)^t \;,
\end{equation}
for all $a\in\Aq$. The map $\dd=\de+\deb:\Aq\to\Omega^1(\CP^2_q)$ satisfies the conditions 
$\Aq\dd\Aq=\Omega^1(\CP^2_q)$ and $\dd a^*=-(\dd a)^*$.
\end{prop}

\begin{proof}
The elements $X\in\Uq\otimes V^{1,0}$ and $Y\in\Uq\otimes V^{0,1}$ given by
\begin{equation}\label{eq:XandY}
X:=(q^{-1}K_2E_2,-K_1K_2[E_1,E_2]_q)^t\;,\qquad
Y:=(K_1K_2[F_2,F_1]_q,K_2F_2)^t \;,
\end{equation}
are invariant, hence by the above discussion the maps $\omega\mapsto\mL{X}\wprod\omega$ and $\omega\mapsto\mL{Y}\wprod\omega$
send $\Omega^{i,j}(\CP^2_q)$ into $\Omega^{i+1,j}(\CP^2_q)$, resp.~$\Omega^{i,j+1}(\CP^2_q)$. The wedge product on zero
forms is diagonal multiplication and on functions these maps become (proportional to) (\ref{eq:ddb}).

Since $S^{-1}$ is anticomultiplicative and $\mL{h}a=a\za S^{-1}(h)$,  
from the covariance of the right canonical action we get that for all $a,b\in\Oq$
\begin{align*}
\mL{X_1}(ab) &=(\mL{X_1}a)b+(a\za K_2^{-2})(\mL{X_1}b) \;,\\
\mL{X_2}(ab) &=(\mL{X_2}a)b+(a\za K_1^{-2}K_2^{-2})(\mL{X_2}b)+q^{-\frac{1}{2}}(q-q^{-1})
(a\za K_1^{-1}K_2^{-2}E_1)(\mL{X_1}b) \;,\\
\mL{Y_1}(ab) &=(\mL{Y_1}a)b+(a\za K_1^{-2}K_2^{-2})(\mL{Y_1}b)-q^{\frac{1}{2}}(q-q^{-1})
(a\za K_1^{-1}K_2^{-2}F_1)(\mL{Y_2}b) \;,\\
\mL{Y_2}(ab) &=(\mL{Y_2}a)b+(a\za K_2^{-2})(\mL{Y_2}b) \;.
\end{align*}
If $a,b$ are indeed in $\Aq$, right $\Kq$-invariance gives
$\mL{X}(ab)=(\mL{X}a)b+a(\mL{X}b)$ and $\mL{Y}(ab)=(\mL{Y}a)b+a(\mL{Y}b)$, so
$\de$ and $\deb$ are derivations on $\Aq$.

Reality is a simple check. {}From $(a^*\za h)^*=a\za S(h)^*$ it follows that
\begin{align*}
\deb a+(\de a^*)^* &=a\za\bigl(S^{-1}(Y_1)-q^{-1}X_2^*,S^{-1}(Y_2)+X_1^*\bigr) \\
&=a\za\bigl(K_1K_2[F_2,F_1]_q+q(K_1K_2)^{-1}[F_1,F_2]_q,(K_2-K_2^{-1})F_2\bigr) \\
&=a\za\bigl((q-q^{-1})F_1F_2,0\bigr) \\
&=0 \;.
\end{align*}

For the generators $p_{ij}=z_i^*z_j$ of $\Aq$ one computes that
$$
\de p_{ij}= - q^{-1} (u^3_i)^*\binom{u^2_j}{u^1_j} \;,\qquad
\deb p_{ij}= q^{-1} \binom{-q^{-\frac{1}{2}}(u^1_i)^*}{q^{\frac{1}{2}}(u^2_i)^*}u^3_j \;, 
$$
and shows that $\dd p_{ij}$ are a generating family for 
$$
\Omega^{1,0}(\CP^2_q)\oplus\Omega^{0,1}(\CP^2_q)=\mathfrak{M}(\sigma_{\frac{1}{2},-\frac{3}{2}})\oplus\mathfrak{M}(\sigma_{\frac{1}{2},\frac{3}{2}})
$$
as a left $\Aq$-module. Indeed, for any $\omega=(v_1,v_2)^t\oplus (w_1,w_2)^t\in\Omega^{1,0}(\CP^2_q)\oplus\Omega^{0,1}(\CP^2_q)$ the coefficients
\begin{align*}
a_{ij}(\omega)&:=-q^{1-2j}\bigl\{q^2v_1(u^2_j)^*+ v_2(u^1_j)^*\bigr\}u^3_i \;,\\
b_{ij}(\omega)&:=q^{5-2j}\bigl\{-q^{\frac{1}{2}}w_1(u^3_j)^*u^1_i+q^{-\frac{1}{2}}w_2(u^3_j)^*u^2_i\bigr\} \;,
\end{align*}
are right $\Kq$-invariant, i.e.~$a_{ij}(\omega), b_{ij}(\omega)\in\Aq$. Since 
$$
\sum\nolimits_jq^{2(a-j)}(u^a_j)^*u^b_j=\sum\nolimits_ju^a_j(u^b_j)^*=\delta_{a,b} \,, \qquad
\sum\nolimits_ju^3_j(u^3_j)^*=\sum\nolimits_jq^{6-2j}(u^3_j)^*u^3_j=1 ,
$$
the following algebraic identities 
\begin{align*}
 \sum\nolimits_{i,j}a_{ij}(\omega)\,\de p_{ij}&=(v_1,v_2)^t \;, &
 \sum\nolimits_{i,j}b_{ij}(\omega)\,\deb p_{ij}&=(w_1,w_2)^t \;, \\
 \sum\nolimits_{i,j}a_{ij}(\omega)\,\deb p_{ij} &=0 \;, &
 \sum\nolimits_{i,j}b_{ij}(\omega)\,\de p_{ij} &=0 \;
\end{align*}
hold. Assembling all together, for any $1$-form $\omega$ one finally gets
$$
\omega= \sum\nolimits_{i,j}\bigl\{a_{ij}(\omega)+b_{ij}(\omega)\bigr\}\dd p_{ij} \;,
$$
proving that the vectors $\dd p_{ij}$ are a generating family for $\Omega^1(\CP^2_q)$ as a left $\Aq$-module.
This concludes the proof.
\end{proof}

\begin{rem}
Using the properties $a \za K_1 = a$, $a \za K_2 = a$ and $a \za E_1 =a \za F_1  = 0$, 
for any element $a\in\Aq$,
simple manipulations in \eqref{eq:ddb} yields
$$
\de a:=-q^{-\frac{3}{2}} (a\za E_2,a\za E_2E_1)^t \;, \qquad
\deb a:=-(a\za F_2F_1, a\za F_2)^t \;.
$$
Also, modulo a proportionality constant the operator $\deb$  
coincides with the one of \cite{DDL08b}.
\end{rem}

\smallskip
We close this section with few remarks on invariant (anti-)selfdual $2$-forms that we shall use later on in the paper when dealing with monopole connections. Since $\Omega^{0,2}(\CP^2_q)\simeq\Sigma_{0,3}$, $\Omega^{2,0}(\CP^2_q)\simeq\Sigma_{0,-3}$, and $\Omega^{1,1}(\CP^2_q)\simeq\Sigma_{1,0}\oplus\Sigma_{0,0}$, by the harmonic decomposition
\begin{gather*}
\Sigma_{0,0} \simeq\bigoplus_{n\in\N}(n,n) \;,\qquad
\Sigma_{0,3} \simeq\bigoplus_{n\in\N}(n,n+3) \;,\qquad
\Sigma_{0,-3}\simeq\bigoplus_{n\in\N}(n+3,n) \;, \\
\Sigma_{1,0} \simeq\bigoplus_{n\in\N}(n+1,n+1)\oplus\bigoplus_{n\in\N}(n,n+3)\oplus (n+3,n) \;,
\end{gather*}
a $2$-form $\omega$ is invariant if{}f $\omega\in\Omega^{1,1}(\CP^2_q)$ and it has the form $\omega=(0,w_4)$ with $w_4\in\C$. By \eqref{nothstar} such a $2$-form is anti-selfdual, for the choice of orientation ($c_3<0$) made above. 


\section{Monopoles}\label{se:mono}
As mentioned in the introduction, by a monopole we mean a line bundle over $\CP^2_q$, that is to say a `rank 1' finitely generated projective module over the coordinate algebra $\A(\CP^2_q)$, endowed with a connection having anti-selfdual curvature. In this section we present some of these connections.

In Sect.~\ref{se:lb} we have described at length the 
isomorphism $\Sigma_{0,N}\to \Aq^{d_N}P_{-N}$ as left $\Aq$-modules
and the isomorphism $\Sigma_{0,N}\to P_N\Aq^{d_N}$ as right $\Aq$-modules.
Any of the two isomorphisms may be used to transport the Grassmannian connection on the
bimodule $\Sigma_{0,N}$. We use the second one because it is notationally simpler.

Recall that the isomorphism $\phi:\Sigma_{0,N}\to P_N\Aq^{d_N}$ is given $\phi(a):=\Psi_Na$, with inverse $\phi^{-1}:P_N\Aq^{d_N}\to\Sigma_{0,N}$, $\phi^{-1}(v)=\Psi^\dag\cdot v$; and  $\Psi_N$ is the column vector with components $\psi_{j,k,l}^N$ given in Prop.~\ref{lemma:linebundles}, $d_N:=\frac{1}{2}(|N|+1)(|N|+2)$ and $P_N:=\Psi_N\Psi_N^\dag$. 

\subsection{The anti-selfdual connections}~\\
The Grassmannian connection on the right 
$\Aq$-module $\E_N:=P_N\Aq^{d_N}$, with respect to the differential calculus of Sect.~\ref{se:diff} 
is the map
$$
\widetilde{\nabla}_{\!N} : \E_N \otimes_{\Aq}  \Omega^n(\CP^2_q)  \to \E_N \otimes_{\Aq}  \Omega^{n+1}(\CP^2_q)  \;,
\qquad \widetilde{\nabla}_{\!N} \omega:=P_N \dd\omega \;.
$$
For its curvature we get
\begin{align*}
\widetilde{\nabla}_{\!N}^2 \,\omega &=P_N\dd\bigl(P_N\dd(P_N\omega)\bigr)=
P_N\dd P_N\wprod \dd(P_N\omega) \\ &=
P_N\dd P_N\wprod \bigl\{P_N\dd\omega+\dd P_N\wprod\omega\bigr\}
=\bigl\{P_N\dd P_N\wprod \dd P_N \bigr\}\wprod\omega \;,
\end{align*}
were we used that $\omega=P_N\omega$ for any $\omega\in \E_N \otimes_{\Aq}  \Omega^n(\CP^2_q)$
and the identity
\begin{equation}\label{eq:ede}
e\dd e=(\dd e)(1-e)
\end{equation}
and so $e(\dd e)e=0$, both valid for any idempotent $e$ (and any differential calculus).

When transported to equivariant maps, the connection 
$$
\nabla_{\!N}:=\phi^{-1}\widetilde{\nabla}_{\!N} \,\phi : \Sigma_{0,N}\otimes_{\Aq}\Omega^n(\CP^2_q)\to \Sigma_{0,N}\otimes_{\Aq}\Omega^{n+1}(\CP^2_q) \;,
$$
is readily found to be given by 
\begin{equation}\label{conn-form}
\nabla_{\!N} \eta  =  \Psi_N^\dag \dd (\Psi_N \eta ) \;,
\end{equation}
and the curvature $\nabla_{\!N}^2=\phi^{-1} \widetilde{\nabla}_{\!N}^2 \,\phi$ becomes the operator of left wedge multiplication by the $2$-form, still denote $\nabla_{\!N}^2$, given by
\begin{equation}\label{eq:curv}
\nabla_{\!N}^2=\Psi_N^\dag(\dd P_N\wprod\dd P_N)\Psi_N \;.
\end{equation}

\begin{lemma}\label{incon}
The connection $\nabla_N$ is left $\Uq$-invariant.
\end{lemma}

\begin{proof}
{}From Prop.~\ref{lemma:psit} we deduce
$$
\Psi_N^\dag\otimes\Psi_N=\sum\nolimits_{\underline{i}}
t^{0,N}_{\uparrow,\underline{i}}
\otimes
(t^{0,N}_{\uparrow,\underline{i}})^* \;. 
$$
In turn, for any $h\in \Uq$,  (\ref{eq:azt}) yields
\begin{align*}
h_{(1)}\az\Psi_N^\dag\otimes h_{(2)}\az\Psi_N
&=\sum\nolimits_{\underline{i}}(h_{(1)}\az t^{0,N}_{\uparrow,\underline{i}})\otimes(S(h_{(2)})^*
\az t^{0,N}_{\uparrow,\underline{i}})^* \\
&=\sum_{\underline{i},\underline{j},\underline{k}}\rho^{0,N}_{\underline{j},\underline{i}}(h_{(1)})
\rho^{0,N}_{\underline{i},\underline{k}}(S(h_{(2)}))t^{0,N}_{\uparrow,\underline{j}}
\otimes(t^{0,N}_{\uparrow,\underline{k}})^* \\
&=\epsilon(h)\sum\nolimits_{\underline{j}}t^{0,N}_{\uparrow,\underline{j}}\otimes
(t^{0,N}_{\uparrow,\underline{j}})^* =\epsilon(h)\Psi^\dag_N\otimes\Psi_N \;.
\end{align*}
Thus, for any $h\in \Uq$, 
\begin{align*}
h\az(\nabla_{\!N}\eta)&=h\az(\Psi^\dag_N\dd\Psi_N\eta)=(h_{(1)}\az\Psi^\dag_N)\dd
(h_{(2)}\az\Psi_N)(h_{(3)}\az\eta)=\Psi_N^\dag\dd\Psi_N(h\az\eta) \\ &=
\nabla_N(h\az\eta) \;.
\end{align*}
This concludes the proof.
\end{proof}

The connection being invariant from Lemma~\ref{incon}, its curvature $\nabla^2_N$ is invariant as well. Then, from the discussion at the end of  Sect.~\ref{se:diff}, the two-form $\nabla^2_N\in\Omega^{1,1}(\CP^2_q)$ is necessarily of the type $\nabla_N^2=(0,w_N)$ for some $w_N\in\R$. Hence by \eqref{nothstar} it is anti-selfdual.

\begin{lemma}
The connection $\nabla_N$ is anti-selfdual, that is to say, its curvature is a (constant) anti-selfdual two-form:
$$
\ast_H\,\nabla_N^2= -\nabla_N^2 \;,
$$
\end{lemma}

For its use later on we compute the constant $w_N$ for $N\geq 0$; a similar computation being possible for $N\leq 0$ as well. By construction
\begin{equation}\label{eq:holZ}
\Psi_N\za K_2=q^{-\frac{N}{2}}\Psi_N \;    \qquad \textup{and} \qquad 
\Psi_N\za h=\epsilon(h)\Psi_N \;\;\; \forall\;h\in\Kq\;.
\end{equation}
Since $z_i\za F_2=0$, it also holds that 
\begin{equation}\label{eq:holA}
\Psi_N\za E_2=0\;,\qquad \Psi^\dag_N\za F_2=0 \;.
\end{equation}
Using the condition $\Psi^\dag_N\Psi_N=1$ and covariance of the action one deduces that
\begin{equation}\label{eq:holB}
(\Psi^\dag_N\za E_2)\Psi_N=0\;,\qquad
\Psi_N^\dag(\Psi_N\za F_2)=0 \;.
\end{equation}
These equations allow one to compute
\begin{align}\label{eq:vvprime}
q^{-\frac{N}{2}}\Psi^\dag_N\dd P_N &=-q^{-\frac{3}{2}}\dbinom{\Psi^\dag_N\za E_2}{\Psi^\dag_N\za E_2E_1}\in\Omega^{1,0}(\CP^2_q)
\;, \nn \\
~& \nn \\ q^{-\frac{N}{2}}(\dd P_N)\Psi_N &=-\dbinom{\Psi_N\za F_2F_1}{\Psi_N\za F_2}\in\Omega^{0,1}(\CP^2_q) \;.
\end{align}
By (\ref{eq:azt}) we have
\begin{align*}
q^{-N}\nabla_{\!N}^2 &=-c_2q^{-\frac{3}{2}s-3}[2]^{-\frac{1}{2}}
\sum\nolimits_{\underline{i}}
q(t^{0,N}_{\uparrow,\underline{i}}\za E_2)(t^{0,N}_{\uparrow,\underline{i}}\za E_2)^*
+q^{-1}(t^{0,N}_{\uparrow,\underline{i}}\za E_2E_1)(t^{0,N}_{\uparrow,\underline{i}}\za E_2E_1)^* \\
&=-c_2q^{-\frac{3}{2}s-3}[2]^{-\frac{1}{2}}
\sum\nolimits_{\underline{i},\underline{j},\underline{k}}
\bigl\{q
\rho^{0,N}_{\uparrow,\underline{j}}(E_2)\rho^{0,N}_{\underline{k},\uparrow}(F_2) \\ 
& \qquad\qquad\qquad\qquad\qquad\qquad\qquad\qquad
+q^{-1}\rho^{0,N}_{\uparrow,\underline{j}}(E_2E_1)\rho^{0,N}_{\underline{k},\uparrow}(F_1F_2)
\bigr\}t^{0,N}_{\underline{j},\underline{i}}(t^{0,N}_{\underline{k},\underline{i}})^* \\
&=-c_2q^{-\frac{3}{2}s-3}[2]^{-\frac{1}{2}}
\sum\nolimits_{\underline{j},\underline{k}}
\bigl\{q\rho^{0,N}_{\uparrow,\underline{j}}(E_2)\rho^{0,N}_{\underline{k},\uparrow}(F_2)
+q^{-1}\rho^{0,N}_{\uparrow,\underline{j}}(E_2E_1)\rho^{0,N}_{\underline{k},\uparrow}(F_1F_2)
\bigr\}\delta_{\underline{j},\underline{k}} \\
&=-c_2q^{-\frac{3}{2}s-3}[2]^{-\frac{1}{2}}
\rho^{0,N}_{\uparrow,\uparrow}(qE_2F_2+q^{-1}E_2E_1F_1F_2) \\
&=-c_2q^{-\frac{3}{2}s-3}[2]^{-\frac{1}{2}}
\inner{0,0,0|qE_2F_2+q^{-1}E_2E_1F_1F_2|0,0,0} \;.
\end{align*}
with $\ket{\uparrow}=\ket{0,0,0}$ the highest weight vector of the representation $\rho^{0,N}$.
Using the vanishing $E_1F_2\ket{0,0,0}=F_2E_1\ket{0,0,0}=0$, one gets that 
\begin{align*}
\inner{0,0,0|E_2E_1F_1F_2|0,0,0} &=\inner{0,0,0|E_2[E_1,F_1]F_2|0,0,0}=
\inner{0,0,0|E_2\smash[t]{\smash[b]{\tfrac{K_1^2-K_1^{-2}}{q-q^{-1}}}}F_2|0,0,0} \\ 
& = \inner{0,0,0|E_2F_2\smash[t]{\smash[b]{\tfrac{qK_1^2-q^{-1}K_1^{-2}}{q-q^{-1}}}}|0,0,0}=
\inner{0,0,0|E_2F_2|0,0,0} \;,
\end{align*}
giving
$ 
w_N  =
-q^{-\frac{3}{2}s-3}c_2[2]^{\frac{1}{2}}q^N\inner{0,0,0|E_2F_2|0,0,0}
$.\\
Since $F_2\ket{0,0,0}=\sqrt{[N]}\ket{0,1,\smash[t]{\smash[b]{\frac{1}{2}}}}$ we finally get 
$w_N = -q^{-\frac{3}{2}s-3}c_2[2]^{\frac{1}{2}}q^N[N]$.
Hence
\begin{equation}\label{curvN1}
\nabla_N^2 = q^{N-1} [N] \, \nabla_1^2 \;,
\end{equation}
where $\nabla_1^2=w_1\in\R$ is an irrelevant normalization constant.

\subsection{Gauged Laplacians}\label{se:gl}~\\
With monopoles connection on the modules $\Sigma_{0,N}$ we can study the corresponding gauged Laplacian operator acting on $\Sigma_{0,N}$. Such an operator describes `excitations moving on the quantum projective space' in the field of a magnetic monopole and in the limit $q \to 1$ provides a model of quantum Hall effect on the projective plane.

We know that on both the modules $\Sigma_{0,N}$ of sections and $\Omega^n(\CP^2_q)$ of forms 
there are $\Aq$-valued Hermitian structures. These can be combined to get a similar Hermitian structure on their tensor products $\Sigma_{0,N}\otimes_{\Aq}\Omega^n(\CP^2_q)$ which, when composed with the restriction to $\Aq$ of the Haar functional $\varphi  : \Aq\to\C$, gives a non-degenerate $\C$-valued  inner product.
Using the latter inner product, one generalized the equation \eqref{eq:astast} for forms to an analogous statement on connections, that is to say
\begin{equation}\label{eq:astconn}
\nabla_N^\dag\eta=\ast_H\,\nabla_{\!N}\!\,\ast_H\eta \;.
\end{equation}
Having this, as usual we define the gauged Laplacian acting on elements $\eta\in\Sigma_{0,N}$ by
$$
\Box_\nabla\eta=\nabla_N^\dag\nabla_N\eta=
\Psi_N^\dag\dd^\dag[\Psi_N\Psi_N^\dag\dd(\Psi_N\eta)] \;.
$$
where $\nabla_{\!N}  \eta  = \Psi_N^\dag\dd(\Psi_N \eta)$ 
is the connection in \eqref{conn-form}. Using the covariance of the action -- remember that 
$\mL{x}(ab)=(\mL{x_{(2)}}a)(\mL{x_{(1)}}b)$ due to the presence of the antipode -- , and equation \eqref{eq:holZ}, straightforward computations yield
$$
\de(\Psi_N\eta)=q^{\frac{N}{2}-\frac{1}{2}}(\mL{X}\Psi_N)\eta+q^{\frac{N}{2}-\frac{1}{2}}\Psi_N(\mL{X}\eta) \;,\qquad
\deb(\Psi_N\eta)=q^{\frac{N}{2}}(\mL{Y}\Psi_N)\eta+q^{\frac{N}{2}}\Psi_N(\mL{Y}\eta) \;.
$$
{}From \eqref{eq:holA} and \eqref{eq:holB}: $\mL{X}\Psi_N=0$
and $\Psi_N^\dag(\mL{Y}\Psi_N)=0$; hence
\begin{align*}
\Box_\nabla\eta &=q^{N/2}\Psi_N^\dag\dd^\dag\bigl(\Psi_N\bigl\{q^{-\frac{1}{2}}\mL{X}\oplus\mL{Y}\bigr\}\eta\bigr) \\
&=q^{N/2}\Psi_N^\dag\bigl\{q^{-\frac{1}{2}}\mL{X^\dag}\oplus\mL{Y^\dag}\bigr\}
\bigl(\Psi_N\bigl\{q^{-\frac{1}{2}}\mL{X}\oplus\mL{Y}\bigr\}\eta\bigr) \;.
\end{align*}
Again, \eqref{eq:holA} and \eqref{eq:holB} yield $\mL{Y^\dag}\Psi_N=0$
and $\Psi_N^\dag(\mL{X^\dag}\Psi_N)=0$; hence using the coproduct, 
$$
\Box_\nabla\eta=q^N\mL{q^{-1}X^\dag X+Y^\dag Y} \,\eta = \eta\za S^{-1}(q^{-1}X^\dag X+Y^\dag Y) \, q^N \;.
$$
Now,
\begin{multline*}
S^{-1}(q^{-1}X^\dag X+Y^\dag Y)=K_2^{-2}(q^{-1}E_2F_2+q^{-2}F_2E_2) \\ +K_1^{-2}K_2^{-2}
\bigl(q[E_2,E_1]_q[F_1,F_2]_q+q^{-2}[F_1,F_2]_q[E_2,E_1]_q \bigr) \;,
\end{multline*}
which, with the commutation rule
\begin{align*}
\bigl[[F_1, F_2]_q, [E_2, E_1]_q\bigr] &=
\bigl[[F_1, [E_2, E_1]_q], F_2\bigr]_q+\bigl[F_1, [F_2, [E_2, E_1]_q]\bigr]_q \\
&=\bigl[[E_2,[F_1,E_1]]_q, F_2\bigr]_q+\bigl[F_1, [[F_2, E_2], E_1]_q  \bigr]_q \\
&=-\bigl[[E_2,\tfrac{K_1^2-K_1^{-2}}{q-q^{-1}}]_q,F_2\bigr]_q-\bigl[F_1,[\tfrac{K_2^2-K_2^{-2}}{q-q^{-1}},E_1]_q\bigr]_q \\
&=-[K_1^2E_2,F_2]_q+[F_1,E_1K_2^{-2}]_q \\
&=K_1^2(q^{-2}F_2E_2-E_2F_2)+(F_1E_1-q^{-2}E_1F_1)K_2^{-2} \;,
\end{align*}
can be rewritten as
\begin{multline*}
S^{-1}(q^{-1}X^\dag X+Y^\dag Y) 
=K_2^{-2}\bigl([2]\tfrac{K_2^2-K_2^{-2}}{q-q^{-1}}+(q+q^{-2})F_2E_2\bigr) \\ +K_1^{-2}K_2^{-2}
(q+q^{-2})[F_1,F_2]_q[E_2,E_1]_q 
-qK_1^{-2}K_2^{-4}(F_1E_1-q^{-2}E_1F_1) \;.
\end{multline*}
By construction the action of $\mathcal{U}_q(su(2))$ is trivial on $\eta$, that is to say 
$\eta\za E_1=\eta\za F_1=0$ and $\eta\za K_1 = \eta$, while
$\eta\za K_2=q^{N/2}\eta$. Then, an intermediate result states that 
$$
\Box_\nabla\eta = \eta\za \left( \tfrac{K_2^2-K_2^{-2}}{q-q^{-1}} \, [2] + (q+q^{-2}) 
\bigl( F_2E_2 -q^{-1}  F_2 F_1 [E_2,E_1]_q \bigr) \right)
$$
Using the commutator $[F_2,E_1]=0$ and again $\eta\za E_1=0$, a 
further simplification comes from the computation
\begin{equation}\label{comp}
\eta\za F_2F_1E_1E_2=\eta\za F_2[F_1,E_1]E_2=-\eta\za F_2\tfrac{K_1^2-K_1^{-2}}{q-q^{-1}}E_2=\eta\za F_2E_2 \;,
\end{equation}
leading to: 
\begin{equation}\label{partlap}
\Box_\nabla\eta=\eta\za \Bigl( [2] [N]+(1+q^{-3}) \bigl([2] F_2E_2 - F_2 F_1 E_2 E_1 \bigr) \Bigr)  \;.
\end{equation}

Next, we relate the gauged Laplacian operator to the Casimir operator in 
\eqref{eq:Cq}. Now, when acting from the right on $\Sigma_{0,N}$, a straightforward  
computation leads to
\begin{multline*}
\eta\za\mathcal{C}_q=
\eta\za\Bigl\{
[\tfrac{1}{3}N]^2+[\tfrac{1}{3}N+1]^2+[\tfrac{2}{3}N+1]^2+ (q^{\frac{1}{3}N + 1}+q^{-\frac{1}{3}N -1})F_2E_2
\\ + 
F_2 F_1 \left( q^{-\frac{1}{3}N+1} [E_1,E_2]_q -  q^{\frac{N}{3}} [E_2,E_1]_q \right)
\Bigr\} \;,
\end{multline*}
which, using again \eqref{comp} reduces to
\begin{multline}\label{partcas}
\eta\za\mathcal{C}_q=
\eta\za\bigl\{
[\tfrac{1}{3}N]^2+[\tfrac{1}{3}N+1]^2+[\tfrac{2}{3}N+1]^2+ (q^{\frac{1}{3}N}+q^{-\frac{1}{3}N}) 
\bigl( [2] F_2E_2 - F_2 F_1 E_2 E_1 \bigr)
\bigr\} \;.
\end{multline}

The generator $L=K_1K_2^2$ of the `structure algebra' $\U_q(\mathfrak{u}(1))$ act on $\Sigma_{0,N}$ as 
$\eta \za L =q^{N/2}\eta$. Then, a comparison of \eqref{partlap} and \eqref{partcas} yields the following
proposition.
\begin{prop}
The gauged Laplacian is related to the Casimir operator by
\begin{multline*}
q^{\frac{3}{2}} \, \frac{ L^{\frac{1}{3}} +L^{-\frac{1}{3}} } {q^{\frac{3}{2}}+q^{-\frac{3}{2}}}
\left( \Box_\nabla + [2]\, \frac{ L - L^{- 1} } {q - q^{-1}} \right) \\ =
\mathcal{C}_q - ({q - q^{-1}})^{-1} \left( ( L^{\frac{1}{3}} - L^{-\frac{1}{3}} )^2 
 +  ( q L^{\frac{1}{3}} - q^{-1} L^{-\frac{1}{3}} )^2 
 + 
( q L^{\frac{2}{3}} - q^{-1} L^{-\frac{2}{3}} )^2  \right)
 \;.
\end{multline*}
or
$$
q^{\frac{3}{2}}\, \frac{q^{\frac{N}{3}}+q^{-\frac{N}{3}}}{q^{\frac{3}{2}}+q^{-\frac{3}{2}}} 
\bigl( \Box_\nabla - [2][N] \bigr) = 
\mathcal{C}_q-[\tfrac{1}{3}N]^2-[\tfrac{1}{3}N+1]^2-[\tfrac{2}{3}N+1]^2 \;.
$$
\end{prop}

\bigskip

The diagonalization of the gauged Laplacian is made simple by the observation that 
for the Casimir left and right action is the same: 
$\mathcal{C}_q\az a=a\za \mathcal{C}_q$ for all $a\in\Oq$, and by the fact that 
with respect to the left action of $\Uq$ there are decompositions:
$$
\Sigma_{0,N}\simeq\bigoplus_{n\in\N}\rho^{(n,n+N)}\quad\mathrm{if}\;N\geq 0\;,\qquad
\Sigma_{0,N}\simeq\bigoplus_{n\in\N}\rho^{(n-N,n)}\quad\mathrm{if}\;N\leq 0\;.
$$
\begin{prop}
The eigenvalues of the gauged Laplacian  $\Box_\nabla$ are given by 
\begin{align*}
\lambda_{n,N}&=
(1+q^{-3})[n][n+N+2]+[2][N] &&\mathrm{if}\;N\geq 0\;,\\
\lambda_{n,N}&=
(1+q^{-3})[n+2][n-N]+[2][N] &&\mathrm{if}\;N\leq 0\;,
\end{align*}
with $n\in\N$. 
\end{prop}

\begin{proof}
Using the above harmonic decomposition of $\Sigma_{0,N}$ 
and the spectrum (\ref{eq:SpCq}) of the Casimir operator, one gets
$$
\lambda_{n,N}:=q^{-\frac{3}{2}}\frac{q^{\frac{3}{2}}+q^{-\frac{3}{2}}}{q^{\frac{N}{3}}+q^{-\frac{N}{3}}}
\big(
[n+\tfrac{1}{3}|N|+1]^2+[n+\tfrac{2}{3}|N|+1]^2
-[\tfrac{1}{3}N+1]^2-[\tfrac{2}{3}N+1]^2\big)
+[2][N] \;.
$$
with $n\in\N$ and for any $N\in\Z$. This expression can be simplified using the identity $[a+b]^2-[b]^2=[a][a+2b]$.
For $N\geq 0$ this becomes
$$
\lambda_{n,N}=q^{-\frac{3}{2}}\frac{q^{\frac{3}{2}}+q^{-\frac{3}{2}}}{q^{\frac{N}{3}}+q^{-\frac{N}{3}}}
[n]\big([n+\tfrac{2}{3}N+2]+[n+\tfrac{4}{3}N+2]\big)+[2][N] \;,
$$
that with a further simplification is the claimed  expression in the proposition. One
proceeds similarly for the case $N\leq 0$.
\end{proof}
It is worth stressing that the spectrum of $\Box_\nabla$  is not  symmetric under
the exchange $N \leftrightarrow -N$, not even when exchanging in addition the parameter as $q \leftrightarrow q^{-1}$. 
A similar phenomenon was already observed in \cite{LRZ07} for gauged Laplacians on the 
standard Podle\'s sphere; the latter can be though of as the quantum projective line $\CP^1_q$.

\section{Quantum characteristic classes}\label{se:qinv}

Classically, topological invariants are computed by integrating powers of the curvature of a connection, the result being independent of the particular connection. 
On the other hand, in order to integrate the curvature of a connection on the quantum projective space 
$\CP^2_q$ one needs `twisted integrals'; the results are not integers any longer but rather $q$-analogues,
as we shall see in Sect.~\ref{se:qtinvcp2}. We start with some general result on equivariant K-theory and K-homology and corresponding Chern-Connes characters.

\subsection{Equivariant K-theory and K-homology}\label{se:eqkt}~\\
Classically, the equivariant topological $K_0$-group of a manifold is the Grothendieck group of the abelian monoid whose elements are equivalence classes of equivariant vector bundles. It has an algebraic version that can be generalized to noncommutative algebras. One has a bialgebra $\U$ and a left $\U$-module algebra $\A$. Equivariant vector bundles are replaced by one sided (left, say)  $\A\rtimes\U$-modules that are finitely generated and projective as left $\A$-modules; these will be simply called  `equivariant projective modules'. Any such a module is given by a pair $(e,\sigma)$, where $e$ is an $N\times N$ idempotent with entries in $\A$, and \mbox{$\sigma:\U\to\mathrm{Mat}_N(\C)$} is a representation and the following compatibility requirement is satisfied (see e.g.~\cite[Sect.~2]{DDL08}), 
\begin{equation}\label{eq:cov}
\left( h_{(1)}\az e \right)  \sigma(h_{(2)})^t=\sigma(h)^t e \;, \qquad \textup{for} \quad h\in\U, 
\end{equation}
with `$\phantom{|}^t$' denoting transposition.

The corresponding module $\E=\A^Ne$ is made of elements  $v=(v_1,\ldots,v_N)\in\A^N$ 
in the range of the idempotent, $ve=v$, with left-module structure given by
$$
(a.v)_i:=av_i\;,\qquad
(h.v)_i:=\sum\nolimits_{j=1}^N(h_{(1)}\az v_j)\sigma_{ij}(h_{(2)})\;, 
\qquad \textup{for}   \quad a\in\A  \quad \textup{and} \quad  h\in\U \,.
$$
 
An equivalence between any  two equivariant modules is simply an invertible left $\A\rtimes\U$-module map between them; $V^{\U}(\A)$ will denote the abelian monoid whose elements are equivalence classes of equivariant projective left modules with operation the direct sum, as usual. The equivariant K-theory group $K_0^{\U}(\A)$ is the Grothendieck group of the abelian monoid $V^{\U}(\A)$. The equivalence of equivariant projective modules can be rephrased in terms of idempotents. What follows is a direct extension of well known results \cite{Bla98}. We give the proof for completeness.

\begin{lemma}\label{lemma:KU}
Two equivariant projective modules $\E=\A^Ne$ and $\E'=\A^{N'}e'$ are equivalent if{}f $e=uv$ and $e'=vu$ for some
$u\in\mathrm{Mat}_{N\times N'}(\A)$ and $v\in\mathrm{Mat}_{N'\times N}(\A)$ satisfying the equivariance conditions
\begin{equation}\label{eq:KU}
(h_{(1)}\az u)\sigma'(h_{(2)})^t=\sigma(h)^tu \;,\qquad
(h_{(1)}\az v)\sigma(h_{(2)})^t=\sigma'(h)^tv \;.
\end{equation}
\end{lemma}

\begin{proof}
Let $w_{\bullet}=\sum_iw_ie_{i\,\bullet}$ be the generic element of $\A^Ne$. If $\pi:\A^Ne\to\A^{N'}e'$ is a left $\A$-module map, then $\pi(w_{\bullet})=\sum_iw_i\pi(e_{i\,\bullet})$, so the map is uniquely determined by the its value on rows of the idempotent $e$, and similarly for $\pi^{-1}$. We call $u\in\mathrm{Mat}_{N\times N'}(\A)$ (resp.~$v\in\mathrm{Mat}_{N'\times N}(\A)$) the matrix with entries $u_{ij}:=\pi(e_{i\,\bullet})_j$ (resp.~$v_{ij}:=\pi^{-1}(e'_{i\,\bullet})_j$). Since $\pi$ maps into the range of $e'$ (resp.~$\pi^{-1}$ maps into the range of $e$), we have the conditions
$$
\sum\nolimits_j\pi(e_{i\,\bullet})_je'_{jk}=\pi(e_{i\,\bullet})_k \;,\qquad \textup{and} \qquad 
\sum\nolimits_j\pi^{-1}(e'_{i\,\bullet})_je_{jk}=\pi^{-1}(e'_{i\,\bullet})_k \;,
$$
that in term of row vectors become 
$$
\sum\nolimits_ju_{ij}e'_{j\,\bullet}=u_{i\,\bullet} \;,\qquad \textup{and} \qquad
\sum\nolimits_jv_{ij}e_{j\,\bullet}=v_{i\,\bullet} \;.
$$
Next, we apply $\pi^{-1}$ to the equation on the left, $\pi$ to the one on the right, and use $\A$-linearity.
Since $u_{ij}$ is the $j$-th component of image, via $\pi$, of the $i$-th row of $e$, and $\pi^{-1}$
is the inverse map, we have $\pi(u_{i\,\bullet})=e_{i\,\bullet}$ and similarly for $\pi^{-1}$. Thus:
$$
\sum\nolimits_ju_{ij}v_{jk}=e_{ik} \;,\qquad
\sum\nolimits_jv_{ij}u_{jk}=e'_{ik} \;.
$$
{}From the definition of $u$, the equivariance properties of $e$, and the fact that $\pi$ is an $\U$-module map, we get
$$
h.u_{i\,\bullet}=h_{(1)}\az u_{i\,\bullet}\,\sigma'(h_{(2)})^t=\pi(h.e_{i\,\bullet})
 =\pi(h_{(1)}\az e_{i\,\bullet}\,\sigma(h_{(2)})^t) =\pi\bigl((\sigma(h)^te)_{i\,\bullet}\bigr)=(\sigma(h)^tu)_{i\,\bullet} \;,
$$
which is the first equivariance condition in \eqref{eq:KU}. Similarly one proves the equivariance  condition for $v$. With this, the `only if' part is proved. 

Next, assume that $e=uv$ and $e'=vu$ for some $u$ and $v$ satisfying \eqref{eq:KU} above.
Then,
$$
(h_{(1)}\az e)\sigma(h_{(2)})^t=
(h_{(1)}\az u)(h_{(2)}\az v)\sigma(h_{(3)})^t=
(h_{(1)}\az u)\sigma'(h_{(2)})^tv=
\sigma(h)^tuv \;,
$$
which means $e$ that satisfies (\ref{eq:cov}). Similarly for $e'$.

We define $\pi:\A^Ne\to\A^{N'}e'$ and $\pi^{-1}:\A^{N'}e'\to\A^Ne$ via the formul{\ae} $\pi(w):=wu$ and $\pi^{-1}(w):=wv$; we need to show that (i) the maps are well defined, (ii) they are one the inverse of the other, (iii) they are left $\A\rtimes\U$-module maps.
Point (iii) is a consequence of the fact that left and right multiplication commute, and of the equivariance conditions for $u$ and $v$.
Point (ii) follows from the identity $\pi^{-1}\pi(w)=we$ (resp.~$\pi\pi^{-1}(w)=we'$), and the fact that $w$ is in the range of $e$
(resp.~$e'$). Finally if $w\in\A^Ne$, we have
$$
\pi(w)e'=w(uvu)=(we)u=wu=\pi(w) \;,
$$
i.e.~$\pi(w)\in\A^{N'}e'$, and similarly for $\pi^{-1}$: the maps are well defined.
\end{proof}

There is a natural map from equivariant K-theory to equivariant cyclic homology given for instance in  \cite{NT03}. We adapt that construction to our situation. One start with $\mathrm{Hom}_{\C}(\U,\A^{n+1})$, the collection of $\C$-linear maps from $\U$ to $\A^{n+1}$, 
and defines operations $b_{n,i}:\mathrm{Hom}_{\C}(\U,\A^{n+1})\to \mathrm{Hom}_{\C}(\U,\A^n)$, for $i=0,\ldots,n$, 
\begin{align}\label{bhom}
b_{n,i}(a_0\otimes a_1\otimes\ldots\otimes a_n)(h)&:=
(a_0\otimes\ldots\otimes a_ia_{i+1}\otimes\ldots\otimes a_n)(h) \;,\qquad\mathrm{if}\;i\neq n\;, \nn \\
b_{n,n}(a_0\otimes a_1\otimes\ldots\otimes a_n)(h)&:=\bigl((h_{(1)}\az a_n)a_0\otimes a_1\otimes\ldots\otimes a_{n-1}\bigr)(h_{(2)}) \;,
\end{align}
and an operation 
$\lambda_n:\mathrm{Hom}_{\C}(\U,\A^{n+1})\to \mathrm{Hom}_{\C}(\U,\A^{n+1})$,
\begin{equation}
\lambda_n(a_0\otimes a_1\otimes\ldots\otimes a_n)(h) :=(-1)^n\bigl((h_{(1)}\az a_n)\otimes a_0\otimes a_1\otimes\ldots\otimes a_{n-1}\bigr)(h_{(2)}) \;.
\end{equation}
They are the face operators and the cyclic operator of equivariant cyclic homology as we 
are going to show.
The maps $b_{n,i}$ make up a presimplicial module -- one checks that $b_{n-1,i}b_{n,j}=b_{n-1,j-1}b_{n,i}$ for all $0\leq i<j\leq n$ --, so that 
$$
b_n:=\sum\nolimits_{i=0}^n(-1)^ib_{n,i}
$$
is a boundary operator \cite{Lod97}.
The Hopf algebra $\U$ acts on $\A^{n+1}$ via the rule
$$
h\aaz (a_0\otimes a_1\otimes\ldots\otimes a_n):=h_{(1)}\az a_0\otimes h_{(2)}\az a_1\otimes\ldots\otimes h_{(n+1)}\az a_n \;,
$$
and $C^{\U}_n(\A)$ will denote the collection of elements $\omega\in\mathrm{Hom}_{\C}(\U,\A^{n+1})$ which are `equivariant', meaning that
$$
(h_{(1)}\aaz\omega)(xh_{(2)})=\omega(hx) \;,
$$
for all $h,x\in\U$. Next, one establishes that the operators $b_{n,i}$ commute with the action of $\U$ (a not completely trivial task for the last one $b_{n,n}$), and it makes sense to consider the complex of equivariant maps. 
The cyclic operator $\lambda_n$ commutes with the action of $\U$ thus it descends to an operator on $C^{\U}_n(\A)$ as well.
Finally, with 
$$
b'_n:=\sum\nolimits_{i=0}^{n-1}(-1)^nb_{n,i} \,,
$$  
it holds that $b_n(1-\lambda_n)=(1-\lambda_{n-1})b'_n$, which says that the boundary operator $b_n$ maps $C^{\U}_n(\A)/\mathrm{Im}(1-\lambda_n)$ into $C^{\U}_{n-1}(\A)/\mathrm{Im}(1-\lambda_{n-1})$. 
The homology of this last complex is called the `$\U$-equivariant cyclic homology' of $\A$ with corresponding homology groups usually denoted $H\!C^{\U}_n(\A)$.

Next, for $\sigma:\U\to\mathrm{Mat}_N(\C)$ a representation as above, consider the set
$$
\mathrm{Mat}_N^\sigma(\A):=\big\{T\in\mathrm{Mat}_N(\A)\; \big|\; 
\left(h_{(1)}\az T \right) \sigma(h_{(2)})^t=\sigma(h)^t \,T\;, \;\; \forall\;h\in\U \,\big\} \;.
$$
This is a subalgebra of $\mathrm{Mat}_N(\A)$; indeed given any two of its elements $T_1,T_2$ one has: 
\begin{multline*}
\left(h_{(1)}\az (T_1T_2)\right) \sigma(h_{(2)})^t=
\left((h_{(1)}\az T_1)(h_{(2)}\az T_2) \right) \sigma(h_{(3)})^t \\ =
(h_{(1)}\az T_1) \, \sigma(h_{(2)})^t \, T_2=\sigma(h)^t \, T_1T_2\;.
\end{multline*}
Moreover, $\sigma$-equivariant $N\times N$ idempotents as in \eqref{eq:cov} are elements of $\mathrm{Mat}^\sigma_N(\A)$.
Due to the definition of $\mathrm{Mat}^\sigma_N(\A)$ there exists a map $\tr_\sigma:\mathrm{Mat}^\sigma_N(\A)^{n+1}\to C^{\U}_n(\A)$ given by
\begin{align*}
\tr_\sigma(T_0\otimes T_1\otimes\ldots\otimes T_n)(h)
&:=\tr\bigl(T_0\dotimes T_1\dotimes\ldots\dotimes T_n \, \sigma(h)^t\bigr) \\
&=\sum_{i_0,i_1,\ldots,i_{n+1}}(T_0)_{i_0i_1}\otimes (T_1)_{i_1i_2}\otimes\ldots\otimes
(T_n)_{i_ni_{n+1}} \, \sigma(h)_{i_0i_{n+1}} \;,
\end{align*}
where $\dotimes$ denotes composition of the tensor product over $\C$ with matrix multiplication.
Also,  
\begin{align*}
(-1)^n\lambda_n\tr_\sigma(T_0\otimes T_1\otimes\ldots\otimes T_n)(h)
&=\tr\bigl(h_{(1)}\az T_n\,\sigma(h)^t\dotimes T_0\dotimes T_1\dotimes\ldots\dotimes T_{n-1}\bigr) \\
&=\tr\bigl(\sigma(h)^t T_n\dotimes T_0\dotimes T_1\dotimes\ldots\dotimes T_{n-1}\bigr) \\
&=\tr_\sigma(T_n\otimes T_0\otimes\ldots\otimes T_{n-1})(h) \;,
\end{align*}
which amounts to say that $\tr_\sigma$ transform the ordinary cyclic operator for the algebra $\mathrm{Mat}^\sigma_N(\A)$
into the `$\U$-equivariant' cyclic operator for $\A$. Since $b_{n,n}=b_{n,0}\lambda_n$, the map $\tr_\sigma$ is a morphism of differential complexes, mapping the complex of the cyclic homology of
$\mathrm{Mat}^\sigma_N(\A)$ to the complex of the $\U$-equivariant cyclic homology of $\A$. This construction is completely analogous to the `non-equivariant' case, cf.~\cite[Cor.~1.2.3]{Lod97}.

At this point, one can repeat verbatim the proof of Thm~8.3.2 in \cite{Lod97}, replacing the ring
$R:=\mathrm{Mat}_N(\A)$ there, with $\mathrm{Mat}^\sigma_N(\A)$ (which is still a matrix ring) and replacing the generalized trace map there,  
with $\tr_\sigma$, to prove the following theorem.
\begin{thm}\label{thm}
A map $\mathrm{ch}^n:K_0^{\U}(\A)\to H\!C^{\U}_n(\A)$ is defined by
$$
\mathrm{ch}^n(e,\sigma):=\tr_\sigma(e^{\otimes n+1}) \;.
$$
\end{thm}
We give in App.~\ref{app:alt} an alternative, more explicit proof for
the cases $n\leq 4$. We stress that  what we denote here $H\!C_n$ and call
cyclic homology is Connes' first version of cyclic homology, i.e.~the homology of Connes' complex
denoted $H^\lambda_n$ in \cite{Lod97}.

\smallskip

Modulo a normalization, the cycle $\mathrm{ch}^n(e,\sigma)(1)$ is the usual Chern-Connes character in cyclic homology (and in fact, no $\sigma$'s in the formul{\ae}). On the other hand, $\mathrm{ch}^n(e,\sigma)((K_1 K_2)^{-4})$ is what we are about to use for 
$\CP^2_q$. In general, one fixes a group-like element $K\in\U$ calling $\eta$
the corresponding automorphism of $\A$, $\eta(a):=K\az a$ for all $a\in\A$. 
Then, one pairs $H\!C^{\U}_\bullet(\A)$ with the
Hochschild cohomology of $\A$ with coefficients in ${_\eta}\A$; the latter is $\A$ itself with bimodule structure $ a' (a) a'' = \eta(a') a a''$. Indeed, the pairing 
$\inner{\,,\,}:\mathrm{Hom}_{\C}(\A^{n+1},\C)\times\mathrm{Hom}_{\C}(\U,\A^{n+1})\to\C$ defined by
\begin{equation}\label{twdprng}
\inner{\tau,\omega}:=\tau\bigl(\omega(K)\bigr) \;,
\end{equation}
when used to compute the dual $b^*_{n,i}:\mathrm{Hom}_{\C}(\A^n,\C)\to\mathrm{Hom}_{\C}(\A^{n+1},\C)$
of the face operators introduced in \eqref{bhom}, yields the formul{\ae}:
\begin{align*}
b_{n,i}^*\tau(a_0,a_1,\ldots,a_n)&=\tau(a_0,\ldots,a_ia_{i+1},\ldots,a_n) \;,\qquad\mathrm{if}\;i\neq n\;,\\
b_{n,n}^*\tau(a_0,a_1,\ldots,a_n)&:=\tau\bigl(\eta(a_n)a_0,a_1,\ldots,a_{n-1}\bigr) \;.
\end{align*}
These are just the face operators of the Hochschild cohomology $H^\bullet(\A,{_\eta}\A)$ of $\A$ with coefficients in 
${_\eta}\A$ (cf.~\cite{Lod97}). Thus, the pairing in \eqref{twdprng} descends to a pairing 
$$
H^n(\A,{_\eta}\A)\times H\!C^{\U}_n(\A)\to\C \,. 
$$

\subsection{The example of $\CP^2_q$}\label{se:qtinvcp2}~\\
As mentioned, for $\CP^2_q$ we take $K=(K_1K_2)^{-4}$ the element implementing the square of the antipode (cf. \eqref{saut}). Now, the Haar state of $\Oq$, satisfies (cf.~Eq.~(11.26) and Eq.~(11.36) in \cite{KS97})
$$
\varphi(ab)=\varphi\bigl( (K\az b\za K ) a\bigr) \;,\qquad \textup{for} \quad a,b\in\Oq\;,
$$
which for $a,b\in\Aq$ results in 
\begin{equation}\label{eq:modprop}
\varphi(ab)=\varphi\bigl( (K\az b) a\bigr)=\varphi\bigl( \eta(b) a\bigr) \;.
\end{equation}
This just means that the restriction of the Haar state of $\Oq$ to $\Aq$ is the representative of a class in $H^0(\Aq,{_\eta}\Aq)$.
An additional zero cocycle is given by the restriction of the counit $\epsilon$ of $\mathrm{SU}_q(3)$, which on $\CP^2_q$
yields the `classical point', that is the character $\chi_0$ in  \eqref{eq:rank}.

On the other hand, with the integral defined in \eqref{int} by using the Haar state as well, an element $[\tau_4]\in H^4(\Aq,{_\eta}\Aq)$ is constructed as 
$$
\tau_4(a_0,\ldots,a_4):=\nint a_0\dd a_1\wprod\ldots\wprod\dd a_4 \;.
$$
Let us check that it is a cocycle. Leibniz rule gives
\begin{multline*}
b^*_5\tau_4(a_0,\ldots,a_5)
 =\nint a_0a_1\dd a_2\wprod a_3\wprod a_4\wprod\dd a_5  
-\nint a_0 (a_1\dd a_2+\dd a_1\,a_2)\wprod a_3\wprod a_4\wprod\dd a_5 \\
 +\nint a_0\dd a_1\wprod (a_2\dd a_3+\dd a_2\,a_3)\wprod a_4\wprod\dd a_5  
-\nint a_0\dd a_1\wprod a_2\wprod (a_3\dd a_4+\dd a_3\,a_4)\wprod\dd a_4 \\
+\nint a_0\dd a_1\wprod a_2\wprod a_3\wprod (a_4\dd a_5+\dd a_4\,a_5)  
 -\nint \eta(a_5)a_0\wprod a_1\dd a_2\wprod a_3\wprod a_4 \\
 =\nint a_0(\dd a_1\wprod\ldots\wprod \dd a_4)a_5
-\nint\eta(a_5)a_0(\dd a_1\wprod\ldots\wprod \dd a_4) \;,
\end{multline*}
which is zero by the modular property \eqref{eq:modprop} of the
Haar state.

\medskip
A $2$-cocycle can be defined in a similar way. Recall that elements of $\Omega^{1,1}(\CP^2_q)$
have the form $\omega=(\alpha,\alpha_4)$, with $\alpha_4\in\Aq$. Let $\pi:\Omega^{1,1}(\CP^2_q)\to\Aq$ be the projection onto the second component $\pi(\omega)=\alpha_4$, and extend it to a projection $\pi:\Omega^2(\CP^2_q)\to\Aq$
by setting $\pi(\omega)=0$ if $\omega\in\Omega^{0,2}$ or $\omega\in\Omega^{2,0}$. 
The map $\pi$ is an $\Aq$-bimodule map. Then, the map
$$
\tau_2(a_0,a_1,a_2):=\varphi\circ\pi(a_0\dd a_1\wprod \dd a_2)
$$
is the representative of a class $[\tau_2]\in H^2(\Aq,{_\eta}\Aq)$. Indeed, by 
the Leibniz rule,
$$
b^*_3\tau_2(a_0,a_1,a_2,a_3)=\varphi\circ\pi\bigl(
a_0(\dd a_1\wprod \dd a_2)a_3-\eta(a_3)a_0(\dd a_1\wprod \dd a_2)\bigr) \;.
$$
Being $\pi$ a bimodule map we get in turn 
$$
b^*_3\tau_2(a_0,a_1,a_2,a_3)=\varphi\bigl(
a_0\pi(\dd a_1\wprod \dd a_2)a_3-\eta(a_3)a_0\pi(\dd a_1\wprod \dd a_2)\bigr) \;,
$$
which is zero by the modular property of the Haar state. 

Both classes $[\tau_4]$ and $[\tau_2]$ will be proven to be not trivial by pairing them with the 
monopole projections \eqref{mon-pro}. Firstly, 
\begin{lemma}
The monopole projections $P_N=\Psi_N\Psi_N^\dag$ in \eqref{mon-pro} are equivariant with respect to the representation $\sigma^N$ of $\Uq$ defined as 
$$
\sigma^N(h):= \begin{cases}
\rho^{0,N}(S(h))^t &\; \mathrm{if}\;N\geq 0\;,\\
\sigma^N(h):=\rho^{-N,0}(S(h))^t &\; \mathrm{if}\;N\leq 0\;. 
\end{cases}
$$
\end{lemma}
\begin{proof}
By Prop.~\ref{lemma:psit}, 
$$
h\az\Psi^\dag_N=\begin{cases}
\Psi_N^\dag\rho^{0,N}(h) &\; \mathrm{if}\;N\geq 0\;,\\
\Psi_N^\dag\rho^{-N,0}(h) &\; \mathrm{if}\;N\leq 0\;,
\end{cases}
$$
for all $h\in\Uq$. Using $(h\az a^*)^*=S(h)^*\az a$, we get $h\az\Psi_N=\sigma^N(h)^t\Psi_N$, and the equivariance of $P_N$ follows.
\end{proof}

{}From the construction of the previous section, we can pair 
$\mathrm{ch}^4(P_N,\sigma^N)$ with $[\tau_4]$ and 
$\mathrm{ch}^2(P_N,\sigma^N)$ with $[\tau_2]$. The pairing of $[\tau_2]$ with $\mathrm{ch}^2(P_N,\sigma^N)$ gives
\begin{align*}
\inner{\tau_2,\mathrm{ch}^2(P_N,\sigma^N)}
 &=\varphi\,\tr\bigl(P_N\pi(\dd P_N\wprod\dd P_N)\sigma^N(K_1^{-4}K_2^{-4})^t\bigr) \\
 &=q^{-2N}\varphi\bigl(\Psi_N^\dag\pi(\dd P_N\wprod\dd P_N)\Psi_N\bigr) \;.
\end{align*}
Since $\Psi_N$ are `functions' on the total space of the bundle, we cannot move them
inside $\pi$ (which is an $\Aq$-bimodule map, not an $\Oq$-bimodule map). Nevertheless,
 -- with a little abuse of notations -- the form $\nabla_N^2=\pi(\nabla_N^2)$ is a constant, and
$$
P_N\dd P_N\wprod\dd P_N=\Psi_N\nabla^2_N\Psi_N^\dag
=\pi(\nabla_N^2)\Psi_N\Psi_N^\dag=\pi(\nabla_N^2)P_N \;.
$$
With this, and using $\Psi_N^\dag  P_N\Psi_N=1$, we come to the final formula
$$
\inner{\tau_2,\mathrm{ch}^2(P_N,\sigma^N)}
 =q^{-2N}\varphi\circ\pi(\nabla_N^2) \;.
$$
In \eqref{curvN1}, we have already shown  that $\nabla_N^2=q^{N-1} [N] \,\nabla_1^2$ 
with $ \nabla_1^2=(0,w_1)$. Thus the corresponding quantum characteristic class
is proportional to 
$q^{-N}[N]$:
$$
\inner{\tau_2,\mathrm{ch}^2(P_N,\sigma^N)}
 = \left( \varphi\circ\pi(\nabla_1^2) \right) q^{-N-1}[N] \,. 
$$ 
At $q=1$ the integral of the curvature is (modulo a global normalization constant)
the \emph{monopole number} of the bundle; it is the same as the first Chern number.

As for the pairing of $\mathrm{ch}^4(P_N,\sigma^N)$ with $\tau_4$, 
$$
\inner{\tau_4,\mathrm{ch}^4(P_N,\sigma^N)}=\nint\tr\bigl(P_N(\dd P_N)^4\sigma^N(K_1^{-4}K_2^{-4})^t\bigr) \;.
$$
Using the modular properties of the Haar state, which means
$$
\nint\tr(\Psi_N V)=\nint \tr\left( V(K_1K_2)^4\az\Psi_N\za (K_1K_2)^{-4} \right)=
q^{-2N}\nint \tr\left( V\,\sigma^N(K_1^4K_2^4)^t\Psi_N \right)
$$
and is valid for any row vector $V$ with entries in $\Oq$, we get:
$$
\inner{\tau_4,\mathrm{ch}^4(P_N,\sigma^N)}=q^{-2N}\nint \Psi^\dag_N(\dd P_N)^4\Psi_N \;.
$$
We need to compute the top form $\Psi^\dag_N(\dd P_N)^4\Psi_N$. {}From the identity \eqref{eq:ede}
it follows that $P_N(\dd P_N)^2\Psi_N=(\dd P_N)^2P_N\Psi_N=(\dd P_N)^2\Psi_N$, and
\begin{multline*}
\Psi^\dag_N(\dd P_N)^4\Psi_N=\Psi^\dag_N(\dd P_N)^2\wprod P_N(\dd P_N)^2\Psi_N \\
=\Psi^\dag_N(\dd P_N)^2\Psi_N\wprod \Psi^\dag_N(\dd P_N)^2\Psi_N=\nabla_N^2\wprod\nabla_N^2 \;,
\end{multline*}
leading to
$$
\inner{\tau_4,\mathrm{ch}^4(P_N,\sigma^N)}=q^{-2N}\nint \nabla_N^2\wprod\nabla_N^2 \;.
$$
Using $\nabla_N^2=q^{N-1} [N] \,\nabla_1^2$ the corresponding quantum characteristic class
is found to be proportional to $[N]^2$:  
$$
\inner{\tau_4,\mathrm{ch}^4(P_N,\sigma^N)} = 
\left(q^{-2} \, \nint \nabla_1^2\wprod\nabla^2_1 \right) [N]^2  \;.
$$
At $q=1$, the integral of the square of the curvature is (modulo a global normalization constant)
the \emph{instanton number} of the bundle. The pairing of a projection with the third
Fredholm module as in Sect.~\ref{se:inpm} does not give the `classical' instanton number, that is $N^2$, but rather the 2nd Chern number which is combination of the instanton number and of the monopole number.

By pairing $\mathrm{ch}^0(P_N,\sigma^N)$ with the Haar state, 
and using its modular property, one gets
\begin{align*}
\inner{\varphi,\mathrm{ch}^0(P_N,\sigma^N)}&=\varphi\bigl(\tr\,P_N\sigma^N(K_1^{-4}K_2^{-4})\bigr)
=q^{-2N}\varphi(\Psi_N^\dag\Psi_N) \\
& =q^{-2N} \;.
\end{align*}
On the other hand, the pairing with the classical point $\chi_0$ yields
\begin{align*}
\inner{\chi_0,\mathrm{ch}^0(P_N,\sigma^N)}&=\tr\,\chi_0(P_N)\sigma^N(K_1^{-4}K_2^{-4})=
\inner{0,0,0|K_1^4K_2^4|0,0,0} \\
&=q^{2N} \;,
\end{align*}
with $\ket{0,0,0}$ the highest weight of the representation $\rho^{0,N}$ if $N\geq 0$ or $\rho^{-N,0}$ if $N\leq 0$.

As a byproduct, we see that if $q$ is transcendental, the equivariant $K_0$-group has (at least) a countable number of generators:
$K_0^{\Uq}(\Aq)\supset\Z^\infty$, i.e.~all $[P_N]$ are independent. Indeed, were the classes $[P_N]$
not independent, there would exist a sequence $\{k_N\}$ of integers -- all zero but for finitely many -- such that $\sum_Nk_Nq^{2N}=0$, 
and $q$ would be the root of a non-zero polynomial with integer coefficients.

The results above are analogous to the ones for the standard Podle\'s sphere
(cf.~Prop.~5.1 and 5.2 in \cite{Wag07}). 
In fact they are instances of the general fact (cf. \cite[Thm.~$3.6$]{NT03}) 
that the equivariant $K_0$ group is a free abelian group and generators are, for the present case, 
in bijection with equivalence classes of irreducible corepresentations of $\Kq$.  


\section{Concluding remarks}

On a four-dimensional manifold (anti)self-dual connections are stationary points (usually minima) of the Yang-Mills action functional, i.e.~they are solutions of the corresponding equations of motion. In dimension greater than four, `generalized' instantons can be defined as solutions of Hermitian Yang-Mills equations.
On $\CP^n$ a basic instanton solution is associated to the canonical (universal) connection on the Stiefel bundle $\mathrm{U}(n)\hookrightarrow \mathrm{U}(n+1)/\mathrm{U}(1)\to\CP^n$. The extension of this construction to \emph{quantum} complex projective spaces -- using the differential calculus in \cite{DD09} -- will be explored in future works, and it should help, in particular, to understand how to generalize Hermitian Yang-Mills equations to noncommutative spaces. 

\appendix

\section{Proof  of Proposition.~{\protect\ref{pr:const}}}\label{app:B}
In this appendix, we determine the most general value of the normalization constants in
order to have a left $\Kq$-covariant product on $V^{\bullet,\bullet}$ which is
i) associative, ii) graded commutative for $q=1$, and iii) it sends real vectors into
real vectors. 

Indeed, as a way of illustration, let us start by considering the cases $V^{0,1}\times V^{1,0}\to V^{1,1}$
and $V^{0,1}\times V^{0,1}\to V^{0,2}$.
As vector spaces $V^{0,1}\simeq V^{1,0}\simeq\C^2$. For $v,w\in\C^2$ we order
the components of $v\otimes w$ as: $v\otimes w=(v_1w_1,v_1w_2,\, v_2w_1,v_2w_2)^t$.
A unitary equivalence $U$ between 
$\sigma_{\frac{1}{2},N}\otimes\sigma_{\frac{1}{2},N'}$ and
$\sigma_{1,N+N'}\oplus\sigma_{0,N+N'}$ is given by
\begin{equation}\label{eq:u}
U=\ma{
  1 & 0 & 0 & 0 \\
  0 & q^{-\frac{1}{2}}[2]^{-\frac{1}{2}} & q^{\frac{1}{2}}[2]^{-\frac{1}{2}} & 0 \\
  0 & 0 & 0 & 1 \\
  0 & q^{\frac{1}{2}}[2]^{-\frac{1}{2}} & -q^{-\frac{1}{2}}[2]^{-\frac{1}{2}} & 0
  }
\end{equation}
It is easy to check that
$U(\sigma_{\frac{1}{2},N}(h_{(1)})v\otimes\sigma_{\frac{1}{2},N'}(h_{(2)})w)=
(\sigma_{1,N+N'}(h)\oplus\sigma_{0,N+N'}(h))U(v\otimes w)$ for all $h\in\Kq$ by doing it explicitly on all generators.  For $h=K_1,K_1K_2^2$ this is trivial.
For $h=E_1$, and omitting the representation symbols, we have
$$
U(E_1v\otimes K_1w+K_1^{-1}v\otimes E_1w)=U\ma{
q^{-\frac{1}{2}}v_2w_2 \\
q^{\frac{1}{2}}v_2w_2 \\
0 \\
q^{-\frac{1}{2}}v_1w_2+q^{\frac{1}{2}}v_2w_1
}=\ma{q^{-\frac{1}{2}}v_1w_2+q^{\frac{1}{2}}v_2w_1 \\ {}[2]^{\frac{1}{2}}v_2w_2 \\ 0 \\ 0 } \;,
$$
and
$$
E_1U(v\otimes w)=[2]^{\frac{1}{2}}
\ma{0 & 1 & 0 & 0 \\ 0 & 0 & 1 & 0 \\ 0 & 0 & 0 & 0 \\ 0 & 0 & 0 & 0}U(v\otimes w)=
\ma{0 & q^{-\frac{1}{2}} & q^{\frac{1}{2}} & 0 \\ 0 & 0 & 0 & {}[2]^{\frac{1}{2}} \\ 0 & 0 & 0 & 0 \\ 0 & 0 & 0 & 0}
(v\otimes w) \;.
$$
Thus $U\Delta(E_1)=E_1U$. Then the statement holds for $F_1$ since $F_1=E_1^*$.
The map $\wprod:V^{0,1}\times V^{1,0}\to V^{1,1}$ is then $v\wprod w=\mathrm{diag}({c_1},{c_1},
{c_1},{c_2})U(v\otimes w)$, where ${c_1},{c_2}\in\R$ are arbitrary for the time being.
When ${c_i}=\pm 1$ we would get partial isometries; 
by composing $U$ with the orthogonal projection onto the last component, we would 
get a partial isometry from $V^{0,1}\times V^{0,1}\to V^{0,2}$.

The general situation is listed in the following proposition.
\begin{prop}\label{wedpro}
The most general left $\Kq$-covariant graded product $\wprod$ on $V^{\bullet,\bullet}$,
sending real vectors to real vectors, is given by
\begin{align*}
V^{0,1}\times V^{0,1} &\to V^{0,2}\;,
   & v\wprod w &:={c_0}\mu_0(v,w)^t\;,\\
V^{0,1}\times V^{1,0} &\to V^{1,1}\;,
   & v\wprod w &:=\bigl({c_1}\mu_1(v,w), \,{c_2}\mu_0(v,w)\bigr)^t \;,\\
V^{0,1}\times V^{2,1} &\to V^{2,2}\;,
   & v\wprod w &:={c_3}\mu_0(v,w)^t\;,\\
V^{0,1}\times V^{1,1} &\to V^{1,2}\;,
   & v\wprod w &:={c_0^1}[3]^{-\frac{1}{2}}\mu_2(v,w)^t+{c_0^2}vw_4 \;,\\
V^{1,0}\times V^{1,0} &\to V^{2,0}\;,
   & v\wprod w &:={c_4}\mu_0(v,w)^t\;,\\
V^{1,0}\times V^{0,1} &\to V^{1,1}\;,
   & v\wprod w &:=\bigl(-{d_0}\mu_1(v,w), \,{d_1}\mu_0(v,w)\bigr)^t \;,\\
V^{1,0}\times V^{1,2} &\to V^{2,2}\;,
   & v\wprod w &:={d_2}\mu_0(v,w)^t\;,\\
V^{1,0}\times V^{1,1} &\to V^{2,1}\;,
   & v\wprod w &:={c_1^1}[3]^{-\frac{1}{2}}\mu_2(v,w)^t+{c_1^2}vw_4 \;,\\
V^{1,2}\times V^{1,0} &\to V^{2,2}\;,
   & v\wprod w &:={d_3}\mu_0(v,w)^t\;,\\
V^{2,1}\times V^{0,1} &\to V^{2,2}\;,
   & v\wprod w &:={d_4}\mu_0(v,w)^t\;,\\
V^{1,1}\times V^{0,1} &\to V^{1,2}\;,
   & v\wprod w &:=-{c_2^1}[3]^{-\frac{1}{2}}\mu_3(v,w)^t+{c_2^2}v_4w \;,\\
V^{1,1}\times V^{1,0} &\to V^{2,1}\;,
   & v\wprod w &:=-{c_3^1}[3]^{-\frac{1}{2}}\mu_3(v,w)^t+{c_3^2}v_4w \;,\\
V^{1,1}\times V^{1,1} &\to V^{2,2}\;,
   & v\wprod w&:={c_4^1}[3]^{-\frac{1}{2}}\mu_4(v,w)+{c_4^2}v_4 w_4 \;. 
\end{align*}
The coefficients $c_i, d_j$ and $c_i^j\in\R$ are arbitrary for the time being; and 
the maps $\mu_i$'s are 
\begin{align*}
\mu_0:\R^2\times\R^2&\to\R\;, &
\mu_0(v,w)&:=[2]^{-\frac{1}{2}}(q^{\frac{1}{2}}v_1w_2-q^{-\frac{1}{2}}v_2w_1)\;,\\
\mu_1:\R^2\times\R^2&\to\R^3\;, &
\mu_1(v,w)&:=\bigl(v_1w_1,[2]^{-\frac{1}{2}}(q^{-\frac{1}{2}}v_1w_2+q^{\frac{1}{2}}v_2w_1), \,v_2w_2\bigr)\;,\\
\mu_2:\R^2\times\R^3&\to\R^2\;, &
\mu_2(v,w)&:=\bigl(qv_1w_2-q^{-\frac{1}{2}}[2]^{\frac{1}{2}}v_2w_1, \,
              q^{\frac{1}{2}}[2]^{\frac{1}{2}}v_1w_3-q^{-1}v_2w_2\bigr)\;,\\
\mu_3:\R^3\times\R^2&\to\R^2\;, &
\mu_3(v,w)&:=\bigl(q^{\frac{1}{2}}[2]^{\frac{1}{2}}v_1w_2-q^{-1}v_2w_1, \,
              qv_2w_2-q^{-\frac{1}{2}}[2]^{\frac{1}{2}}v_3w_1\bigr)\;,\\
\mu_4:\R^3\times\R^3&\to\R\;, &
\mu_4(v,w)&:=qv_1w_3-v_2w_2+q^{-1}v_3w_1 \;.
\end{align*}
On the other hand, there is no need to to specify the multiplication rule by elements in $V^{0,0}=V^{0,2}=V^{2,0}=V^{2,2}=\C$,
being scalars.
\end{prop}

\smallskip

When all the coefficients $c_i$'s and $d_i$'s are equal to $\pm 1$ and $(c_i^1)^2+(c_i^2)^2=1$, the maps in the proposition are partial isometries. In general, not all the choices give an associative product. Associativity fixes the value of the parameters $d_i$'s and $c_i^j$'s.

\begin{lemma}
We have
\begin{subequations}\label{eq:dueto}
\begin{align}
\mu_4(\mu_1(v,v'), w)&=\mu_0(v,\mu_2(v',w))                 \label{eq:duetoA}\;,\\
\mu_0(\mu_2(v,w), v')&=\mu_0(v,\mu_3(w,v'))                 \label{eq:duetoB}\;,\\
\mu_0(\mu_3(w,v), v')&=\mu_4(w,\mu_1(v,v'))                 \label{eq:duetoC}\;,\\
\mu_2(v,\mu_1(v',v''))&=[2] \, \mu_0(v,v')v''+v\,\mu_0(v',v'') \label{eq:duetoD}\;,\\
\mu_3(\mu_1(v,v'), v'')&=\mu_0(v,v')v''+[2] \,v\,\mu_0(v',v'') \label{eq:duetoE}\;,
\end{align}
\end{subequations}
for all $v,v',v''\in\C^2$ and $w\in\C^3$.
\end{lemma}

\begin{proof}
By direct computation. Both sides of (\ref{eq:duetoA}) are equal to
$$
qv_1v'_1w_3-[2]^{-\frac{1}{2}}(q^{-\frac{1}{2}}v_1v'_2+q^{\frac{1}{2}}v_2v'_1)w_2+q^{-1}v_2v'_2w_1 \;,
$$
both sides of (\ref{eq:duetoB}) are equal to
$$
q^{\frac{3}{2}}[2]^{-\frac{1}{2}}v_1w_2v'_2-v_1w_3v'_1-v_2w_1v'_2+q^{-\frac{3}{2}}[2]^{-\frac{1}{2}}v_2w_2v'_1 \;,
$$
both sides of (\ref{eq:duetoC}) are equal to
$$
qv_1w_2v'_2-[2]^{-\frac{1}{2}}v_2(q^{-\frac{1}{2}}w_1v'_2+q^{\frac{1}{2}}w_2v'_1)+q^{-1}v_3w_1v'_1 \;,
$$
both sides of (\ref{eq:duetoD}) are equal to
$$
\binom{
q^{\frac{1}{2}}[2]^{-\frac{1}{2}}v_1v'_1v''_2
+q^{\frac{3}{2}}[2]^{-\frac{1}{2}}v_1v'_2v''_1
-q^{-\frac{1}{2}}[2]^{\frac{1}{2}}v_2v'_1v''_1
}{\rule{0pt}{18pt}
q^{\frac{1}{2}}[2]^{\frac{1}{2}}v_1v'_2v''_2
-q^{-\frac{3}{2}}[2]^{-\frac{1}{2}}v_2v'_1v''_2
-q^{-\frac{1}{2}}[2]^{-\frac{1}{2}}v_2v'_2v''_1
} \;,
$$
both sides of (\ref{eq:duetoE}) are equal to
$$
\binom{
q^{\frac{1}{2}}[2]^{\frac{1}{2}}v_1v'_1v''_2
-q^{-\frac{3}{2}}[2]^{-\frac{1}{2}}v_1v'_2v''_1
-q^{-\frac{1}{2}}[2]^{-\frac{1}{2}}v_2v'_1v''_1
}{\rule{0pt}{18pt}
q^{\frac{1}{2}}[2]^{-\frac{1}{2}}v_1v'_2v''_2
+q^{\frac{3}{2}}[2]^{-\frac{1}{2}}v_2v'_1v''_2
-q^{-\frac{1}{2}}[2]^{\frac{1}{2}}v_2v'_2v''_1
} \;.
$$
\end{proof}

\begin{prop}
The map $\wprod$ in Prop.~\ref{wedpro} is a graded associative product on $V^{\bullet,\bullet}$ if
\begin{subequations}\label{subeq:coef}
\begin{align}
{d_0}&={s_1}q^{\frac{1}{2}{s_2}}{c_1} \;, &
{d_1}&={s_1}q^{-\frac{3}{2}{s_2}}{c_2} \;, \label{eq:coefA} \\
{c^1_1}&=\frac{\sqrt{[3]}}{[2]}{\frac{c_0}{c_1}} \;, &
{c^1_3}&={s_1}q^{-\frac{1}{2}{s_2}}\frac{\sqrt{[3]}}{[2]}{\frac{c_0}{c_1}} \;, \label{eq:coefB} \\
{c^2_1}&=-\frac{1}{[2]}{\frac{c_0}{c_2}} \;, &
{c^2_3}&=-{s_1}q^{\frac{3}{2}{s_2}}\frac{1}{[2]}{\frac{c_0}{c_2}} \;, \label{eq:coefC} \\
{c^1_2}&=-{s_1}q^{-\frac{1}{2}{s_2}}\frac{\sqrt{[3]}}{[2]}{\frac{c_4}{c_1}} \;, &
{c^1_4}&=-\frac{\sqrt{[3]}}{[2]}{\frac{c_4}{c_1}} \;, \label{eq:coefD} \\
{c^2_2}&=-{s_1}q^{\frac{3}{2}{s_2}}\frac{1}{[2]}{\frac{c_4}{c_2}} \;, &
{c^2_4}&=-\frac{1}{[2]}{\frac{c_4}{c_2}} \;, \label{eq:coefE} \\
{c^1_5}&=-{s_1}q^{-\frac{1}{2}{s_2}}\frac{\sqrt{[3]}}{[2]}{\frac{c_3c_4}{|c_1|^2}} \;, &
{c^2_5}&=-{s_1}q^{\frac{3}{2}{s_2}}\frac{1}{[2]}{\frac{c_3c_4}{|c_2|^2}} \;, \label{eq:coefF} \\
{d_2} &={\frac{c_3c_4}{c_0}} \;, &
{d_3} &={\frac{c_3c_4}{c_0}} \;, \label{eq:coefG} \\
{d_4} &={c_3} \;, \label{eq:coefH} \;,
\end{align}
\end{subequations}
where $\{{s_1},{s_2}\}\in\{\pm 1\}$ are arbitrary signs.
The algebra is graded commutative in the $q\to1$ limit if{}f ${s_1}=1$.
\end{prop}

\begin{proof}
We seek solutions with all parameters different from zero
The product is graded by construction. We impose the condition:
$$
(v\wprod v')\wprod v''=
v\wprod(v'\wprod v'') \;.
$$
If one of the three vectors $v,v',v''$ is a scalar, i.e.~an element of
$V^{0,0}$, $V^{0,2}$, $V^{2,0}$ or $V^{2,2}$, the equality follows by the
bilinearity of $\wprod$. The non-trivial cases are when all three vectors
are not scalars. If the total degree is greater than $4$, we get $0=0$.
The remaining non-trivial cases are
firstly $(v,v',v'')\in V^{0,1}\times V^{0,1}\times V^{1,0}$
and (two) permutations:
\begin{align*}
(v\wprod v')\wprod v'' &={c_0}\mu_0(v,v')v''  \;,\\
v\wprod (v'\wprod v'') &={c_1}{c^1_1}[3]^{-\frac{1}{2}}\mu_2(v,\mu_1(v',v''))
                         +{c_2}{c^2_1}v\mu_0(v',v'') \;,\\[8pt]
(v\wprod v'')\wprod v' &=-{c_1}{c^1_3}[3]^{-\frac{1}{2}}\mu_3(\mu_1(v,v''),v')
                         +{c_2}{c^2_3}\mu_0(v,v'')v' \;,\\
v\wprod (v''\wprod v') &=-{d_0}{c^1_1}[3]^{-\frac{1}{2}}\mu_2(v,\mu_1(v'',v'))
                         +{d_1}{c^2_1}v\mu_0(v'',v') \;,\\[8pt]
(v''\wprod v)\wprod v' &={d_0}{c^1_3}[3]^{-\frac{1}{2}}\mu_3(\mu_1(v'',v),v')
                         +{d_1}{c^2_3}\mu_0(v'',v)v' \;,\\
v''\wprod (v\wprod v') &={c_0}v''\mu_0(v,v') \;.
\end{align*}
Using (\ref{eq:duetoD}-\ref{eq:duetoE}) one sees that associativity is just the conditions 

(\ref{eq:coefA}-\ref{eq:coefC}) respectively. 
\\
Then,  $(v,v',v'')\in V^{1,0}\times V^{1,0}\times V^{0,1}$
and (two) permutations:
\begin{align*}
(v\wprod v')\wprod v'' &={c_4}\mu_0(v,v')v'' \;, \\
v\wprod (v'\wprod v'') &=-{d_0c^1_2}[3]^{-\frac{1}{2}}\mu_2(v,\mu_1(v',v''))+{d_1c^2_2}v\mu_0(v',v'') \;, \\[8pt]
(v\wprod v'')\wprod v' &={d_0c^1_4}[3]^{-\frac{1}{2}}\mu_3(\mu_1(v,v''),v')+{d_1c^2_4}\mu_0(v,v'')v' \;, \\
v\wprod (v''\wprod v') &={c_1}{c^1_2}[3]^{-\frac{1}{2}}\mu_2(v,\mu_1(v'',v'))+{c_2}{c^2_2}v\mu_0(v'',v') \;,\\[8pt]
(v''\wprod v)\wprod v'&=-{c_1}{c^1_4}[3]^{-\frac{1}{2}}\mu_3(\mu_1(v'',v),v')+{c_2}{c^2_4}\mu_0(v'',v)v' \;,\\
v''\wprod (v\wprod v')&={c_4}v''\mu_0(v,v') \;.
\end{align*}
Using again (\ref{eq:duetoD}-\ref{eq:duetoE}) one sees that associativity is just  the condition (\ref{eq:coefD}-\ref{eq:coefE}) respectively. 
\\
Finally, $(v,v',v'')\in V^{0,1}\times V^{1,0}\times V^{1,1}$ and (five) permutations: 
\begin{align*}
\{v\wprod v'\}\wprod v''&={c_2}{c^2_5}\mu_0(v,v')w_4+{c_1}{c^1_5}[3]^{-\frac{1}{2}}\mu_4(\mu_1(v,v'),v'') \;,\\
v\wprod \{v'\wprod v''\}&={c_3}{c^2_2}\mu_0(v,v'w_4)+{c_3}{c^1_2}[3]^{-\frac{1}{2}}\mu_0(v,\mu_2(v',v'')) \;,\\[8pt]
\{v\wprod v''\}\wprod v'&={c^2_1d_3}\mu_0(vw_4,v')+{d_3c^1_1}[3]^{-\frac{1}{2}}\mu_0(\mu_2(v,v''),v') \;,\\
v\wprod \{v''\wprod v'\}&={c_3}{c^2_4}\mu_0(v,w_4v')-{c_3}{c^1_4}[3]^{-\frac{1}{2}}\mu_0(v,\mu_3(v'',v')) \;,\\[8pt]
\{v'\wprod v\}\wprod v''&={d_1c^2_5}\mu_0(v',v)w_4-{d_0c^1_5}[3]^{-\frac{1}{2}}\mu_4(\mu_1(v',v),v'') \;,\\
v'\wprod \{v\wprod v''\}&={d_2c^2_1}\mu_0(v',vw_4)+{d_2c^1_1}[3]^{-\frac{1}{2}}\mu_0(v',\mu_2(v,v'')) \;,\\[8pt]
\{v'\wprod v''\}\wprod v&={d_4c^2_2}\mu_0(v'w_4,v)+{d_4c^1_2}[3]^{-\frac{1}{2}}\mu_0(\mu_2(v',v''),v) \;,\\
v'\wprod \{v''\wprod v\}&={d_2c^2_3}\mu_0(v',w_4v)-{d_2c^1_3}[3]^{-\frac{1}{2}}\mu_0(v',\mu_3(v'',v)) \;,\\[8pt]
\{v''\wprod v\}\wprod v'&={d_3c^2_3}\mu_0(w_4v,v')-{d_3c^1_3}[3]^{-\frac{1}{2}}\mu_0(\mu_3(v'',v),v') \;,\\
v''\wprod \{v\wprod v'\}&={c_2}{c^2_5}w_4\mu_0(v,v')+{c_1}{c^1_5}[3]^{-\frac{1}{2}}\mu_4(v'',\mu_1(v,v')) \;,\\[8pt]
\{v''\wprod v'\}\wprod v&={d_4c^2_4}\mu_0(w_4v',v)-{d_4c^1_4}[3]^{-\frac{1}{2}}\mu_0(\mu_3(v'',v'),v) \;,\\
v''\wprod \{v'\wprod v\}&={d_1c^2_5}w_4\mu_0(v',v)-{d_0c^1_5}[3]^{-\frac{1}{2}}\mu_4(v'',\mu_1(v',v)) \;.
\end{align*}
Using (\ref{eq:dueto}), associativity for the first couple is shown to be equivalent to (\ref{eq:coefF}),  for the second couple to the second condition in (\ref{eq:coefG}), for the third couple to the first condition in (\ref{eq:coefG}), for the fourth couple to (\ref{eq:coefH}). For the last two 
couples the associativity is automatically satisfied.

Graded anticommutativity in the $q\to1$ limit is a simple check, based on the observation
that for $q=1$: $\mu_0$ is antisymmetric, $\mu_1$ is symmetric, $\mu_2(v',v'')=-\mu_3(v'',v')$,
and $\mu_4$ is symmetric. This concludes the proof.
\end{proof}

\section{Some general facts on calculi}\label{ap:cal}

In this appendix we review the material that leads to the construction of the differential calculus 
$(\Omega^\bullet(\CP^2_q), \dd)$ on $\CP^2_q$, in particular the fact that it is enough to define things for one-forms and then extend them in a natural and unique (by universality) way. 

Recall that a first order differential calculus on a unital $*$-algebra $\A$ is a pair
$(\Omega^1\A,\dd )$, where $\Omega^1A$ is an $\A$-bimodule giving
the space of one-forms and $\dd : \A \to \Omega^1\A$ is the exterior differential -- a linear map
satisfying the Leibniz rule,
$$
\dd (a b)=a(\dd b)+(\dd  a)b \qquad \text{for all} \quad a,b \in \A.
$$
We also assumes that $\Omega^1\A= \A (\dd \A)$ and that we are dealing with 
a $*$-calculus with the $*$-structure on $\Omega^1\A$ given (uniquely) by 
$(\dd a)^*=-\dd (a^*)$ for all $a \in A$. 

To go beyond one forms a differential graded $*$-algebra $(\Omega^\bullet=\oplus_k\Omega^k\A ,\dd)$ is defined as follows.
The set $\Omega^\bullet\A$ is a graded algebra, with $\Omega^0\A=\A$ and graded product denoted $\wprod:\Omega^k\A \times \Omega^l\A \to \Omega^{k+l}\A$. 
The \emph{universal} differential calculus of $\A$ is obtained when $\wprod$ is the tensor product, $\Omega^1\A$ is the kernel of the multiplication map $m(a\otimes b)=a b$ and $\dd a=a\otimes 1-1\otimes a$.  With a  general associative (unital) graded multiplication $\wprod$, the calculus will always be a quotient of the universal calculus by a differential ideal.

Element in $\Omega^k\A$ (that is \emph{$k$-forms}) are (sums of\hspace{2pt}) products of $k$ 1-forms: 
$\omega=v_1\wprod\ldots\wprod v_k$.
The differential is extended to $\Omega^k\A$ by requiring that its square vanishes, $\dd^2=0$, and that it is a graded derivation of degree 1, that is $\dd : \Omega^k\A \to \Omega^{k+1}\A$ and
$$
\dd(\omega\wprod\omega')=\dd\omega\wprod\omega'+(-1)^{\mathrm{dg}(\omega)}\omega\wprod\dd\omega' \;,
$$
for any two forms $\omega$ and $\omega'$.
By universality, these two properties uniquely identify $\dd$. It is the map given on $1$-forms by
$$
\dd(a\dd b)=\dd a\wprod\dd b=-\dd (\dd a\,b)
$$ 
on any product of $k$ $1$-forms by 
$$ 
\dd(v_1\wprod\ldots\wprod v_k)=\sum_{i=1}^k(-1)^{i-1}v_1\wprod\ldots\wprod\dd v_i\wprod\ldots\wprod v_k \;,
$$ 
and extended by linearity. Finally, a $*$-structure on a differential calculus is given by a graded involution which anti-commutes with the differential. Given a $*$-structure on $\A$, a $*$-structure on
$(\Omega^\bullet\A,\dd)$ is uniquely defined firstly on $1$-forms by
$(a\dd b)^*=-(\dd b^*)a^*$, and then by induction on $k$-forms by
$$
(v\wprod w)^* =(-1)^{k-1} w^*\wprod v^*\;,
$$
where $v$ is a $1$-form and $w$ is a $k-1$ form. With this, one easily check, again by induction, that $\dd\omega^* =-(\dd\omega)^*$ for any form $\omega$. 

An equivalent way to give a $*$-differential calculus $(\Omega^\bullet\A,\dd)$ is via a differential  double complex $(\Omega^{\bullet,\bullet}\A,\de,\deb)$. That is to say, starting with a bigradation $\Omega^k\A:=\bigoplus_{i+j=k}\Omega^{i,j}\A$ with a bigraded product,
the following statements are equivalent:
\begin{itemize}
\item[i)] there is a graded derivation $\dd:\Omega^{\bullet}\A\to\Omega^{\bullet+1}\A$ satisfying
$$
\dd^2=0 \qquad \textup{and} \qquad \dd\omega=-(\dd\omega^*)^*\,;
$$
\item[ii)] there are two graded derivations $\de:\Omega^{\bullet,\bullet}\A\to
\Omega^{\bullet+1,\bullet}\A$ and
$\deb:\Omega^{\bullet,\bullet}\A\to\Omega^{\bullet,\bullet+1}\A$, satisfying 
$$
\de^2=\deb^2=\de\deb+\deb\de=0 \qquad \textup{and} \qquad 
\deb\omega=-(\de\omega^*)^*\,.
$$
\end{itemize}


\section{An alternative proof of Theorem \ref{thm} for $n\leq 4$}\label{app:alt}

We give an alternative proof that the map $\mathrm{ch}^n:K_0^{\U}(\A)\to H\!C^{\U}_n(\A)$ 
given by
\begin{equation}\label{eq:eccc}
\mathrm{ch}^n(e,\sigma)(h):=\tr(e^{\dotimes n+1}\sigma(h)^t) \;,
\end{equation}
is well-defined. This proof is valid for $n\leq 4$ and for it we use the results of 
Lemma~\ref{lemma:KU} stating the equivalence of equivariant projective modules in terms of idempotents.

Thus, we aim at proving that for any $[(e,\sigma)]\in K_0^{\U}(\A)$, the expression in (\ref{eq:eccc}) is an equivariant cyclic $n$-cycle, and that if $(e,\sigma)$
is in the same class of $(e',\sigma')$,  the difference $\mathrm{ch}^n(e,\sigma)-\mathrm{ch}^n(e',\sigma')$ is a boundary. In fact, when $n$ is odd (\ref{eq:eccc}) 
$$
\mathrm{ch}^n(e,\sigma)=\tfrac{1}{2}(1-\lambda_n)\mathrm{ch}^n(e,\sigma) \;,
$$
and $\mathrm{ch}^n(e,\sigma)$ reduces in $C^{\U}_n(\A)/\mathrm{Im}(1-\lambda_n)$ to the trivial equivariant cyclic cocycle $0$. 

For $n$ even,  one finds that $b_n\mathrm{ch}^n(e,\sigma)=\mathrm{ch}^{n-1}(e,\sigma)$ which vanishes in $C^{\U}_n(\A)/\mathrm{Im}(1-\lambda_n)$ since $n-1$ is now odd and the class of $\mathrm{ch}^{n-1}(e,\sigma)$ is zero as before. Thus, $\mathrm{ch}^n(e,\sigma)$ is an equivariant cyclic cocycle. 

For the last part of the proof, $\mathrm{ch}^0(e,\sigma)-\mathrm{ch}^0(e',\sigma')$ is the boundary of
the function
$$
\tr\bigl(u\dotimes v\sigma(h)^t\bigr) \; ;
$$
$\mathrm{ch}^2(e,\sigma)-\mathrm{ch}^2(e',\sigma')$ is the boundary 
(modulo the image of $1-\lambda_2$) of
$$
\tr\bigl(e^{\dotimes 2}\dotimes u\dotimes v\sigma(h)^t\bigr)+\tr\bigl(u\dotimes e'^{\dotimes 2}\dotimes v\sigma(h)^t\bigr)
+\tfrac{1}{2}\tr\bigl(u\dotimes v\dotimes u\dotimes v\sigma(h)^t\bigr) \; ;
$$
and $\mathrm{ch}^4(e,\sigma)-\mathrm{ch}^4(e',\sigma')$ is the boundary (modulo the image of $1-\lambda_4$) of
\begin{multline*}
\tr\bigl(e^{\dotimes 4}\dotimes u\dotimes v\sigma(h)^t\bigr)+\tr\bigl(u\dotimes e'^{\dotimes 4}\dotimes v\sigma(h)^t\bigr) \\
-\tfrac{1}{2}\tr\bigl(e^{\dotimes 2}\dotimes u\dotimes v\dotimes u\dotimes v\sigma(h)^t\bigr)
-\tfrac{1}{2}\tr\bigl(u\dotimes e'^{\dotimes 2}\dotimes v\dotimes u\dotimes v\sigma(h)^t\bigr) \\
+\tfrac{1}{2}\tr\bigl(e^{\dotimes 2}\dotimes u\dotimes e'^{\dotimes 2}\dotimes v\sigma(h)^t\bigr) 
+\tfrac{1}{6}\tr\bigl(u\dotimes v\dotimes u\dotimes v\dotimes u\dotimes v\sigma(h)^t\bigr) \;.
\end{multline*}
This concludes the proof. In general, for an arbitrary high $n$, it is not an easy task to write 
$\mathrm{ch}^n(e,\sigma)-\mathrm{ch}^n(e',\sigma')$ explicitly as the boundary of something else.


\section{Irreducible $*$-representations of $\CP^2_q$}\label{app-IR}
In this appendix we prove the result mentioned at the end of Sect.~\ref{se:inpm}, namely  
that any non-trivial irreducible $*$-representation (irrep for short) of $\Aq$ (by \emph{bounded} operators on a Hilbert space)
is unitarily equivalent to one of the representations 
$\chi_0$, $\chi_1$ or $\chi_2$ described in the paper.

\begin{lemma}\label{lemma:app1}
Let $\pi:\Aq\to\B(\HH)$ be a $*$-representation such that
$\ker\pi(p_{11})=\{0\}$. Then there exists $Z_i\in\B(\HH)$,
$i=1,2,3$, such that $Z_1$ is positive and
\begin{equation}\label{eq:piP}
\pi(P)=\begin{pmatrix}
Z_1^2 & qZ_2Z_1 & qZ_3Z_1 \\
qZ_1Z_2^* & Z^*_2Z_2 & Z_2^*Z_3 \\
qZ_1Z_3^* & Z_3^*Z_2 & Z^*_3Z_3
\end{pmatrix} \;.
\end{equation}
With this, a $*$-representation $\tilde{\pi}:\Sq\to\B(\HH)$ is defined by
$\tilde{\pi}(z_i)=Z_i$ and it satisfies $\ker\tilde{\pi}(z_1)=\{0\}$.
Moreover, $Z^*_iZ_j=\pi(p_{ij})$, so that the restriction of $\tilde{\pi}$
to $\Aq$ is exactly $\pi$. The $*$-representation $\tilde{\pi}$ is
irreducible if $\pi$ is irreducible. 
\end{lemma}
\begin{proof}
Recall that the generators $p_{ij}$ of $\Aq$ satisfy the commutation relations
in Sect.~\ref{sec:2.3}, the quadratic relations $\sum_jp_{ij}p_{ji}=p_{ii}$,
$i=1,2,3$, and the real structure is given $p_{ij}=p_{ji}^*$. They are related
to the generators $z_i$ of $\Sq$ by $p_{ij}=z_i^*z_j$.

Since $\pi$ is a $*$-representation,
$\pi(p_{ii})=\sum_j\pi(p_{ij})\pi(p_{ji})=\sum_j\pi(p_{ji})^*\pi(p_{ji})$
is a sum of positive operators and then it is positive as well. The positive operators
$A_i=\pi(p_{ii})^{\frac{1}{2}}$ are well defined and are mutually
commuting since $p_{ii}$ are mutually commuting.

We have $(q^4p_{11}+q^2p_{22}+p_{33})p_{ii}=p_{ii}=\sum_jp_{ij}p_{ji}$, for all $i=1,2,3$.
Using the commutation rules between $p_{ij}$ and $p_{ji}$ yields
\begin{align*}
q^2p_{11}p_{22}+p_{11}p_{33} &= p_{21}p_{12}+q^{-2}p_{31}p_{13} \;, \\
q^4p_{11}p_{22}+p_{22}p_{33} &= q^2p_{21}p_{12}+(q^{-2}-1)p_{31}p_{13}+q^{-2}p_{32}p_{23} \;, \\
q^4p_{11}p_{33}+q^2p_{22}p_{33} &= p_{31}p_{13}+p_{32}p_{23} \;.
\end{align*}
By solving these equations we get
\begin{equation}\label{irrep:eq1}
p_{i1}p_{1i} = q^2p_{11}p_{ii} \;, \qquad i=2,3\;, 
\end{equation}
and a third relation for $p_{32}p_{23}$ that we don't need.
Similarly, using the commutation rule $[p_{33},p_{23}]=(1-q^2)(p_{21}p_{13}+p_{22}p_{23})$,
the relation $(q^4p_{11}+q^2p_{22}+p_{33})p_{23}=p_{23}=\sum_ip_{2i}p_{i3}$ can be rewritten as
\begin{equation}\label{irrep:eq2}
p_{11}p_{23}=q^{-2}p_{21}p_{13} \; .
\end{equation}
Writing $T=V|T|$ for the polar decomposition of a bounded operator $T$, for the 
polar decomposition of $\pi(p_{1i})$, $i=2,3$, from \eqref{irrep:eq1}
we get $|\pi(p_{1i})|^2=q^2A_1^2A_i^2=(qA_1A_i)^2$.
Therefore $\pi(p_{1i})=qV_iA_1A_i$. Applying $\pi$ to both sides of the equation \eqref{irrep:eq2} yields
\begin{equation}\label{irrep:eq3}
A_1\pi(p_{23})A_1=A_1A_2V_2^*V_3A_3A_1 \;,
\end{equation}
where we used the fact that $p_{11}$ and $p_{23}$ commute, and so do $A_1$
and $\pi(p_{23})$.

For any bounded operator $T$, it holds that $\ker(T)^\perp=\overline{\mathrm{range}(T^*)}$.
Being $A_1$ positive with $\ker(A_1)=\{0\}$, we have $\overline{\mathrm{range}(A_1)}=\HH$.
Since the range of $A_1$ is dense in $\HH$, we can simplify the $A_1$ factor
on the right of both sides of \eqref{irrep:eq3}, and since $\ker(A_1)=\{0\}$
we can simplify the $A_1$ on the left too. We get $\pi(p_{23})=A_2V_2^*V_3A_3$.

Setting $A_1=:Z_1$, $V_2A_2=:Z_2$ and $V_3A_3=:Z_3$ concludes the proof of \eqref{eq:piP}.

Next, we  show that $\tilde{\pi}(z_i):=Z_i$ defines a $*$-representation
of $\Sq$, i.e.~that the $Z_i$'s satisfy the commutation rules of the algebra $\Sq$. Then the rest 
of the proof is straightforward: if the relations of $\Sq$ are satisfied, some algebraic manipulation  immediately gives from \eqref{eq:piP} that $Z_i^*Z_j=\pi(p_{ij})$ with $\ker Z_1=\ker\pi(p_{11})=\{0\}$
by hyphothesis. Also, since $\tilde{\pi}(\Sq) \HH\supset\pi(\Aq) \HH$, if
the latter space is dense in $\HH$, so must be the former: so if $\pi$ is
irreducible, $\tilde{\pi}$ is irreducible too.

Among the defining relations of $\CP^2_q$ we have $p_{11}p_{1i}=q^2p_{1i}p_{11}$
for $i=2,3$. Applying $\pi$ to both sides and using \eqref{eq:piP} it becomes
$Z_1(Z_1^2Z_i-q^2Z_iZ_1^2)=0$. This reduces to $Z_1^2Z_i=q^2Z_iZ_1^2$,  
since $\ker Z_1=\{0\}$ allowing to simplify the factor $Z_1$ on the left and get. Being $Z_1$ 
positive it can be diagonalized, and by previous relation $Z_i$ intertwines the eigenspace
of $Z_1^2$ with eigenvalue $\lambda\geq 0$ with the eigenspace of $Z_1^2$
with eigenvalue $q^2\lambda$, i.e.~it maps the eigenvalue $\sqrt{\lambda}$
of $Z_1$ to $q\sqrt{\lambda}$ ($-q\sqrt{\lambda}$ is excluded by the positivity of $Z_1$). This proves
\begin{equation}\label{eq:Z1Zi}
Z_1Z_i=qZ_iZ_1 \;,\qquad i=2,3\;.
\end{equation}
Since $Z_1=Z_1^*$, by conjugating we get $Z_1^*Z_i=qZ_iZ_1^*$.

The defining relation $p_{12}p_{13}=qp_{13}p_{12}$ yields $Z_1^2(Z_2Z_3-qZ_3Z_2)=0$, using \eqref{eq:piP}
and \eqref{eq:Z1Zi}. Again $Z_1^2$ can be simplified having kernel $\{0\}$ and so $Z_2Z_3=qZ_3Z_2$.
Similarly from $p_{21}p_{13}=q^3p_{13}p_{21}$ we deduce
$Z_2^*Z_3=qZ_3Z_2^*$ -- and by conjugation $Z_3^*Z_2=qZ_2Z_3^*$. 
From $Z_1=Z_1^*$ we get the relation $[Z_1^*,Z_1]=0$;
from $q^{-2}p_{21}p_{12}-p_{12}p_{21}=(1-q^2)p_{11}^2$ we get
$[Z_2^*,Z_2]=(1-q^2)Z_1Z_1^*$, from $q^{-2}p_{31}p_{13}-p_{13}p_{31}
=(1-q^2)(p_{11}^2+p_{12}p_{21})$ we get $[Z_3^*,Z_3]=(1-q^2)(Z_1Z_1^*+Z_2Z_2^*)$.
Using previous formulas for $[Z_2^*,Z_2]$ and $[Z_3^*,Z_3]$, 
the tracial relation $\tr_q\pi(P)=1$  gives the spherical relation $\sum_iZ_iZ_i^*=1$.
With this, all the defining relations of $S^5_q$ are satisfied and the proof is complete.
\end{proof}

It is in general not true that a representation on a Hilbert space of a subalgebra of a given algebra
 can be extended to a representation of the full algebra on \emph{the same}
Hilbert space:  one needs at least to extend the Hilbert space. 
The point here is that we can extend the representation
from $\Aq$ to $\Sq$ without enlarging the Hilbert space.

By previous lemma any irrep $\pi:\Aq\to\B(\HH)$ with
$\ker\pi(p_{11})=\{0\}$ is the restriction of an irrep
$\tilde{\pi}:\Sq\to\B(\HH)$ with $\ker\pi(z_1)=\{0\}$. On the other hand any
irrep $\pi:\Aq\to\B(\HH)$ with $\pi(p_{1i})=0$
for all $i=1,2,3$ is the pullback of an irrep
of the standard Podle\'s sphere, the $*$-algebra morphism $\Aq\to\A(\CP^1_q)$
being the map
$$
p_{1j},p_{j1}\mapsto 0\;\forall\;j\,,\qquad
p_{22}\mapsto A\,,\qquad
p_{23}\mapsto B^*\,,\qquad
p_{32}\mapsto B\,,\qquad
p_{33}\mapsto 1-q^2A\,,
$$
where $A$ and $B$ are the original generators of Podle\'s \cite{Pod87}.
Next lemma shows that only these two cases are possible.

\begin{lemma}\label{lemma:app2}
For any irrep $\pi$ of $\Aq$ it is either
$\ker\pi(p_{11})=\{0\}$ or $\pi(p_{1i})=0$ for all $i=1,2,3$.
\end{lemma}

\begin{proof}
Suppose there is a vector $v$ such that $\pi(p_{11})v=0$. The element
$p_{11}$ commutes with the generators $p_{ij}$ for all $i,j>1$, while
$p_{11}p_{1i}=q^2p_{1i}p_{11}$ and 
$p_{11}p_{i1}=q^{-2}p_{i1}p_{11}$ for all $i>1$.
This means that for any polynomial $a$ in the generators
$\pi(p_{11})\pi(a)v\propto \pi(a)\pi(p_{11})v=0$ and by linearity
$\pi(p_{11})\pi(a)v=0$ for any $a\in\Aq$.
We conclude that the kernel of $\pi(p_{11})$ carries a subrepresentation
of $\pi$, but $\pi$ is irreducible thus we have only two cases:
$\ker\pi(p_{11})=\{0\}$ or $\pi(p_{11})=0$. In the latter case
$$
\textstyle{\sum_i\pi(p_{i1})^*\pi(p_{i1})
=\pi\left(\sum_ip_{1i}p_{i1}\right)=\pi(p_{11})=0}\;.
$$
A sum of positive operators $\pi(p_{i1})^*\pi(p_{i1})$ is zero
iff each one is zero, and for a bounded operator
$T$ the condition $T^*T=0$ implies $T=T^*=0$. Thus $\pi(p_{11})=0$ implies that 
$\pi(p_{12})=\pi(p_{13})=0$ as well.
\end{proof}

We need to recall the representation theory of quantum spheres.
Irreps of quantum orthogonal spheres $\A(S^n_q)$
are classified in \cite{HL04} -- and in particular we are interested in $\Sq$ --,
while the case of the standard Podle\'s sphere $\CP^1_q$ is in \cite{Pod87}.
The generators $z_i$ of $\Sq$ are related to the $x_i$'s used in
\cite{HL04} by $x_i=z^*_i$ and by the replacement $q\to q^{-1}$.
We collect some results of \cite{HL04,Pod87} into the following
proposition, adapted to our notations.

\begin{prop}\label{prop:app}
Any non-trivial irrep of $\A(\CP^1_q)$ is
unitarily equivalent to one of the following two representations (cf.~\cite[Prop.~4.I]{Pod87},
with parameter $c=0$). The first is one-dimensional, $\chi_0(p_{ij})=\delta_{i3}\delta_{j3}$
(and $\chi_0(1)=1$), the second is just the map $\chi_1$ in \eqref{eq:chionep}.

Any irrep of $\Sq$ with $z_1$ not in the
kernel is unitarily equivalent to one of the following family, parametrized
by $\lambda\in U(1)$ (cf.~\cite{HL04}, eq.~(3.10), $n=2$):
\begin{align*}
\psi^5_\lambda(z_1)\ket{k_1,k_2} &:=\lambda q^{k_1+k_2}\ket{k_1,k_2} \;,\\
\psi^5_\lambda(z_2)\ket{k_1,k_2} &:=q^{k_1}\sqrt{1-q^{2(k_2+1)}}\ket{k_1,k_2+1} \;,\\
\psi^5_\lambda(z_3)\ket{k_1,k_2} &:=\sqrt{1-q^{2(k_1+1)}}\ket{k_1+1,k_2} \;.
\end{align*}
Note that $\ker\psi^5_\lambda(z_1)=\{0\}$.
For $\lambda=1$ we get the representation $\psi^5_{\lambda=1}=\chi_2$ in \eqref{eq:chitwo}.
\end{prop}

Notice that restricted to $\Aq$ all the representations $\psi^5_\lambda$ are
unitarily equivalent. Indeed, let $U$ be the unitary transformation
$U\ket{k_1,k_2}:=\lambda^{k_1+k_2}\ket{k_1,k_2}$; one easily checks that
$U\psi^5_\lambda(a)U^*=\psi^5_{\lambda=1}(a)=\chi_2(a)$ for all $a\in\Aq$
(indeed it is enough to prove it for $a=p_{ij}$).
By Lemma~\ref{lemma:app1} and Lemma~\ref{lemma:app2}
any non-trivial irrep of $\Aq$ is either the pullback of one of $\A(\CP^1_q)$ -- and then unitarily
equivalent to $\chi_0$ or $\chi_1$ -- or it is the restriction of an
irrep $\tilde{\pi}$ of $\Sq$ such that $\ker\tilde{\pi}(z_1)=\{0\}$
-- and by Prop.~\ref{prop:app} it is unitarily equivalent to (the restriction of) $\chi_2$.


\bigskip\smallskip

\begin{center}
\textsc{Acknowledgments}
\end{center}
GL was partially supported by the `Italian project Cofin06 - Noncommutative geometry, quantum groups and applications'.

\medskip


\providecommand{\bysame}{\leavevmode\hbox to3em{\hrulefill}\thinspace}

\end{document}